\documentclass[11pt, reqno]{amsart}
\usepackage{mathrsfs}
\usepackage{amssymb}
\usepackage{amsmath,amsthm,amssymb,amsfonts,dsfont,mathtools}
\usepackage{epsfig}
\usepackage{url}
\usepackage[all]{xy}

\theoremstyle{plain}
\newtheorem{thm}{Theorem}[section]
\newtheorem{prop}[thm]{Proposition}
\newtheorem{lem}[thm]{Lemma}

\newtheorem{cor}[thm]{Corollary}
\newtheorem{conj}[thm]{Conjecture}

\theoremstyle{definition}

\newtheorem{rem}[thm]{Remark}
\newtheorem{defn}[thm]{Definition}
\newtheorem{eg}[thm]{Example}
\newtheorem{subtitle}[thm]{}
\newtheorem{ex}{Exercise}[section]
\numberwithin{equation}{section}

\def\a{\alpha}
\def\b{\beta}

\def\g{\gamma}
\def\G{\Gamma}

\def\n{\,\vert\,}

\def\w{\omega}
\def\W{\Omega}
\def\li{\langle}
\def\ri{\rangle}

\def\ms{\medskip}

\def\ti{\tilde}
\def\p{\partial}

\def\sec{{\rm sec\/}}

\def\R{\mathds{R} }
\def\C{\mathds{C}}

\def\H{\mathbb{H}}
\def\N{\mathbb{N}}
\def\Z{\mathbb{Z}}

\def\comp{\ensuremath\mathop{\scalebox{.6}{$\circ$}}}

\newcommand{\beq}{\begin{equation}}
\newcommand{\eeq}{\end{equation}}
\newcommand{\beg}{\begin{eg}}
\newcommand{\eeg}{\end{eg}}
\newcommand{\bthm}{\begin{thm}}
\newcommand{\ethm}{\end{thm}}
\newcommand{\bprop}{\begin{prop}}
\newcommand{\eprop}{\end{prop}}
\newcommand{\bcor}{\begin{cor}}
\newcommand{\ecor}{\end{cor}}
\newcommand{\blem}{\begin{lem}}
\newcommand{\elem}{\end{lem}}
\newcommand{\bca}{\begin{cases}}
\newcommand{\eca}{\end{cases}}
\newcommand{\brem}{\begin{rem}}
\newcommand{\erem}{\end{rem}}
\newcommand{\bconj}{\begin{conj}}
\newcommand{\econj}{\end{conj}}
\newcommand{\bpm}{\begin{pmatrix}}
\newcommand{\epm}{\end{pmatrix}}
\newcommand{\bbm}{\begin{bmatrix}}
\newcommand{\ebm}{\end{bmatrix}}
\newcommand{\bvm}{\begin{vmatrix}}
\newcommand{\evm}{\end{vmatrix}}
\newcommand{\bdefn}{\begin{defn}}
\newcommand{\edefn}{\end{defn}}
\newcommand{\bsub}{\begin{subtitle}}
\newcommand{\esub}{\end{subtitle}}
\newcommand{\bex}{\begin{ex}}
\newcommand{\eex}{\end{ex}}
\newcommand{\ben}{\begin{enumerate}}
\newcommand{\een}{\end{enumerate}}

\begin{document}

\title{On non-elliptic symplectic manifolds$^\ast$}
\thanks{$^\ast$The work is supported by PRC grant NSFC 11401514(Fang), 11771377,12171417 (Wang), and the Natural Science Foundation of Jiangsu Province
(BK20191435)(Fang).}

\author{Shouwen Fang}
\address{School of Mathematical Science\\
Yangzhou University\\ Yangzhou, Jiangsu 225002, China }
\email{shwfang@163.com}

\author{Hongyu Wang$^\dagger$}
\thanks{$^\dagger$E-mail: hywang@yzu.edu.cn.}
\address{School of Mathematical Sciences\\
Yangzhou University\\ Yangzhou, Jiangsu 225002, China}
\email{hywang@yzu.edu.cn}

\dedicatory{In Memory of Professor Shiing-Shen Chern (1911-2004)}
\ms
\hskip 3in 

\begin{abstract}
Let $M$ be a closed symplectic manifold of dimension $2n$ with non-ellipticity. We can define an almost K\"ahler structure on $M$ by using the given symplectic form.
Hence, we have a $\G=\pi_1(M)$-invariant almost K\"ahler structure on the universal covering, $\ti M$, of $M$.
Using Darboux coordinate charts, we globally deform the given almost K\"ahler structure on $\ti M$ off a Lebesgue measure zero subset to obtain a $\G$-invariant (measurable) Lipschitz K\"ahler flat structure
on $\ti M$ which is a singular K\"ahler structure and $\G$-homotopy equivalent to the given almost K\"ahler structure with Lipschitz condition, restricted to a $\G$-invariant open dense submanifold of $\ti M$, the K\"ahler flat metric is quasi isometric to the given almost K\"ahler metric. Analogous to Teleman's $L^2$-Hodge decomposition on PL manifolds or Lipschitz Riemannian manifolds, we give a $L^2$-Hodge
decomposition theorem on $\ti M$ with respect to the (measurable) Lipschitz K\"ahler flat metric. As done in K\"ahler case, using an argument of Gromov, we give a vanishing theorem for $L^2$ harmonic $p$-forms, $p\not=n$
(resp. a non-vanishing theorem for $L^2$ harmonic $n$-forms) on $\ti M$, then the signed Euler characteristic satisfies $(-1)^n\chi(M)\geq0$ (resp. $(-1)^n\chi(M)>0$). Similarly, for any closed even dimensional Riemannian manifold $(M, g)$, we can construct a $\G$-invariant (measurable) Lipschitz
K\"ahler flat structure on the universal covering, $(\ti M, \ti g)$, of $(M, g)$ which is a singular K\"ahler structure and $\G$-homotopy equivalent to and quasi-isometric to the metric $\ti g$.
As an application, as done in smooth case, using Gromov's method we show that the Chern-Hopf conjecture holds true in closed even dimensional Riemannian manifolds with nonpositive curvature (resp. strictly negative curvature), it gives a positive answer to a Yau's
problem due to S. S. Chern and H. Hopf.
\end{abstract}
\maketitle
\section{Introduction}
In K\"ahler geometry, M. Gromov \cite{Gr2} considered K\"ahler hyperbolicity; N. Hitchin \cite{Hit}, J. Jost and K. Zuo \cite{JZ}, J. Cao and F. Xavier \cite{CX}, B. L. Chen and X. K. Yang \cite{CY} considered K\"ahler parabolicity. Similarly, Q. Tan, H. Y. Wang and J. R. Zhou \cite{TWZ2}, T. Huang and Q. Tan \cite{HuangT} considered symplectic parabolicity whose underlying manifolds satisfy the hard Lefshetz condition, and J. Jost and Y. L. Xin considered pinched curvature case \cite{JX}.

Let $(M, \w)$ be a closed symplectic manifold. Let $J$ be an $\w$-compatible almost complex structure, i.e., $J^2=-\text{id}$, $\w(J\cdot, J\cdot)=\w(\cdot,\cdot)$, and $g(\cdot,\cdot)=\w(\cdot, J\cdot)$ is a $J$-compatible (that is, $g(J\cdot, J\cdot)=g(\cdot, \cdot)$) Riemannian metric on $M$. The triple $(\w,J,g)$ is called an almost K\"ahler structure on $M$. Notice that any one of the pairs $(\w,J)$, $(J,g)$ or $(g,\w)$ determines the other two. An almost-K\"ahler metric $(\w,J,g)$ is K\"ahler if and only if $J$ is integrable. We have the following definition (cf. M. Gromov \cite{Gr1,Gr2}, N. Hitchin \cite{Hit}, Jost-Zuo \cite{JZ}, Cao-Xavier \cite{CX} or Tan-Wang-Zhou \cite{TWZ2}):

\bdefn (1) A closed almost K\"ahler manifold $(M,\w,J,g)$ is called {\it symplectic hyperbolic\/} if the lifting $\ti\w$ of $\w$ to the universal covering $(\ti M,\ti\w,$ $\ti J,\ti g)$ is $d$(bounded), that is, there is a 1-form $\b$ on $\ti M$ such that $\ti\w=d\b$, where $$\|\b\|_{L^\infty} =\sup_{x\in \ti M} |\b(x)|_{\ti g} <+\infty.$$

\ms
(2)  A closed almost K\"ahler manifold $(M,\w,J,g)$ is called {\it symplectic parabolic\/} if the lifting $\ti\w$ of $\w$ to the universal covering $(\ti M,\ti\w, \ti J,\ti g)$ is $d$ (sublinear), that is, there is a 1-form $\b$ on $\ti M$ such that $\ti\w=d\b$, where $$|\b (x)|_{\ti g}  \leq C(\ti \rho(x_0,x)+1),\quad C>0$$ is a constant, $\ti \rho(x_0,x)$ is the distance function on $\ti M$ relative to a base point $x_0\in\ti M$, with respect to the metric $\ti g$.
\edefn

\beg
\begin{itemize}
\item If $M^{2n}$ is homotopy equivalent to a closed Riemannian manifold with negative sectional curvature and having convex boundary (if any), then $M^{2n}$ is K\"ahler (symplectic) hyperbolic provided it admits some K\"ahler (symplectic) structure (see M. Gromov \cite{Gr2}).
\item If $M^{2n}$ is a closed K\"ahler (symplectic) manifold such that $\pi_1(M^{2n})$ is word-hyperbolic in the sense of M. Gromov \cite{Gr1} and $\pi_2(M^{2n})=0$, then $M^{2n}$  is K\"ahler (symplectic) hyperbolic.
\item Let $M^{2n}$ be a closed Riemannian manifold with Anosov geodesic flow. If $M^{2n}$ is homotopy equivalent to a closed  K\"ahler (symplectic) manifold, then $M^{2n}$ is K\"ahler (symplectic) hyperbolic (see X. Cheng \cite{Chg}).
\item Let $M^{2n}$ be a closed Riemannian manifold of non-positive curvature. If $M^{2n}$ is homeomorphic to a K\"ahler (symplectic) manifold, then  $M^{2n}$ is K\"ahler (symplectic) parabolic (see Cao-Xavier \cite{CX}, Jost-Zuo \cite{JZ}).
\item The nilmanifolds in\cite[Chapter 2]{TO}, are aspherical and symplectic parabolic with no hard Lefschetz property.
\item A closed K\"ahler manifold $(M^{2n}, \w)$ with $\pi_1(M)$ being CAT(0) or automatic, then $M^{2n}$ is K\"ahler parabolic (see Chen-Yang \cite{CY}).
\end{itemize}
\eeg

It is clear that for any non-elliptic symplectic manifold $M$, (that is, symplectic hyperbolic or parabolic), the fundamental group of $M$, $\pi_1(M)$, is infinite.

The space of spherical classes \cite{GM}, $\Pi(M)\subset H_2 (M; \Z)$, is defined by the image of the Hurewicz homomorphism $$H: \pi_2(M)\to H_2(M; \Z).$$ Note that if $f: S^2 \to M$, then there is induced $$f_*: H_2 (S^2) \to H_2(M;\Z).$$ Thus we have $f_*([S^2]) \in H_2(M;\Z)$, where $[S^2]$ is the fundamental class of $S^2$. This determines a natural transformation $$H: \pi_2(M)\to H_2(M;\Z),$$ given by $[f] \mapsto f_*([S^2])$.

\bdefn (R. Gompf \cite{Gom})
A closed symplectic manifold $(M,\w)$ is called {\it symplectically aspherical\/} if $\w$ vanishes on all spherical classes.
\edefn

Recall the definition of aspherical manifolds (cf. \cite{Dav,Luck}):

\bdefn
A manifold is called {\it aspherical\/} if it is connected, and its universal covering is contractible.
\edefn

\brem (1) R. Gompf \cite{Gom} pointed out that symplectically aspherical manifolds may not be aspherical.

(2) M. W. Davis \cite{Dav} used the fundamental group at infinity to get the following results:

In each dimension $n\geq 4$, there are closed, aspherical manifolds $M^n$ with universal covering $\ti M^n$ not homeomorphic to the Euclidean space $\R^n$.
\erem

\bconj{\rm(}Chern-Hopf conjecture, cf. Chern \cite[p.123]{Chern}, \cite[Problem III]{Chern2} or \cite{Dav,Luck,Yau}{\rm)}
If $M^{2n}$ is an aspherical closed manifold of dimension $2n$, then
$$ (-1)^n \chi (M^{2n}) \geq 0.$$

\ms
If $M^{2n}$ is a closed Riemannian manifold of dimension $2n$ with sectional curvature $\sec (M^{2n})$, then
$$ (-1)^n \chi (M^{2n}) >0 \ ({\rm or\/}\ \chi (M^{2n}) >0), \ \text{if }\ \sec (M^{2n}) <0\ ({\rm or\/}\ >0);
$$
$$ (-1)^n \chi (M^{2n}) \geq 0 \ ({\rm or\/}\ \chi (M^{2n})\geq0), \ \text{if }\ \sec (M^{2n}) \leq 0\ ({\rm or\/}\ \geq 0);
$$
$$  \chi (M^{2n}) =0, \ \text{if }\ \sec (M^{2n}) \equiv 0.
$$
\econj

In original version of the Chern-Hopf conjecture the statement for aspherical manifolds did not appear. Notice that any Riemannian manifold with nonpositive sectional curvature is aspherical by Cartan-Hadamard theorem
(cf. \cite{Au,Cha}). The Chern-Hopf conjecture first appeared in print in Chern's paper \cite[p.123]{Chern} on the general Gauss-Bonnet Theorem. This question were asked by H. Hopf in 1932 (cf. M. Berger \cite[pp.544-545]{Ber} or M. W. Davis \cite[p.314]{Dav} or
S. T. Yau \cite[Problem 10]{Yau}). Underlying Hopf's question is the question of whether the conjecture follows directly from the Gauss-Bonnet Theorem, that is, does the nonpositivity of sectional curvature imply that the sign
of the Gauss-Bonnet integrand is $(-1)^n$? I. M. Singer \cite{Sing} had proposed to settle this problem by looking at the universal cover of $M$. He pointed out that if the $L^2$ harmonic forms on the universal cover are all zero except in the middle dimension, then one can apply the $\G$-index theorem for coverings (see M. Atiyah \cite{At}) to prove the statement in the affirmative.

The Chern-Hopf conjecture holds true in two dimensions by the Gauss-Bonnet formula immediately. In four dimensions, one can still check that positive (resp. negative) sectional curvature implies that the Gauss-Bonnet integrand is pointwise positive (cf. J. Milnor \cite[\S4]{Mil} or S. S. Chern \cite[\S4]{Chern}). However, in higher dimensions, it is known that the sign of the sectional curvature does not determine the sign of the Gauss-Bonnet integrand, see R. Geroch \cite{Ger}. M. Gromov \cite{Gr2} gave a positive answer to the Chern-Hopf conjecture when the metric is K\"ahler since K\"ahler metric has the hard Lefschetz property. R. Charney and M. W. Davis \cite{CD} investigated the conjecture in the context of piecewise Euclidean manifolds having ``nonpositive curvature" in the sense of Gromov's CAT(0) inequality (cf. \cite{CY,DK,Gr1}). Tan-Wang-Zhou \cite{TWZ2}, Huang-Tan \cite{HuangT} gave a positive answer to Chern-Hopf conjecture when the closed symplectic manifold $(M, \w)$ is symplectic parabolic with the hard Lefschetz condition.
\vskip 0.1cm
In this paper, we study non-elliptic symplectic manifolds and even dimensional closed non-positive manifolds. Supposed that $(M, \w, J, g_J)$ is a closed non-elliptic almost K\"ahler manifold. Then its universal cover $(\ti M, \ti\w, \ti J, g_{\ti J})$ is a $\G=\pi_1$-equivariant simply connected almost K\"ahler manifold. Using Darboux coordinate charts, we globally deform the given almost K\"ahler structure $(\ti\w, \ti J, g_{\ti J})$ on $\ti M$ off a Lebesgue measure zero subset to get a family of $\G$-invariant (measurable) Lipschitz almost K\"ahler structures $(\ti\w, \ti J(t), g_{\ti J}(t))$ on $\ti M$, $0\leq t\leq 1$ (cf. Drutu-Kapovich \cite{DK}), where $(\ti\w, \ti J(0), g_{\ti J}(0))$ is a $\G$-invariant (measurable) Lipschitz K\"ahler structure, and $(\ti\w, \ti J(1),$ $g_{\ti J}(1))$$=(\ti\w, \ti J, g_{\ti J})$. We will prove that with those given $\G$-invariant (measurable) Lipschitz almost K\"ahler metrics $g_{\ti J}(t)$,$0\leq t\leq 1$, the corresponding $L^2_k$ norms of $p$-form on $\ti M^{2n}$ are quasi isometry, that is, with equivalent norms (cf. Gromov \cite{Gr2}). As done in K\"ahler case, we will prove vanishing and non-vanishing theorems for $\G$-invariant (measurable) Lipschitz K\"ahler flat structure $(\ti\w, \ti J(0), g_{\ti J}(0))$ which is a singular K\"ahler structure (hence, for $\G$-invariant almost K\"ahler structure $(\ti\w, \ti J, g_{\ti J})$) since $[L_{\ti\w}, \Delta_{g_{\ti J}(0)}]=0$ which implies that $g_{\ti J}(0)$ has the hard Lefschetz property, where
$\Delta_{g_{\ti J}(0)}=d^{*_0}d+dd^{*_0}$, $d^{*_0}=-*_0d*_0$, and $*_0$ is the Hodge star operator with respect to the (measurable) Lipschitz K\"ahler flat metric $g_{\ti J}(0)$.
Thus, we have the following similar result (cf. M. Gromov \cite[main theorem]{Gr2}):

\bthm \label{main}
(1) If $(M^{2n},\w)$ is a $2n$-dimensional closed symplectic hyperbolic manifold, then
$$  (-1)^n \chi (M^{2n}) >0.$$

(2) If $(M^{2n},\w)$ is a $2n$-dimensional closed symplectic parabolic manifold, then
$$  (-1)^n \chi (M^{2n}) \geq 0.$$
\ethm

\brem (1) The theorem above drops the condition that the closed symplectic manifold $(M^{2n}, \w)$ has the hard Lefschetz property (see Tan-Wang-Zhou \cite[Theorem 1.5]{TWZ2} or Huang-Tan \cite[Theorem 1.5]{HuangT}), hence it gives a positive answer to a question
in \cite[Question 1.7]{TWZ2}.

(2) It is easy to see that every closed non-elliptic symplectic manifold is symplectically aspherical. It is natural to ask whether symplectically aspherical manifolds are always non-elliptic symplectic.
\erem
 Suppose that $(M, g)$ is a closed Riemannian manifold of even dimension with non positive curvature. Then the universal covering $$\pi:(\ti M,\ti g)\to(M, g)$$ is a complete simply connected manifold with deck transformation group $\G=\pi_1(M)$. By Cartan-Hadamard theorem, there is a diffeomorphism $\Phi:\ti M\to\R^{\dim M}\cong\C^{\frac{1}{2}\dim M}$.  Hence, as done in symplectic case, we can construct a family of $\G$-invariant (measurable) Lipschitz almost K\"ahler  structures $(\ti\w', \ti J'(t), \ti g'(t))$, $0\leq t\leq 1$, on $\ti M$, where $\ti g'(1)=\ti g$ is a smooth $\G$-invariant Riemannian metric on $\ti M$, and $\ti g'(0)$ is a $\G$-invariant (measurable) Lipschitz K\"ahler flat metric on $\ti M$ which is a singular K\"ahler metric on $\ti M$. Notice that $\ti\w'$ is a $\G$-invariant symplectic form on $\ti M$ off a Lebesgue measure zero subset (in fact, $\ti\w'$ is a symplectic form on an open dense submanifold $\overset{\circ}{\ti M}$, $\ti M\setminus \overset{\circ}{\ti M}$ is a Lebesgue measure zero set), $\ti g'(0)$ is quasi isometric to $\ti g'(1)$ on $\overset{\circ}{\ti M}$, and
 $$
\frac{\ti\w'^n}{n!}=dvol_{\ti g'(t)}=dvol_{\ti g},\quad 0\leq t\leq 1,
 $$
 is smooth $\G$-invariant volume density on $\ti M$. Note that $\ti\w'$ is defined on the open dense submanifold $\overset{\circ}{\ti M}$, in general, $\ti\w'$ can not be extended to $\ti M$. Using the same method of proving Theorem \ref{main}, we obtain the following theorem:
\bthm\label{mc}
(1) If $(M^{2n}, g)$ is a $2n$-dimensional closed Riemannian manifold with strictly negative sectional curvature, then
$$  (-1)^n \chi (M^{2n}) >0.$$

(2) If $(M^{2n}, g)$ is a $2n$-dimensional closed Riemannian manifold with non-positive sectional curvature, then
$$  (-1)^n \chi (M^{2n}) \geq 0.$$
\ethm
\brem
Theorem \ref{mc} shows that the Chern-Hopf conjecture holds true in closed Riemannian manifolds with nonpositive curvature (resp. strictly negative
curvature). In particular, it gives a positive answer to a Yau's problem \cite[Problem 10]{Yau} due to S. S. Chern and H. Hopf:

Let $M$ be a $2n$-dimensional closed oriented Riemannian manifold of strictly negative sectional curvature $K<0$. Then the signed Euler characteristic satisfies $ (-1)^n \chi (M^{2n}) >0.$
\erem
\vskip 0.1cm
The remainder of the paper is organized as follows:

In Section 2, we devote to constructing (measurable) Lipschitz K\"ahler flat structures which can be regarded as singular K\"ahler structures on closed symplectic manifolds.

In Section 3, as done in smooth case, we will introduce $L^2$-Hodge theory on the universal covering of $M$ with respect to the (measurable) Lipschitz K\"ahler flat metric.

In Section 4, similar to K\"ahler case, we prove the vanishing theorem of $L^2$-Betti numbers of symplectic parabolic manifolds.

In Section 5, we will investigate the non-vanishing theorem for $L^2$-cohomology on symplectic hyperbolic manifolds, and give a proof of Theorem \ref{main}.

Finally, in Section 6, by using the technique developed in Section 2--5, we will prove Theorem \ref{mc}.
\section{Measurable K\"ahler flat structures on closed symplectic manifolds}
This section is devoted to studying measurable K\"ahler flat structures on closed symplectic manifolds.

We first describe closed Lipschitz manifolds and some basic properties. Lipschitz structure on a topological manifold $M^n$ of dimensional $n$ is a maximal atlas $\mathscr{U}=\{U_\alpha,\phi_\alpha\}_{\alpha\in\Lambda}$ where $$\phi_\alpha:U_\alpha\to V_\alpha=\phi_\alpha(U_\alpha)\subseteq\R^n$$
is a homeomorphism from an open set $U_\alpha\subset M^n$ onto an open subset of $\R^n$, and the changes of coordinates $\phi_\beta\comp\phi_\alpha^{-1}$ are Lipschitz functions, that is,
$$|\phi_\beta\comp\phi_\alpha^{-1}(x)-\phi_\beta\comp\phi_\alpha^{-1}(y)|\leq c_{\alpha\beta}|x-y|$$ for any $x, y\in\phi_\alpha(U_\alpha\cap U_\beta)$ with $c_{\alpha\beta}$ a constant (cf. N. Teleman \cite{Tele3} or J. Luukkainen and J. V\"ais\"al\"a \cite{LV}).

Let $M^n$ be a Lipschitz manifold as above. A $L^2$-form $\sigma$ of degree $p$ on $M^n$ is, by definition, a system $\sigma=\{\sigma_\alpha\}_{\alpha\in\Lambda}$ where each $\sigma_\alpha$ is a complex $L^2$-differential form
of degree $p$ on the open subset $V_\alpha$ and they are required to satisfy the compatibility conditions
\beq \label{a1}
(\phi_\b\comp\phi_\a^{-1})^*\sigma_\b=\sigma_\a,
\eeq
the pull-back $(\phi_\b\comp\phi_\a^{-1})^*$ is defined component-wise (cf. N. Teleman \cite{Tele3} or J. Luukkainen and J. V\"ais\"al\"a \cite{LV}). This definition makes sense in view of the following result (see e.g. H. Whitney \cite[p. 272]{Wh}) or J. Cheeger \cite{Cheeger}):
\bprop[Rademacher] \label{} Let $f:U\to \R$ be a Lipschitz function defined on an open subset $U$ of $\R^n$. Then \\
(1) the partial derivatives $\frac{\partial}{\partial x_k}f$ exist almost everywhere on $U$, $1\leq k\leq n$,\\
(2) $\frac{\partial}{\partial x_k}f$ are measurable and $L^{\infty}$-bounded.
\eprop
It is well-known that on any Lipschitz manifold, $L^2$-forms, currents and exterior derivatives can be defined (see D. Sullivan \cite{Su1,Sul}). The space of $L^2$-forms of degree $p$ on $M^n$ will be denoted by $L^2\Omega^p(M^n)$. If $\sigma$ is any $L^2$-form of degree $p$ on an open subset $U$ of $\R^n$, then is said (classically) to have distributional exterior derivative $d\sigma$ in $L^2$ if there exists an $L^2$-form denoted $d\sigma$ in
$L^2\Omega^{p+1}(U)$ such that for any $C^\infty$ $(n-p-1)$-form $\varphi$ with compact support in $U$:
$$\int_U\sigma\wedge d\varphi=(-1)^{p+1}\int_Ud\sigma\wedge\varphi.$$
If $\sigma=\{\sigma_\a\}_{\a\in\Lambda}\in L^2\Omega^p(M^n)$, and if $d\sigma_\a\in L^2\Omega^{p+1}(\R^n)$ for any $\a\in\Lambda$, then we say that $\sigma$ has distributional exterior derivative $d\sigma=\{d\sigma_\a\}_{\a\in\Lambda}$ in $L^2\Omega^{p+1}(M^n)$. Of course, in order to check that this definition is correct, it remains to verify that the forms $d\sigma_\a$ satisfy the compatibility condition:
\beq \label{a2}
(\phi_\b\comp\phi_\a^{-1})^*d\sigma_\b=d\sigma_\a.
\eeq
The relation \ref{a2} follows from the following proposition:
\bprop{\rm(}cf. N. Teleman \cite[Proposition 1.2]{Tele3} {\rm)} \label{} For any Lipschitz mapping $f:U_1\to U_2$, where $U_1$ and $U_2$ are relatively compact open sets in $\R^n$, and for any form $\sigma\in L^2\Omega^p(U_2)$, the
form $f^*\sigma$ belongs to $L^2\Omega^p(U_1)$, and $d(f^*\sigma)=f^*(d\sigma)\in L^2\Omega^{p+1}(U_1)$.
\eprop
We introduce the spaces
\beq \label{a3}L^2_d\Omega^p(M^n)=\{\sigma|\sigma\in L^2\Omega^p(M^n), d\sigma\in L^2\Omega^{p+1}(M^n)\}.\eeq
The exterior derivative $d$ satisfies $d^2=0$, therefore
\begin{align*} \label{}
L^2_d\Omega^*(M^n)\coloneqq\{0&\to L^2_d\Omega^0( M^n)\xrightarrow{d}L^2_d\Omega^1( M^n)\xrightarrow{d}\cdots\nonumber\\
&\xrightarrow{d}L^2_d\Omega^{p}( M^n)\xrightarrow{d}L^2_d\Omega^{p+1}( M^n)\xrightarrow{d}\cdots\}
\end{align*}
is a cohomology complex.

A Riemannian metric on a Lipschitz manifold $M^n$ is collection $g=\{g_\a\}_{\a\in\Lambda}$ (cf. Drutu-Kapovich \cite{DK} or N. Teleman \cite{Tele3} or J. Luukkainen and J. V\"ais\"al\"a \cite{LV}), where $g_\a$ is a Riemannian metric on $V_\a=\phi_\a(U_\a)\subseteq\R^n$,
with measurable components, which satisfy the compatibility conditions:
\beq \label{a4}
(\phi_\b\comp\phi_\a^{-1})^*g_\b=g_\a.
\eeq
In addition, the Riemannian metrics $g_\a$ are required to define $L^2$-norm on $V_\a$ which are equivalent to the standard $L^2$-form, that is, there should exist a positive constant $c_\a>1$ such that, for any smooth form $\sigma$ with compact support in $V_\a$,
\beq \label{a5}
c_\a^{-1}||\sigma||\leq||\sigma||_\a\leq c_\a||\sigma||,
\eeq
here $||\cdot||$ and $||\cdot||_\a$ denote the usual $L^2$-forms:
\begin{align}\label{a6}
||\sigma||^2&=\int_{V_\a}\sigma\wedge\overline{*_{\R^n}\sigma},\nonumber\\
||\sigma||_\a&=\int_{V_\a}\sigma\wedge\overline{*_{\a}\sigma},
\end{align}
where $*_{\R^n}$ and $*_{\a}$ are the Hodge star operators of the Euclidean metric and the metric $g_\a$ respectively (cf. I. Chavel \cite{Cha}).

Any such Riemannian metric $g$ will be called a Lipschitz Riemannian metric on $M^n$ (cf. N. Teleman \cite{Tele3} or J. Luukkainen and J. V\"ais\"al\"a \cite{LV}).

We show now that any closed Lipschitz manifold $M^n$ have Lipschitz Riemannian metrics. For a closed Lipschitz manifold $M^n$, we choose a finite Lipschitz  subatlas $\{U_i, \phi_i\}_{1\leq i\leq N}$ on $M^n$ such that all $U_i$ have boundaries of Lebesgue measure zero (Hausdorff dimension = $n-1$). We partition $M^n$ by the sets:
\begin{align*}
&T_1=U_1,\\
&T_i=U_i  \setminus \bigcup\limits_{r=1}^{i-1}U_r,\quad 2\leq i\leq N,
\end{align*}
and we transfer on each $T_i$ the standard Euclidean metric of $\R^n$ via the coordinate map $\phi_i$. Then $\{T_i\}_{1\leq i\leq N}$ is a disjoint partition of $M^n$,
$$M^n=\bigcup\limits_{i=1}^N T_i.$$
Let
$$W_1=U_1,\quad W_i=U_i  \setminus \bigcup_{r=1}^{i-1}\overline{U}_r,\quad 2\leq i\leq N.$$
Then
$$ \overset{\circ}{ M}{}^{n}=\bigcup_{i=1}^N W_i$$
is an open dense submanifold of $M^n$, and $M^n\setminus\overset{\circ}{ M}{}^{n}$ has Hausdorff dimension $\leq n-1$ (Lebesgue measure zero).

For any Lipschitz Riemannian metric $g=\{g_\a\}_{\a\in\Lambda}$, and for any $$\sigma=\{\sigma_\a\}_{\a\in\Lambda}\in L^2\Omega^p(M^n),\quad 0\leq p\leq n,$$ the form $*_g\sigma$ defined by
\beq \label{a7}
*_g\sigma=\{*_{g_\a}\sigma_\a\}_{\a\in\Lambda}
\eeq
is an $L^2$-form of complementary degree on $M^n$. We define by means of it the $L^2$-form of $\sigma$:
\beq \label{a8}
||\sigma||^2_g=\int_{M^n}\sigma \wedge \overline{*_g\sigma}.
\eeq

This is a norm on $L^2\Omega^p(M^n)$, $0\leq p\leq n$, which makes $L^2\Omega^p(M^n)$ a Hilbert space. It is not hard to see that two different Lipschitz Riemannian metrics on $M^n$ define equivalent norms on $L^2\Omega^p(M^n)$
(cf. N. Teleman \cite{Tele3} or J. Luukkainen and J. V\"ais\"al\"a \cite{LV}).

As done in smooth case (cf. I. Chavel \cite{Cha}), we are able to consider Hodge theory on Lipschitz manifolds (cf. N. Teleman \cite[\S4]{Tele3}).

Let $g=\{g_\a\}_{\a\in\Lambda}$ be a Lipschitz metric on manifold $M^n$, and let $*_g$ be Hodge star operator. This operator is an isometry $$*_g:L^2\Omega^p(M^n)\to L^2\Omega^{n-p}(M^n).$$ The operator $d^*$ acting on $L^2\Omega^p(M^n)$ is introduced as in the smooth case $d^*=(-1)^{n(p+1)+1}*_gd*_g$. Its domain of $d^*$ is the space
\beq\label{a9}
L^2_{d^*}\Omega^p(M^n)=*_gL^2_d\Omega^{n-p}(M^n)=L^2\Omega^p(M^n).
\eeq

The following spaces $L^2_1\Omega^p(M^n)$ are the Lipschitz analogues of the Sobolev spaces of forms of degree $p$ and order one:
\beq \label{a10}
L^2_1\Omega^p(M^n)=\{\sigma\in L^2\Omega^p(M^n)|d\sigma\in L^2\Omega^{p+1}(M^n),d*_g\sigma\in L^2\Omega^{n-p+1}(M^n)\}.
\eeq

It is easy to check that $L^2_1\Omega^p(M^n)$ is a Hilbert space under the diagonal norm $||\cdot||_1$:
\beq \label{a11}
||\sigma||^2_1=||\sigma||^2+||d\sigma||^2+||d^*\sigma||^2.
\eeq

We define now the space of $(d,d^*)$-harmonic
\beq \label{a12}
\mathcal{H}^p(M^n)=\{\sigma\in L^2_1\Omega^p(M^n)| d\sigma=0=d^*\sigma\}.
\eeq

We also define $L^2$ de Rham cohomologies of $M^n$:
\beq \label{a13}
H^p(M^n, \C)=\ker d|_{L^2_d\Omega^p(M^n)}/dL^2_d\Omega^{p-1}(M^n).
\eeq

Any harmonic form is a cocycle in $L^2$ de Rham cohomologies $H^*(M^n, \C)$:
\begin{align}\label{a14}
\chi^p:\mathcal{H}^p(M^n)&\to H^p(M^n, \C),\nonumber\\
\a\quad&\mapsto \quad [\a]
\end{align}
is called Hodge homomorphism.
\bthm {\rm(}cf. N. Teleman \cite[Theorem 4.1]{Tele3}{\rm)}\label{thm}
For any closed oriented Lipschitz Riemannian manifold $(M^n, g)$ and for any degree $p$:\\
(1) the Hodge homomorphism $$\chi^p:\mathcal{H}^p(M^n)\to H^p(M^n, \C)$$
is an isomorphism;\\
(2) $*_g:\mathcal{H}^p(M^n, g)\to \mathcal{H}^{n-p}(M^n, g)$ is an isomorphism;\\
(3) there is a strong Hodge decomposition
$$L^2 \W^p (M^n) =\mathcal{H}^p(M^n, g)\oplus dL^2_d\Omega^{p-1}(M^n)\oplus d^*L^2_d\Omega^{p+1}(M^n).$$
\ethm
\brem
The key lemma of the proof of Theorem \ref{thm} can be proved by a little convolution argument.
\erem
For a closed PL manifold $(M^n, g)$, since there exists Stokes' formula on $(M^n, g)$ (cf. N. Teleman \cite[Proposition 3.2]{Tele1}, for PL-forms over a cell complex (see D. Sullivan \cite{Su1}) PL distributional adjoint operator of $d$ is $d^*=\pm*_gd*_g$ (cf. N. Teleman \cite[Chapter II]{Tele2}) by using the method of harmonic analysis (cf. E. M. Stein \cite{St}).

For a closed Lipschitz manifold $(M^n, g)$, a theory of signature operators
(which is Dirac operator $d+d^*$ when $n=4m$) has been developed by N. Teleman \cite[Lipschitz Hodge theory]{Tele3} and \cite{Tele4}, $d^*$ is formally defined as $d^*=\pm*_gd*_g$, and it has sufficient good properties to mimic the usual linear elliptic analysis. Here $*_g$ is the Hodge star operator with respect to the metric $g$. As done in smooth case, for Lipschitz manifolds we have the Rellich type lemma (cf. \cite[Theorem 7.1]{Tele3})
by using basic properties of classical elliptic pseudodifferential operators \cite[Chapter I]{Shu2}. Notice that D. Sullivan \cite{Sul} proved the following result:

``If $n\neq4$, any topological $n$ manifold admits a Lipschitz structure. Also, any two Lipschitz structures are equivalent by a homeomorphism isotopic to the identity."

It is worth to remarking that we have the diagram
\begin{displaymath}\xymatrix@R=1ex{
 DIFF\ar[dr] & &\\
&LIP\ar[r]&TOP\\
PL\ar[ur] &&\\}
\end{displaymath}
where DIFF means differential structure, PL piecewise linear structure, LIP locally Lipschitz structure, and TOP topological structure (cf. J. Luukkainen and J. V\"ais\"al\"a \cite{LV}). Thus, for any closed smooth manifold there exist Lipschitz structures which are with equivalent norms (cf. N. Teleman \cite{Tele3} or J. Luukkainen and J. V\"ais\"al\"a \cite{LV}).

\vskip 0.1cm
In the remainder of this section, we devote to constructing (measurable) Lipschitz almost K\"ahler structures on closed almost K\"ahler manifolds which are with equivalent norms, that is, with a Lipschitz condition. But
these structures are not, in general, Lipschitz structures. Our method is similar to the construction of Lipschitz Riemannian metrics on closed Lipschitz manifolds.

Suppose that $(M^{2n}, \w)$ is a closed symplectic manifold of dimension $2n$, this implies that $M^{2n}$ is a closed Lipschitz manifold. Let $J$ be an $\w$-compatible almost complex structure on $M^{2n}$, and $g_J(\cdot,\cdot)=\w(\cdot,J\cdot)$. Then $(\w, J, g_J)$ is an almost K\"ahler structure on $M^{2n}$. We choose a Darboux's coordinate atlas $\{U'_\a, \varphi_\a\}_{\a\in\Lambda}$ such that
$$\varphi_\a:U'_\a\to\varphi_\a(U'_\a)\subseteq\R^{2n}\cong\C^n$$
is a diffeomorphism and $\varphi^*_\a\w_0=\w|_{U'_\a}$, where
$$\w_0=\frac{\sqrt{-1}}{2}\sum_{i=1}^ndz_i\wedge d\bar{z}_i$$
is the standard K\"ahler form on $\C^n$ (more details, see McDuff-Salamon \cite{MS}). We may assume that, without loss of generality, $\{o\}\in\varphi_\a(U'_\a)$ is contractible and relatively compact, strictly pseudoconvex domain (cf. L. H\"omander \cite{Hom}).

Let $\mathcal{A}(\w)$ be the set of all associated metrics with respect to the given symplectic form on $M^{2n}$, $\mathcal{H}$ the set of all metrics with the same volume element $\frac{\w^n}{n!}$ which is totally geodesic in the set $\mathcal{M}$ of all Riemannian metrics on $M^{2n}$. Note that on a compact $M^{2n}$, the set $\mathcal{M}$ of all Riemannian metrics may be given a Riemannian metric (cf. D. Ebin \cite{Ebin}). For symmetric tensor fields
$S,T$ of type $(2,0)$
$$\langle S,T\rangle_g=\int_{M^{2n}}S_{ij}T_{kl}g^{ik}g^{jl}dvol_g.$$
The geodesics in $\mathcal{M}$ are the curves $ge^{S_t}=g_t$ is computed by $g_t(x,y)=g(x,e^{S_t}y)$, where $e^{S_t}$ acts on $y$ as a tensor field of type $(1,1)$. The $\mathcal{A}(\w)$  is totally geodesics in $\mathcal{H}$ in the sense that if $h$ is $J$-anti-invariant symmetric $(2,0)$-tensor field, where $J$ is an $\w$-compatible almost complex structure on $M^{2n}$ and $g(\cdot,\cdot)=\w(\cdot,J\cdot)$, then $g_t=ge^{ht}$ lies in $\mathcal{A}(\w)$. In particular, two metrics in $\mathcal{A}(\w)$ may jointed by a unique geodesic (cf. D. E. Blair \cite[Propositions in p. 304, p. 307]{Blair}).

On the Darboux's coordinate charts, one can make the local deformations of $\w$-compatible almost complex structures (cf. D. E. Blair \cite{Blair}, also see Tan-Wang-Zhou-Zhu \cite{TWZ1, TWZ2, TWZZ}).
Let $g_{\a,0}(\cdot,\cdot)=\w(\cdot,\varphi_\a^*J_{st}\cdot)$ on $U'_\a$, where $\varphi_\a^*J_{st}=J_{\a,0}$ is integrable, $J_{st}$ is the standard complex structure on $\C^n$.
Define
$$
g_\a(t)(\cdot,\cdot)=g_\a(0)e^{th_\a}(\cdot,\cdot)=g_\a(0)(\cdot,e^{tg^{-1}_\a(0)h_\a}\cdot), g_\a(0)=g_{\a,0},
$$
that is,
\beq \label{a15}
g_\a(t)(x)(X,Y)=g_\a(0)(X,Y)+g_\a(0)(X,\sum^{\infty}_{k=1}\frac{t^k(g^{-1}_\a(0)h_\a)^k(Y)}{k!}),
\eeq
where $U'_\a\supset\supset U_\a$, $x\in U_\a$, $X,Y\in T_xU_\a$ and $h_\a$ is $J_\a(0)$-anti-invariant symmetric $(2,0)$-tensor field. We will show that if $U_{\a_1}\cap U_{\a_2}\neq\emptyset$, then $$g_{\a_1}(1)=g_{\a_2}(1)=g_J|_{U_{\a_1}\cap U_{\a_2}}$$ and, in general,
$$(\varphi_{\a_2}\comp\varphi_{\a_1}^{-1})^*g_{\a_2}(0)\neq g_{\a_1}(0)\quad \text{on}\quad U_{\a_1}\cap U_{\a_2}.$$

Since $(M^{2n},\w,J,g_J)$ is a closed almost K\"ahler manifold, without loss of generality, we may assume that $inj(M^{2n},g_J)\geq 3$. Hence, for $\forall p\in M^{2n}$, there exists a Darboux's coordinate chart
$$B(p,3)\supset\supset U_p\supset B(p,1)=\{x\in M^{2n}|\rho_{g_J}(p,x)<1\}$$ and a diffeomorphism $$\psi_p:B(p,3)\to \C^n\cong\R^{2n}$$ such that $\psi_p(U_p)$ is a strictly pseudoconvex domain, by Darboux Theorem $\w|_{U_P}=\psi^*_p\w_0$, and $\psi_p(p)=\{o\}$ in $\C^n$,
where $\rho_{g_J}$ is the distance function defined by the metric $g_J$ and $\w_0$ is the standard K\"ahler form on $\C^n$. It is easy to see that all derivatives of $\psi_p$ (measured with respect to the metric $g_J$ and its
Levi-Civita connection) are bounded. Let $$J_{p,0}=\psi^*_pJ_{st},\quad g_{p,0}(\cdot,\cdot)=\w(\cdot, J_{p,0}\cdot)\quad \text{on}\quad U_p.$$ Then $(\w|_{U_p}, J_{p,0}, g_{p,0})$ is a K\"ahler flat structure on $U_p$ since
$\w|_{U_p}=\psi^*_p\w_0$. Since $\overline{U}_p$ is compact, it is easy to see that
$$||\psi_p||_{C^k(U_p,g_J)}+||\psi^{-1}_p||_{C^k(\psi_p(U_p),(\psi^{-1}_p)^*g_J)}\leq C(M^{2n},p,\w,J,k)<\infty, $$ where $k$ is a nonnegative integer.
Hence
\begin{align} \label{a16}
||J_{p,0}||_{C^k(U_p,g_J)}, ||g_{p,0}||_{C^k(U_p,g_J)},||g^{-1}_{p,0}||_{C^k(U_p,g_J)}\nonumber\\
\leq C(M^{2n},p,\w,J,k)<+\infty.
\end{align}
Thus, one can construct a family of almost K\"ahler structures on $U_p$. Define
$J_p(t)=J_p(0)e^{th_p}$, $J_p(0)=J_{p,0}$, $t\in[0,1]$, $h_p$ is a $J_p(0)$-anti-invariant symmetric $(2,0)$ tensor field on $U_p$ which will be defined later. Notice that let $S_p=-J_p(0)\comp J$, then it is a symmetric positive definite and symplectic matrix on $U_p$, and $$0<C_1(M^{2n},p,\w,J,k)\leq ||S_p||_{C^k(U_p,g_J)}\leq C_2(M^{2n},p,\w,J,k),$$ where $k$ is a nonnegative integer
(cf. McDuff-Salamon \cite[Lemma 2.5.5]{MS}), that is, $$S_pJ_p(0)S_p=J_p(0),\quad \text{and}\quad J|_{U_p}=J_p(0)S_p.$$ Let $h_p=g_p(0)\ln S_p$ which is (2,0) tensor field on $U_p$, where $g_p(0)=g_{p,0}$. Thus the $C^k$-norms of $S_p$ and $h_p$ with respect to $g_J$ and its Levi-Civita connection are
bounded from above by $C(M^{2n},p,\w,J,k)$ which depends only on $M^{2n}$, $p$, $\w$, $J$, and $k$. It is easy to see that $J_p(t)$ are $\w$-compatible for $t\in[0,1]$. If $U_{p_1}\cap U_{p_2}\neq\emptyset$, for any $x\in U_{p_1}\cap U_{p_2}$,
$$J_{p_1}(1)(x)=J_{p_1}(0)S_{p_1}(x)=J|_{U_{p_1}}(x)=J_{p_2}(1)(x).$$
Let $$g_p(t)(\cdot,\cdot)=\w(\cdot,J_p(t)\cdot), \quad t\in[0,1],$$ which are $J_p(t)$-invariant Riemannian metric on $U_p$. If $U_{p_1}\cap U_{p_2}\neq\emptyset$, for any $x\in U_{p_1}\cap U_{p_2}$,
$$g_{p_1}(1)(x)=g_{p_2}(1)(x)=g_J(x).$$
Notice that $\forall x\in U_p$, $X,Y\in T_xU_p$,
\beq \label{a17}
g_p(t)(x)(X,Y)=g_p(0)(x)(X,Y)+g_p(0)(x)(X,\sum_{k=1}^\infty\frac{t^k(g^{-1}_p(0)h_p)^k}{k!}Y).
\eeq
Thus, $(\w|_{U_p},J_p(t),g_p(t))$, $t\in[0,1]$, is a family of almost K\"ahler structures on $U_p$. If $t=0$, then $(\w|_{U_p},J_p(0),g_p(0))$ is a K\"ahler flat structure on $U_p$; if $U_{p_1}\cap U_{p_2}\neq\emptyset$,
then for any $x\in U_{p_1}\cap U_{p_2}$,
$$(\w,J_{p_1}(1),g_{p_1}(1))|_x=(\w,J_{p_2}(1),g_{p_2}(1))|_x=(\w,J,g_J)|_x.$$
Since $[0,1]\times M^{2n}$ is compact, hence, as \eqref{a16}-\eqref{a17}, for $t\in[0,1]$,
\begin{align}\label{a18}
||J_p(t)||_{C^k(U_p,g_J)},||g_p(t)||_{C^k(U_p,g_J)},||g^{-1}_p(t)||_{C^k(U_p,g_J)}\nonumber\\
\leq C(M^{2n},p,\w,J,k)<\infty,
\end{align}
where $C(M^{2n},p,\w,J,k)$ depends only on $M^{2n}$, $p$, $\w$, $J$, and $k$. Let $*_{p,t}$ be the Hodge star operator with respect to the metric $g_p(t)$ on $U_p$, therefore
$$||*_{p,t}||_{C^k(M^{2n},p,\w,J)}, \quad t\in[0,1],$$ are bounded by $C(M^{2n},p,\w,J)<+\infty$. Hence, $$\text{coefficients of }g(t),*_{p,t} \text{ are in }L^{\infty}_k(U_p, g_J)\cap C^{\infty}(U_p)$$ with respect to the smooth metric $g_J$ and its Levi-Civita connection.

On $U_p\subset M^{2n}$, let
$$S_p(t)=-J_p(t)J|_{U_p}=-J|_{U_p}J_p(t),\quad t\in[0,1],$$
which are symmetric positive definite symplectic matrices, and $$0<C_1(M^{2n},p,\w,J,k)\leq ||S_p||_{C^k(U_p,g_J)}\leq C_2(M^{2n},p,\w,J,k)<+\infty,$$ where $k$ is a nonnegative integer (cf. McDuff-Salamon \cite[Chapter 2, Lemma 2.5.5]{MS}). If $t=1$, $S_p(1)$ is $I_{2n}$. If $t=0$, $S_p(0)=S_p$, and then
$$J|_{U_p}=J_p(t)S_p(t), J_p(t)=J|_{U_p}S_p(t),$$ where $S_p(t)$ is a symmetric positive definite matrix.
Let
$$h_p(t)=g_p(t)\ln S_p(t)\quad (resp.\quad h'_p(t)=g_p(1)\ln S_p(t)),\quad t\in [0,1].$$
Hence, for $x\in U_p$, $X,Y\in T_xU_p$, $t\in [0,1]$,
$$g_p(t)(x)(X,Y)=g_p(1)(x)(X,Y)+g_p(1)(x)(X,\sum^{\infty}_{k=1}\frac{(t-1)^k(g^{-1}_p(1)h'_p(t))^kY}{k!}),$$ and $$g_p(1)(x)(X,Y)=g_p(t)(x)(X,Y)+g_p(t)(x)(X,\sum^{\infty}_{k=1}\frac{(1-t)^k(g^{-1}_p(t)h_p(t))^kY}{k!}).$$
By (\ref{a18}), it is easy to see that
\beq \label{a19}
g_p(t)(x)(X,X)\leq C(M^{2n},p,\w,J)g_p(1)(x)(X,X),
\eeq
and
\beq \label{a20}
g_p(1)(x)(X,X)\leq C(M^{2n},p,\w,J)g_p(t)(x)(X,X),
\eeq
where $t\in [0,1]$ and $C(M^{2n},p,\w,J)>1$ depending only on $M^{2n}$, $p$, $\w$, and $J$.

Since $(M^{2n},\w)$ is a closed symplectic manifold of dimension $2n$, one can find an $\w$-compatible almost complex structure $J$ on $M^{2n}$ and set $g_J(\cdot,\cdot)=\w(\cdot,J\cdot)$ such that
$(\w,J,g_J)$ is an almost K\"ahler structure on $M^{2n}$, and choose a finite Darboux's coordinate subatlas $\{U_{\a_i},\psi_{\a_i}\}_{1\leq i\leq N}$ where $$B(p_i,3)\supset\supset U_{\a_i}\supset B(p_i,1), \{p_1,...,p_N\}\subset M^{2n},$$ since $M^{2n}$ is closed.
$$||\psi_{\a_i}||_{C^k(U_{\a_i},g_J)}+||\psi^{-1}_{\a_i}||_{C^k(\psi_{\a_i}(U_{\a_i}),(\psi^{-1}_{\a_i})^*g_J)}\leq C(M^{2n},\w,J,k)<\infty, $$ where $k$ is a nonnegative integer.
Hence, by using the local deformation of $\w$-compatible almost complex structures on Darboux's coordinate chart \cite{TWZ1, TWZ2}, one can make a global deformation of $\w$-compatible almost complex structures on $M^{2n}$ off a Lebesgue measure zero subset. In fact, one can find a partition of $M^{2n}$ (as done in Lipschitz manifolds): Let $W_1=U_{\a_1}$, $T_1=U_{\a_1}$, and as $2\leq i\leq N$,
$$W_i=U_{\a_i}-\bigcup^{i-1}_{j=1}\overline{U}_{\a_j},\quad T_i=U_{\a_i}-\bigcup^{i-1}_{j=1}U_{\a_j},$$
where $W_i$, $1\leq i\leq N$, are open, $$W_i\subset T_i\subset U_{\a_i},\quad T_i\cap T_j=\emptyset\quad \text{for}\quad i\neq j.$$ Moreover, $T_i$, $1\leq i\leq N$, constitutes a disjoint partition of $M^{2n}$, that is, $$M^{2n}=\bigcup\limits^N_{i=1}T_i.$$ Then let
$$\overset{\circ}{M} {}^{2n}(\w,J,P)\coloneqq\bigcup^{N}_{i=1}W_i.$$
$M^{2n}\setminus\overset{\circ}{M} {}^{2n}(\w,J,P)$ has Hausdorff dimension $\leq 2n-1$ with Lebesgue measure zero. We transfer on each $T_i$ the almost K\"ahler metrics $(\varphi_{\a_i}^{-1})^*g_{\a_i}(t)$, $t\in[0,1]$, on $\psi_{\a_i}(U_{\a_i})\subset\R^{2n}$ via the coordinate map $\varphi_{\a_i}$, that is,
$$g_t=\{g_{\a_i}(t)|_{T_i}\}_{1\leq i\leq N}, \quad 0\leq t\leq 1,$$
is a family of measurable Riemannian metrics on $M^{2n}$ (cf. Drutu-Kapovich \cite[p.34]{DK}) which are with Lipschitz condition by \eqref{a19}-\eqref{a20}. Restricted to the open dense submanifold $\overset{\circ}{M} {}^{2n}(\w,J,P)$, $g_t$, $t\in[0,1]$, are smooth metrics which are quasi-isometric due to \eqref{a19}-\eqref{a20}. Hence, $$\text{coefficients of }g_t\text{ are in }C^\infty(\overset{\circ}{M} {}^{2n})\cap L_k^2(M^{2n},g_J),$$ where $k$ is a nonnegative integer. By $g_t$ and $\w$, one can define $g_t$, $\w$-compatible almost complex structures $$J_t=\{J_{\a_i}(t)|_{T_i}\}_{1\leq i\leq N}, 0\leq t\leq1$$
(cf. McDuff-Salamon \cite[Chapter 2, Lemma 2.5.5]{MS}). Thus, by global $\w$-compatible almost complex structures on $M^{2n}$ off the Lebesgue measure zero subset $M {}^{2n}\backslash\overset{\circ}{M} {}^{2n}(\w,J,P)$, one can define a family of almost K\"ahler structures $(\w,J_t,g_t)$, $t\in[0,1]$, on $\overset{\circ}{M} {}^{2n}(\w,J,P)$, where $(\w,J_0,g_0)$ is a K\"ahler flat structure on $\overset{\circ}{M} {}^{2n}(\w,J,P)$ and $(\w,J_1,g_1)=(\w,J,g_J)$ on $M^{2n}$. It is easy to see that $g_t$, $t\in[0,1]$, are a family of measurable almost K\"ahler metrics which are with equivalent norms by Lipschitz condition \eqref{a19}-\eqref{a20}. Hence, $(M^{2n},\w,J_t,g_t)$, $t\in[0,1]$, are called a family of (measurable) Lipschitz almost K\"ahler manifolds.
Let $*_t$, $t\in[0,1]$, be the Hodge star operator defined by the Riemannian metric $g_t$, $t\in[0,1]$, on $\overset{\circ}{M} {}^{2n}(\w,J,P)$. Hence, by (\ref{a17}) and (\ref{a18}),
$$\text{coefficients of }*_t\in C^\infty(\overset{\circ}{M} {}^{2n}(\w,J,P))\cap L^{\infty}_k(M^{2n})$$ with respect to the smooth metric $g_J$ and its Levi-Civita connection.

If $t=1$, $g_1=g_J$ is a smooth almost K\"ahler metric on $M^{2n}$; if $t=0$, $g_0$ is a K\"ahler flat structure on $\overset{\circ}{M} {}^{2n}(\w,J,P)$. By (\ref{a15}), $g_t$ is a family of smooth metrics, $t\in[0,1]$, on $\overset{\circ}{M} {}^{2n}(\w,J,P)$. Note that, in general, $(\w, J_t,g_t)$, $t\in[0,1)$ are not almost K\"ahler structures on $M^{2n}$ which satisfy Lipschitz condition \eqref{a19}-\eqref{a20}. Let $\nabla^1_t$
be the second canonical connection on $\overset{\circ}{M} {}^{2n}(\w,J,P)$ with respect to the metric $g_t$ satisfying
$$\nabla^1_tg_t=0,\quad\nabla^1_tJ_t=0,\quad\nabla^1_t\w=0.$$
By (\ref{a15}) and (\ref{a17})-(\ref{a20}), it is easy to see that
\begin{align} \label{a21}
C^{-1}g_t(X,X)\leq g_1(X,X)\leq Cg_t(X,X),  \nonumber\\
\forall X\in T\overset{\circ}{M} {}^{2n}(\w,J,P), t\in[0,1],
\end{align}
where $C>1$ is a constant depending only on $M^{2n}$, $\w$, and $J$. Note that (\ref{a21}) implies that $(\w,J_t,g_t)$, $0\leq t\leq 1$, is a family of almost K\"ahler structures on $\overset{\circ}{M} {}^{2n}(\w,J,P)$ which are
quasi-isometric by \eqref{a21}, in particular,
$(\w,J_1,g_1)=(\w,J,g_J)$ is an almost K\"ahler structure on $M^{2n}$, $(\w,J_0,g_0)$ is a K\"ahler flat structure on $\overset{\circ}{M} {}^{2n}(\w,J,P)$, $g_t$, $t\in[0,1]$ are equivalent (measurable) Lipschitz Riemannian metrics on $M^{2n}$ which satisfy Lipschitz condition \eqref{a19}-\eqref{a20} (cf. N. Teleman \cite{Tele3} or J. Luukkainen and J. V\"ais\"al\"a \cite{LV}). Notice that $g_t$ is not Lipschitz metric on $M^{2n}$, $t\in[0,1)$. By using $g_t$, $t\in[0,1]$, one can define an inner product $(\cdot, \cdot)_{g_t}$ on $$\Lambda^pT^*\overset{\circ}{M} {}^{2n}(\w,J,P), \quad 0\leq p\leq 2n,$$ with respect to the metric $g_t$, $t\in[0,1]$, as follows:
$$(\a,\b)_{g_t}dvol_{g_t}=\a\wedge*_t\b\in\Lambda^{2n}T^*\overset{\circ}{M} {}^{2n}(\w,J,P), \forall \a,\b\in\Omega^p(\overset{\circ}{M} {}^{2n}),$$
where $*_t$ is the Hodge star operator with respect to the metric $g_t$.

Notice that the volume density on $\overset{\circ}{M} {}^{2n}(\w,J,P)$
$$dvol_{g_t}=\frac{\w^n}{n!}=dvol_{g_J},\quad t\in[0,1],$$
which is a smooth volume form on $M^{2n}$. In fact, if a volume form $dvol_{g_t}$, $t\in[0,1]$, defined by the (measurable) Lipschitz almost K\"ahler metric $g_t$ is smooth (or the volume density $\sqrt{\det g_t}$ is constant), then $g_t$ is able to be smoothing (cf. N. Teleman \cite{Tele1}). In particular, $g_0$ and $g_1$ are equivalent (measurable) Lipschitz Riemannian metrics on $M^{2n}$ satisfying Lipschitz condition \eqref{a19}-\eqref{a20} (cf. N. Teleman \cite{Tele3} or J. Luukkainen and J. V\"ais\"al\"a \cite{LV}). Hence, one can define an inner product on $\Omega^p(M^{2n})$ due to \eqref{a21} by using (measurable) Lipschitz Riemannian metrics $g_t$ since $M^{2n}\setminus\overset{\circ}{M} {}^{2n}(\w,J,P)$ has Hausdorff dimension $\leq 2n-1$ with Lebesgue measure zero.
$$\langle\a,\b\rangle_{g_t}\coloneqq\int_{M^{2n}}(\a,\b)_{g_t}\frac{\w^n}{n!}\coloneqq\int_{\overset{\circ}{M}{}^{2n}}(\a,\b)_{g_t}\frac{\w^n}{n!},$$  where $\forall \a,\b\in\Omega^p(M^{2n}),0\leq p\leq 2n$,$t\in[0,1]$.
Thus, by the same reason, one can define $L^2_k(t)$-norm on $\Omega^p(M^{2n})$ (note that since $M^{2n}$ is closed by \eqref{a18}-\eqref{a21} $g_t\in C^\infty(\overset{\circ}{M} {}^{2n})\cap L_k^\infty(M^{2n}, g_J)$, $k$ nonnegative integer) with respect to the (measurable) Lipschitz Riemannian metrics $g_t$ on $M^{2n}$,
\begin{align} \label{a22}
||\a||^2_{L^2_k(t)}\coloneqq\sum^k_{i=0}\langle(\nabla_t)^i\a|_{\overset{\circ}{M} {}^{2n}(\w,J,P)},(\nabla_t)^i\a|_{\overset{\circ}{M} {}^{2n}(\w,J,P)}\rangle_{g_t},  \nonumber\\
t\in[0,1], \forall \a\in\Omega^p(M^{2n}), 0\leq p\leq 2n ,
\end{align}
where $\nabla_t$ is the Levi-Civita connection with respect to the metrics $g_t$ on $\overset{\circ}{M} {}^{2n}(\w,J,P)$. Hence, we define Hilbert spaces of $\Omega^p(M^{2n})$, $t\in[0,1]$, $0\leq p\leq 2n$ as follows:

\bdefn  Let $(M^{2n},\w,J,g_J)$ be a closed almost K\"ahler manifold of dimension $2n$. Suppose that for any partition of $M^{2n}$, let $(\w,J_t,g_t)$, $t\in[0,1]$, be a family of (measurable) Lipschitz almost K\"ahler structures on $M^{2n}$ by using Darboux's coordinate charts, where $(\w,J_1,g_1)=(\w,J,g_J)$, and $(\w,J_0,g_0)$ is a (measurable) Lipschitz K\"ahler flat structure on $M^{2n}$ which can be regarded as a singular K\"ahler structure on $M^{2n}$. Since $$\text{coefficients of }*_t, g_t\text{ are in }C^\infty(\overset{\circ}{M} {}^{2n})\cap L^2_k(M^{2n},g_J)$$ for $k$ nonnegative integer, $M^{2n}\setminus\overset{\circ}{M} {}^{2n}$ has Lebesgue measure zero, one can define $L^2_k\Omega^p(M^{2n})(t)$ being the completion of $\Omega^p(M^{2n})$$\otimes_{\R}\C$ with respect to the norm
$||\cdot||_{L^2_k(t)}$.
\edefn

Since $$\text{coefficients of }*_t\text{ are in }C^\infty(\overset{\circ}{M} {}^{2n})\cap L^{\infty}_k(M^{2n},g_J), $$by (\ref{a21}), $L^2_k\Omega^p(M^{2n})(t)$, $0\leq p\leq 2n$, are equivalent for $t\in[0,1]$. Since $(\w, J_1,g_1)=(\w, J, g_J)$ is an almost K\"ahler structure on $M^{2n}$ which is quasi-isometric to $g_t$, $t\in[0,1)$, on $\overset{\circ}{M} {}^{2n}$, it is not hard to see that $\Omega^p(M^{2n})\otimes_{\R}\C$ is dense in $L^2_k\Omega^p(M^{2n})(1)$ (cf. T. Aubin \cite[Theorem 2.6 and Remark 2.7]{Au} or \cite{Au0}). By the definition of $L^2_k\Omega^p(M^{2n})(t)$, $t\in[0,1]$, and (\ref{a21}), it is easy to see that $L^2_k\Omega^p(M^{2n})(t)$ and $L^2_k\Omega^p(M^{2n})(1)$, for $t\in[0,1]$, satisfy
\begin{align} \label{a23}
C^{-1}(k)||\a||_{L^2_k(t)}\leq ||\a||_{L^2_k(1)}\leq C(k)||\a||_{L^2_k(t)},\nonumber\\ \a\in\Omega^p(M^{2n}), 0\leq p\leq 2n,
\end{align}
where $k$ is a nonnegative integer and $C(k)>1$ is a constant depending only on $k$, $M^{2n}$, $\w$, $J$. Thus, by \eqref{a23}, $L^2_k\Omega^p(M^{2n})(0)$ and $L^2_k\Omega^p(M^{2n})(1)$ are also called quasi isometry (cf. Drutu-Kapovich \cite[\S 5.1]{DK}).

In summary, we have the following lemma:
\blem\label{lem2.6} Let $(M^{2n},\w,J,g_J)$ be a closed almost K\"ahler manifold of dimension $2n$. Suppose that $(\w,J_t,g_t)$, $t\in[0,1]$, is a family of (measurable) Lipschitz almost K\"ahler structures on $M^{2n}$ with Lipschitz condition \eqref{a19}-\eqref{a20}, where $(\w, J_1,g_1)=(\w, J, g_J)$,
 $(\w,J_0,g_0)$ is a (measurable) Lipschitz K\"ahler flat structure. $$\text{coefficients of }g_t\in C^\infty(\overset{\circ}{M} {}^{2n})\cap L^\infty_k(M^{2n}, g_J),$$ where $k$ is a nonnegative integer.
Moreover, $(\w,J_t,g_t)$, $t\in[0,1]$, are almost K\"ahler structures on $\overset{\circ}{M} {}^{2n}(\w,J,P)$ which are quasi isometry. The $L^2$-norm $||\cdot||_{L^2_k(t)}$ on $\Omega^p(M^{2n})$, $t\in[0,1]$,
 are equivalent, that is, for any $\a\in\Omega^p(M^{2n})$, $0\leq p\leq 2n$, $t\in[0,1]$,
$$ C^{-1}||\a||_{L^2_k(t)}\leq ||\a||_{L^2_k(1)}\leq C||\a||_{L^2_k(t)},$$
where $C>1$ is a constant depending only on $k$, $M^{2n}$, $\w$, $J$. In particular, $L^2_k\Omega^p(M^{2n})(t)$ can be approximated by $\Omega^p(M^{2n})$ with respect to $L^2$-norm $L^2_k\Omega^p(M^{2n}, g_J)$, where $k$ is a nonnegative integer.
\elem
\begin{proof}
Since background metric $g(1)=g_J$ is a smooth Riemannian metric on $M^{2n}$, by the construction of $g(t)$, $t\in [0,1]$, and \eqref{a18}-\eqref{a21}, $$\text{coefficients of }*_t, g(t)\in C^\infty(\overset{\circ}{M} {}^{2n})\cap L_k^\infty(M^{2n}, g_J),$$ where $k$ is a nonnegative integer.
By \eqref{a22} and Definition 2.5, the smooth $\Omega^p(M^{2n})$ is dense in $L_k^2\Omega^p(M^{2n})(t)$ with respect to the metric $g(t)$. Since restricted to $\overset{\circ}{M} {}^{2n}$, $g(t)$ and $g(1)=g_J$ are quasi isometric, $M^{2n}\setminus\overset{\circ}{M} {}^{2n}(\w,J,P)$ has Lebesgue measure zero, hence for any $\a\in\Omega^p(M^{2n})$, $0\leq p\leq 2n$, $t\in[0,1]$, $$ C^{-1}||\a||_{L^2_k(t)}\leq ||\a||_{L^2_k(1)}\leq C||\a||_{L^2_k(t)},$$
where $C>1$ is a constant depending only on $k$, $M^{2n}$, $\w$, $J$. Thus $L^2$-space $L_k^2\Omega^p(M^{2n})(t)$ and $L_k^2\Omega^p(M^{2n})(1)$ are quasi isometric for $0\leq p\leq 2n$, $t\in[0,1]$ and $k$ nonnegative integer. In particular, $L^2_k\Omega^p(M^{2n})(t)$ can be approximated by $\Omega^p(M^{2n})$ with the norm $L^2_k\Omega^p(M^{2n})(1)$, where $0\leq p\leq 2n$ and $k$ is a nonnegative integer. This completes the proof of Lemma \ref{lem2.6}.
\end{proof}

Let $*_t$ be the Hodge star operator with respect to the (measurable) Lipschitz almost K\"ahler metric $g_t$ on $M^{2n}$, $t\in[0,1]$, let $d^{*_t}$, $t\in[0,1]$, be $L^2$-adjoint operator of $d$ in the sense of distributions with respect to the metric $g_t$. Since $$\text{coefficients of }*_t\in C^\infty(\overset{\circ}{M} {}^{2n})\cap L_k^\infty(M^{2n}, g_J),$$ $-*_td*_t$ is well defined on $L_1^2\Omega^p(M^{2n})(t)$ for $t\in[0,1]$. Similar to the smooth case, we have the following lemma:
\blem\label{lem2.7} Let $d^{*_t}$ be $L^2$-adjoint operator of $d$ with respect to the metric $g_t$, then $d^{*_t}=-*_td*_t$, $(d^{*_t})^{*_t}=d$ in the sense of distributions.
Thus, Hodge-Laplacian $\Delta_t=dd^{*_t}+d^{*_t}d$, $t\in[0,1]$, in $L_2^2\Omega^p(M^{2n})(t)$, with respect to the (measurable) Lipschitz metric $g_t$ on $M^{2n}$ is essentially self adjoint.
\elem
\begin{proof} Since due to Lemma \ref{lem2.6} $\Omega^p(M^{2n})$ is dense in $L_k^2\Omega^p(M^{2n})(t)$ with respect to the norm $L^2_k\Omega^p(M^{2n})(1)$, $0\leq p\leq 2n$, $t\in[0,1]$ and $k$ nonnegative integer,
it is enough to choose $\a\in \Omega^{p-1}(M^{2n})\otimes_{\R}\C$ and $\b\in\Omega^p(M^{2n})\otimes_{\R}\C$. Then
$$\a\wedge*_t\b\in\Omega^{2n-1}(\overset{\circ}{M} {}^{2n}(\w,J,P))\bigcap L^{\infty}_1\Omega^{2n-1}(M^{2n})(1).$$
Since coefficients of $*_t$ are in $L^{\infty}_1(M^{2n}, g_J)$, $M^{2n}\setminus\overset{\circ}{M} {}^{2n}(\w,J,P)$ has Lebesgue measure zero,  $$*_t\b\in L^2_1\Omega^{2n-p}(M^{2n})(t),$$
hence,
$$\a\wedge*_t\b\in L^2_1\Omega^{2n-1}(M^{2n})(t).$$
Since $L^2_1\Omega^{2n-p}(M^{2n})(t)$ and $L^2_1\Omega^{2n-p}(M^{2n})(1)$ are quasi-isometric or with equivalent norms, note that $g(1)=g_J$ is a smooth Riemannian metric on $M^{2n}$,
$$*_t\b\in L^2_1\Omega^{2n-p}(M^{2n})(1),\quad d(\a\wedge*_t\b)\in L^2\Omega^{2n}(M^{2n})(1).$$ By Lemma \ref{lem2.6}, there is a smooth sequence $\tilde{\b}_k\in\Omega^{2n-p}(M^{2n})$ such that
$$||*_t\b-\tilde{\b}_k||_{L^2_1(1)}\to0\quad \text{as}\quad k\to\infty.$$ Therefore,
$$\int_{M^{2n}}d(\a\wedge*_t\b)=\lim_{k\to\infty}\int_{M^{2n}}d(\a\wedge\tilde{\b}_k)=0.$$
Thus,
\begin{align} \label{a24}
0&=\int_{M^{2n}}d(\a\wedge*_t\b)  \nonumber\\
&=\int_{M^{2n}}d\a\wedge*_t\b+(-1)^{p-1}\int_{M^{2n}}\a\wedge d(*_t\b)\nonumber\\
&=\int_{M^{2n}}(d\a,\b)_{g_t}dvol_{g_t}-\int_{M^{2n}}(\a,-*_td*_t\b)_{g_t}dvol_{g_t}\nonumber\\
&=\langle d\a, \b\rangle_{g_t}-\langle\a,-*_td*_t\b\rangle_{g_t}.
\end{align}
Hence, by (\ref{a24}), formal adjoint $d^{*_t}$,  $t\in[0,1]$, is $-*_td*_t$, $(d^{*_t})^{*_t}=d$ in Hilbert spaces. Thus, it is easy to see that Hodge-Laplacian $$\Delta_t=dd^{*_t}+d^{*_t}d, t\in[0,1],$$
with respect to the metric $g_t$ is essentially self-adjoint (cf. Section 3 or P. R. Chernoff \cite{Ch}).
This completes the proof of Lemma \ref{lem2.7}.
\end{proof}

For the (measurable) Lipschitz almost K\"ahler metrics $g_t$, $t\in[0,1]$, as done in smooth manifolds, one can define $d+d^{*_t}$-harmonic, $\a\in L^2\Omega^p(M^{2n})(t)$, $0\leq p\leq 2n$, is called $d+d^{*_t}$-harmonic iff $d\a=0$ and $d^{*_t}\a=0$ in the sense of distributions. Similarly, by Lemma \ref{lem2.7} one can also define $\Delta_t=dd^{*_t}+d^{*_t}d$-harmonic, $t\in[0,1]$, $\a\in L^2\Omega^p(M^{2n})(t)$, $0\leq p\leq 2n$, is called $\Delta_t$-harmonic iff $\Delta_t\a=0$ in the sense of distributions.

By Lemma \ref{lem2.7}, as done in smooth case, it is easy to obtain the following lemma:
\blem[cf. Section 3]\label{lem2.8} Let $(M^{2n},\w,J,g_J)$ be a closed almost K\"ahler manifold. One can construct a family of (measurable) Lipschitz almost K\"ahler structures $(\w,J_t,g_t)$, $t\in[0,1]$, on $M^{2n}$, where $(\w,J_1,g_1)=(\w,J,g_J)$ is smooth
on $M^{2n}$, and $(\w,J_0,g_0)$ is a (measurable) Lipschitz K\"ahler flat structure on $M^{2n}$ which is smooth on $\overset{\circ}{M} {}^{2n}(\w,J,P)$. If a $L^2$ $p$-form $\a$ on $M^{2n}$ is $\Delta_t$-harmonic with respect to the metric $g_t$, $t\in[0,1]$, in the sense of distributions, if and only if $\a$ is $d+d^{*_t}$-harmonic with respect to the metric $g_t$ in the sense of distributions.
\elem

By Lemma \ref{lem2.8}, define $\Delta_t$-harmonic spaces as follows:
\begin{align} \label{a25}
\mathcal{H}^p(M^{2n},g_t)&\coloneqq\{\ker\Delta_t|_{L^2\Omega^p(M^{2n})(t)}\}\nonumber\\
&=\{\a\in L^2\Omega^p(M^{2n})(t)|d\a=0=d^{*_t}\a\},
\end{align}
in the sense of distributions, note that, for any $\a\in\mathcal{H}^p(M^{2n},g_t)$, then $*_t\a|_{\overset{\circ}{M} {}^{2n}(\w,J,P)}$ is smooth.

The Laplacian $\Delta_t=dd^{*_t}+d^{*_t}d$, $t\in[0,1]$, is essentially self adjoint, the $k$th Hilbert space of $p$-forms on $M^{2n}$ is equivalent to the following definition:
$$\langle\a,\b\rangle_{L^2_k(t)}\coloneqq\langle(1+\Delta_t)^{\frac{k}{2}}\a,(1+\Delta_t)^{\frac{k}{2}}\b\rangle_{L^2(t)}=\langle\a,(1+\Delta_t)^{k}\b\rangle_{L^2(t)}.$$
Hence, we have the following definition equivalent $L^2_k(M^{2n})(t)$-norm (cf. \eqref{a22}):
\beq\label{a26}
||\a||_{L^2_k(t)}\coloneqq\sqrt{\langle\a,(1+\Delta_t)^{k}\a\rangle_{L^2(t)}}.
\eeq
By Lemma \ref{lem2.8}, we have Hodge-Kodaira decomposition (cf. Theorem \ref{thm}):
\beq \label{a27}
L^2 \W^p (M^{2n})(t) =\mathcal{H}^p (M^{2n},g_t)\oplus\overline{d\Omega^{p-1}(M^{2n})_t}\oplus\overline{d^{*_t}\Omega^{p+1}(M^{2n})_t},
\eeq
where $$\mathcal{H}^p (M^{2n},g_t)=\ker\Delta_t|_{L^2\W^p(M^{2n})(t)},\quad t\in[0,1],$$ in the sense of distributions, $\overline{d\Omega^{p-1}(M^{2n})_t}$ and $\overline{d^{*_t}\Omega^{p+1}(M^{2n})_t}$ are closure of
$d\Omega^{p-1}(M^{2n})$ and $d^{*_t}\Omega^{p+1}(M^{2n})$ with respect to $L^2(t)$-norm respectively. By the norm $||\cdot||_{L^2_k(t)}$, it is easy to obtain that
$$||\a||_{L^2_k(t)}=||\a||_{L^2(t)}\quad \text{if}\quad \a\in\mathcal{H}^p (M^{2n},g_t).$$
By Lemma \ref{lem2.6}, $L^2_k\W^p(M^{2n})(t)$, $t\in[0,1)$, and $L^2_k\W^p(M^{2n})(1)$ are quasi isometry or with equivalent norms (cf. \eqref{a23}), hence, if $\a\in\mathcal{H}^p (M^{2n},g_t)$,
$$||\a||_{L^2_k(1)}\leq C||\a||_{L^2(t)},$$ where $t\in[0,1)$, $k$ is a nonnegative integer and $C$ depends on $k$. Notice that $(\w,J_1,g_1)=(\w,J,g_J)$ is the smooth
almost K\"ahler structure, by Sobolev embedding theorem (cf. T. Aubin \cite{Au}), if $k-l>\frac{n}{2}$, then $\a\in C^l(\Lambda^pT^*M^{2n})$, hence $\a\in\W^p(M^{2n})$. Therefore,
\beq \label{a28}
\mathcal{H}^p (M^{2n},g_t)\subset\W^{p}(M^{2n}),\quad t\in[0,1],0\leq p\leq 2n.
\eeq

One can define a Hilbert cochain complex $L^2_{l-*}\Omega^*(M^{2n})(t)$ as follows (cf. Br\"uning-Lesch \cite{BL}, W. L\"uck \cite{Luck}, N. Teleman \cite{Tele3}).
\begin{align*}
0\to L^2_{l}\Omega^0(M^{2n})(t)&\xrightarrow{d}L^2_{l-1}\Omega^1(M^{2n})(t)\xrightarrow{d}\cdots\\
&\xrightarrow{d}L^2_{l-2n}\Omega^{2n}(M^{2n})(t)\to 0.
\end{align*}
As done in smooth manifolds (cf. T. Aubin \cite{Au}), one can define
\beq \label{a29}
\mathrm{H}^p (M^{2n},g_t,\C)=\left.\ker d|_{L^2\Omega^p(M^{2n})(t)}\middle/\overline{dL^2\Omega^{p-1}(M^{2n})(t)}\right.,
\eeq
$t\in[0,1]$, $0\leq p\leq 2n$, in the sense of distribution. It is easy to see that the Hodge isomorphism
\begin{align}\label{a30}
\chi^p(t):\mathcal{H}^p (M^{2n},g_t)&\to\mathrm{H}^p (M^{2n},g_t,\C),\quad 0\leq p\leq 2n, t\in[0,1],\nonumber\\
\alpha&\mapsto[\alpha]
\end{align}
is an isomorphism.

Since by \eqref{a28} $\mathcal{H}^p (M^{2n},g_t)\subset\W^{p}(M^{2n})$, then
$$\mathrm{H}^p (M^{2n},g_t,\C)=\left.\ker d|_{\Omega^p(M^{2n})}\middle/\overline{d\Omega^{p-1}(M^{2n})(t)}\right.,$$
$0\leq p\leq 2n$, $\overline{d\Omega^{p-1}(M^{2n})(t)}$ is the closure of $d\Omega^{p-1}(M^{2n})(t)$ with respect to the norm $L_1^2\Omega^p(M^{2n})(t)$.
Since $L_1^2\Omega^p(M^{2n})(t)$ and $L_1^2\Omega^p(M^{2n})(1)$ are quasi isomorphic,
$$d\Omega^{p-1}(M^{2n})(t)\cong d\Omega^{p-1}(M^{2n})(1),$$
then
\beq \label{a31}
\mathrm{H}^p (M^{2n},g_t,\C)\cong\mathrm{H}^p (M^{2n},g_J,\C)
\eeq

Define
$$b^p(M^{2n})(t)\coloneqq\dim\mathrm{H}^p (M^{2n},g_t,\C),\quad 0\leq p\leq 2n, t\in[0,1].$$
By (\ref{a31}), for $0\leq p\leq 2n$, $t\in[0,1)$,
\begin{align}\label{a32}
b^p(M^{2n})(t)&=b^p(M^{2n})(1)=\dim\mathrm{H}^p (M^{2n},g_J,\C)\nonumber\\
&=b^p(M^{2n})\quad(\text{p-th Betti number of}\quad M^{2n}).
\end{align}
Hence, the Euler characteristic of $M^{2n}$ is the following equality:
\begin{align}\label{a33}
\chi(M^{2n})(t)&=\sum_{p=0}^{2n}(-1)^pb^p(M^{2n})(t)\nonumber\\
&=\sum_{p=0}^{2n}(-1)^pb^p (M^{2n})\nonumber\\
&=\chi(M^{2n}),\quad t\in[0,1).
\end{align}

In terms of (\ref{a27})-(\ref{a33}), we have the following theorem (cf. Theorem \ref{thm}, N. Teleman \cite[Theorem 4.1]{Tele3}):
\bthm \label{thm2.9}
Let $(M^{2n}, \w, J, g_J)$ be a closed almost K\"ahler manifold. By the deformation of $\w$-compatible almost complex structures in terms of Darboux's coordinate charts, one can construct a family of (measurable) Lipschitz almost K\"ahler structures
$(\w,J_t,g_t)$, $t\in[0,1]$, on $M^{2n}$. Restricted to the open dense submanifold, $\overset{\circ}{M} {}^{2n}(\w,J,P)$, of $M^{2n}$, $(\w,J_t,g_t)$, $t\in[0,1]$, is a family of almost K\"ahler structures which are quasi-isometric, in particular, $(\w,J_1,g_1)=(\w,J,g_J)$ is an almost K\"ahler structure on $M^{2n}$, $(\w,J_0,g_0)$ is a K\"ahler flat structure on $\overset{\circ}{M} {}^{2n}(\w,J,P)$. Then, for any degree, $0\leq p\leq 2n$,\\
(1) there are the strong Kodaira-Hodge decompositions
$$L^2 \W^p (M^{2n})(t) =\mathcal{H}^p(M^{2n}, g_t)\oplus \overline{d\Omega^{p-1}(M^{2n})_t}\oplus \overline{d^{*_t}\Omega^{p+1}(M^{2n})_t},t\in[0,1],$$
and $\mathcal{H}^p(M^{2n}, g_t)$ is smooth on $M^{2n}$;\\
(2) the Hodge homomorphism
\begin{align*}
\chi^p:\mathcal{H}^p(M^{2n},g_t)&\to \mathrm{H}^p(M^{2n}, g_t, \C),\quad t\in[0,1],\\
\alpha&\mapsto[\alpha]
\end{align*}
is an isomorphism;\\
(3) Since $$\mathrm{H}^p(M^{2n}, g_1, \C)\cong\mathrm{H}^p(M^{2n}, g_t, \C) \quad\text{for}\quad t\in[0,1),$$ $$b^p(M^{2n})=b^p(M^{2n},g_J)=b^p (M^{2n},g_t), t\in[0,1),$$
\begin{align*}
\chi(M^{2n})&=\sum_{p=0}^{2n}(-1)^pb^p (M^{2n},g_J)\nonumber\\
&=\sum_{p=0}^{2n}(-1)^pb^p(M^{2n},g_t),\quad t\in[0,1).
\end{align*}
\ethm

Finally, we have the following remark:
\brem (1) Since $\Delta_t$, $t\in[0,1]$, are elliptic operators with $L^{\infty}_k$, $k$ being a nonnegative integer, coefficients on the smooth manifold $M^{2n}$, $\mathcal{H}^p (M^{2n},g_t)$, $0\leq p\leq 2n$, $t\in[0,1)$, are smooth. In general, by using different methods, E. de Giorgi \cite{Giorgi}, J. Nash \cite{Nash} proved the following result: The weak solution of elliptic (parabolic) operators with measurable coefficients are continuous.

(2) For $(\w, J_1, g_1)=(\w,J,g_J)$, one has
$$\W^k(M^{2n})=\underset{i+j=k}\oplus\W^{i,j}(M^{2n}),$$
see D. Yan \cite{Yan}, or L. S. Tseng and S. T. Yau \cite{TY}. Hence, one can define $$\mathcal{H}^{i,j} (M^{2n},g_J)=\{\alpha\in\W^{i,j}(M^{2n})|d\alpha=0=d^*\alpha, d^*=-*_{g_J}d*_{g_J}\}.$$ In general (cf. T. Huang \cite{Huang1,Huang2}), $$\underset{i+j=k}\oplus\mathcal{H}^{i,j} (M^{2n},g_J)\subsetneq\mathcal{H}^k(M^{2n},g_J).$$
\erem

\section{$L^2$-Hodge theory on the universal covering of closed symplectic manifolds}
Let $(M^{2n}, \w)$ be a closed, $2n$-dimensional symplectic manifold with infinite fundamental group $\pi_1(M^{2n})$. Let $J$ be an $\w$-compatible almost complex structure, and $g_J(\cdot,\cdot) =\w(\cdot, J\cdot)$ the associated Riemannian metric. Then $(\w, J, g_J)$ is an almost K\"ahler structure on $M^{2n}$. The universal covering, $\ti M^{2n}$, of $M^{2n}$ has the lifted almost K\"ahler structure $(\ti\w,\ti J, \ti g_J)$. Since the fundamental group $\G=\pi_1(M^{2n})$ is infinite which is regarded
as the deck transformation group, and acts on the universal covering $\ti M^{2n}$ by deck transformation, the lifted  almost K\"ahler structure $(\ti\w,\ti J, \ti g_J)$ is $\G$-invariant (cf. I. Chavel \cite{Cha}).

Let $\W^p_c (\ti M^{2n})$ be the space of smooth $p$-forms on $\ti M^{2n}$ with compact support. There are natural inner product and norm on it, given by
\beq \label{ba}
\li \a, \b\ri_{L^2} =\int_{\ti M^{2n}} \a \wedge *_{\ti g_J} \b,
\eeq
\beq \label{bb}
\| \a \|_{L^2} =\sqrt{\li \a,\a\ri_{L^2} } ,
\eeq
where $\a,\b\in \W^p_c (\ti M^{2n})$. Let $L^2 \W^p (\ti M^{2n})$ be the Hilbert space completion of $\W^p_c (\ti M^{2n})$ with respect to this inner product defined by Riemannian metric $\ti g_J$. Let $\Delta_{\ti g_J}=dd^*+d^*d$ be the Hodge-Laplace operator on $(\ti M^{2n}, \ti g_J)$, where $d^*=-*_{\ti g_J}d*_{\ti g_J}$, $*_{\ti g_J}$ is Hodge star operator with respect to the Riemannian metric $\ti g_J$. Define the space of $L^2$-integrable harmonic $p$-forms on $(\ti M^{2n}, \ti g_J)$ as
\beq \label{bc}
\mathcal{H}^p_{(2)} (\ti M^{2n}) \coloneqq\{ \a\in L^2 \W^p (\ti M^{2n}) \n \Delta_{\ti g_J} \a=0 \},
\eeq
see W. L\"uck \cite{Luck} or P. Pansu \cite{Pan}.

\bprop[cf. J. Dodziuk \cite{Dod},also see \cite{Gr2,Luck}] \label{bd} Let $(M^{2n},\w)$ be a $2n$-dimensional closed non-elliptic symplectic manifold. Then there exists an almost K\"ahler structure $(\w,J,g_J)$ on $M$, and the universal covering $(\ti M^{2n}, \ti\w, \ti J, \ti g_J)$ is also an almost K\"ahler manifold. Similar to the compact case, we have Hodge-Kodaira decomposition of reduced $L^2$-cohomology
\beq \label{be}
L^2 \W^p (\ti M^{2n}) =\mathcal{H}^p_{(2)} (\ti M^{2n})\oplus\overline{dL^2\Omega^{p-1}(\ti M^{2n})}\oplus\overline{d^*L^2\Omega^{p+1}(\ti M^{2n})},
\eeq
where $$\mathcal{H}^p_{(2)} (\ti M^{2n})\subset\Omega^p(\ti M^{2n})\cap L^2\Omega^p(\ti M^{2n}),$$ $\overline{dL^2\Omega^{p-1}(\ti M^{2n})}$ and $\overline{d^*L^2\Omega^{p+1}(\ti M^{2n})}$ are closure of $dL^2\Omega^{p-1}(\ti M^{2n})$ and \linebreak
$d^*L^2\Omega^{p+1}(\ti M^{2n})$ with respect to $L^2$-norm respectively.
\eprop

We can define the analytic $p$-th $L^2$-Betti number by the von Neumann dimension of the finitely generated Hilbert $\G$-module $\mathcal{H}^p_{(2)} (\ti M^{2n})$, see F. J. Murray and J. von Neumann \cite{MN}, M. Atiyah \cite{At}, P. Pansu \cite{Pan} and M. Shubin \cite{Shu}.

Recall that a Hilbert space $\mathcal{H}$ with a unitary action of a discrete group $\G$ is called a {\it $\G$-module\/} if $\mathcal{H}$ is isomorphic to a $\G$-invariant subspace in the space of $L^2$-functions on $\G$ with values in some Hilbert space $\mathcal{H}$. To each  $\G$-module $\mathcal{H}$, one assigns the von Neumann dimension, also called the $\G$-dimension, $\dim_\G \mathcal{H} \in [0,+\infty]$, which is a nonnegative real number or $+\infty$. The precise definition is not important for the moment, but the following properties convey the idea of $\dim_\G \mathcal{H}$ as some kind of size of the ``quotient space" $\mathcal{H}/\G$ (see P. Pansu \cite{Pan}):

\ben
\item $\dim_\G \mathcal{H}\geq 0$, and $\dim_\G \mathcal{H}=0$ if and only if $\mathcal{H}=0$.
\item If $\G$ is a finite group, then $\dim_\G \mathcal{H} =\dim\mathcal{H} /|\G|$.
\item If $\mathcal{H}_0$ is isomorphic to a dense subspace of $\mathcal{H}$, then $\dim_\G \mathcal{H}_0 =\dim_\G \mathcal{H}$.
\item $\dim_\G \mathcal{H}$ is additive in the following sense. Given an exact sequence of Hilbert $\G$-modules
 $$ 0\to \mathcal{H}_1 \to \mathcal{H}_2\to \mathcal{H}_3 \to 0, $$
 one has $\dim_\G \mathcal{H}_2 =\dim_\G \mathcal{H}_1+ \dim_\G \mathcal{H}_3$.
\item If $\mathcal{H}$ equals the space of $L^2$-functions $\G\to H$, then $\dim_\G \mathcal{H} =\dim H$.
\item Continuity: if $\{ \mathcal{H}_j \}_{j=1}^\infty$ is a decreasing sequence of $\G$-invariant subspaces, then
 $$\dim_\G (\bigcap_j\mathcal{H}_j) =\lim_{j\to \infty} \dim_\G \mathcal{H}_j$$
 \item If $\G' \subset \G$ is a subgroup of finite index $d$, then any unitary representation $\mathcal{H}$ of $\G$ becomes a representation of $\G'$, and $\dim_{\G'} \mathcal{H} =d\, \dim_\G \mathcal{H}$.
 \item Normalization: $\dim_\G (\ell^2 (\G)) =1$.
\een
\vskip 0.1cm

We now return to the universal covering $$\pi: (\ti M^{2n},\ti \w, \ti J,\ti g_J) \to (M^{2n},\w, J,g_J),$$ where the fundamental group  of $(M^{2n},\w)$,  $\G= \pi_1(M^{2n})$, is infinite which acts on the universal covering $\ti M^{2n}$ as a deck transformation. We may choose a fundamental domain $F\subset \ti M^{2n}$ for the action of $\G$, where $F$ is an open submanifold of $\ti M^{2n}$ with boundary $\p F$ of Hausdorff dimension $2n-1$, such that $\g_i F\cap \g_j F =\emptyset$ whenever  $\g_i \not= \g_j$, and $\ti M^{2n} =\bigcup_{\g\in \G} \g \bar F$ (cf. M. Atiyah \cite{At}).  Then we have a decomposition $\ti M^{2n} \cong \G \times F$ up to a subset of measure zero, given by $$\g x \mapsto (\g, x) \in \G \times F.$$  Let $$\overset{\circ}{{\ti M}}  {}^{2n}\coloneqq\bigcup_{\g\in \G} \g F.$$ Then
$$\ti M^{2n} \setminus \overset{\circ}{{\ti M}}  {}^{2n} =\bigcup_{\g\in \G} \g (\p F), \quad \text{and}\quad \overset{\circ}{{\ti M}}  {}^{2n} \cong \G \times F.$$
It is easy to see that the action of $\G$ on $\W^p (\overset{\circ}{{\ti M}}  {}^{2n} )$ becomes the action given by $\g_2(\g_1,\a) =(\g_2 \g_1, \a)$, where $\g_1,\g_2 \in\G$, $\a\in \W^p(F)$. This gives a unitary isomorphism
$$ U^p_F: L^2 \W^p(\overset{\circ}{{\ti M}}  {}^{2n}) \to \ell^2(\G) \otimes L^2 \W^p (F), $$
and the action of $\g$ on $L^2 \W^p(\overset{\circ}{{\ti M}}  {}^{2n})$ becomes $L_\g \otimes {\rm I\/}$ on $\ell^2(\G) \otimes L^2 \W^p (F)$, where $L_\g$ is regular left representation of $\g$ on $\ell^2(\G)$ and $I$ is an identity map on $L^2 \W^p (F)$. Hence $$L^2 \W^p(\ti M^{2n}) \cong L^2 \W^p(\overset{\circ}{{\ti M}}  {}^{2n})$$ is a $\G$-module. Thus, in terms of heat kernel, we define the analytic $p$-th $L^2$-Betti number as follows (cf. \cite{At,Luck,Pan,Shu}):
\begin{align*}
b^p_{(2)} (\ti M^{2n}) &\coloneqq \dim_\G \mathcal{H}^p_{(2)} (\ti M^{2n})  \\
&= \lim_{t\to +\infty} \int_F {\rm tr\/}_\C (e^{-t \Delta_{\ti g_J}}(x,x))  d{\rm vol\/}_{\ti g_J} \\
&= \dim_\G \mathcal{H}^p_{(2)} (\overset{\circ}{{\ti M}}  {}^{2n}) =b^p_{(2)} (\overset{\circ}{{\ti M}}  {}^{2n}),
\end{align*}
where $e^{-t \Delta_{\ti g_J}}(x,y)$ is the heat kernel on $\ti M^{2n}$ with respect to the metric $\ti g_{J}$.

\ms Note that $L^2$-Betti numbers share many properties of ordinary Betti numbers. For instance:
\ben
\item Homotopy invariance (cf. M. Gromov \cite[1.1.E]{Gr2}).  $L^2$-Betti numbers are not homotopy invariants for complete noncompact manifolds. For example, Hyperbolic space $\H^{2n}$ and Euclidean space $\R^{2n}$ are homotopy equivalent and diffeomorphic, but their middle dimensional $L^2$-Betti numbers are different.  Yet, the $L^2$-Betti numbers are invariant under bi-Lipschitz homeomorphisms, and more generally, invariant under bi-Lipschitz homotopy equivalence. In particular, singular $L^2$-Betti numbers are $\G$-equivariant homotopy invariants (cf. \cite{Dod,CG1,CG2,CG3}).
\item Finite covering: if $\hat M^{2n} \to M^{2n}$ is a $d$-fold covering, then $b^p_{(2)} (\hat M^{2n}) =d\,  b^p_{(2)} (M^{2n})$.
\item Harmonic forms: $b^p_{(2)} (\ti M^{2n})$ vanishes if and only if
$\mathcal{H}^p_{(2)} (\ti M^{2n})$ is trivial.
\item Continuity: if $\hat M_j \to M$ is a $d_j$-fold covering, and if the sequence $\hat M_j \to M$ converges to the universal covering $\ti M$ in the following sense: every loop in $M$ lifts to an open path in some $\hat M_j$, then $$b_{(2)}^i (\ti M) =\lim_{j\to \infty} \frac1{d_j} b^i (\hat M_j).$$ However, not every manifold admits such a tower of finite coverings. It is the case, if and only if the fundamental group of $M$ is residually finite.

 The actual definition involves reduced $L^2$ cohomology of the universal covering $\ti M$ of $M$. Although the main idea is presented in Murray-von Neumann's theory of type II factors \cite{MN}, the concept is due to M. Atiyah \cite{At} and I. M. Singer \cite{Sing}. $L^2$-Betti numbers are useful tools for topology, as shown by J. Cheeger and M. Gromov \cite{CG3}, see also W. L\"uck \cite{Luck}.
\een

Let $$\pi: (\ti M^{2n},\ti\w, \ti J, \ti g_J) \to (M^{2n},\w,J,g_J)$$ be the universal covering, where $(M^{2n},\w)$ is a  closed symplectic manifold with infinite fundamental group. By Atiyah's $\G$-index theorem in \cite{At}, we have the following theorem:

\bthm \label{bh}
$$ \chi_{(2)} (\ti M^{2n}) \coloneqq \sum_{p=0}^{2n} (-1)^p b_{(2)}^p (\ti M^{2n}) =\chi(M^{2n}),$$
where $b_{(2)}^p (\ti M^{2n}) =\dim_\G \mathcal{H}_{(2)}^p (\ti M^{2n}) $.
\ethm

Property (4) above and Theorem \ref{bh} imply that a positive answer to the Chern-Hopf conjecture would follow from a vanishing theorem for $\mathcal{H}_{(2)}^i (\ti M^{2n})$ when $i\not= n$, and a nonvanishing theorem for $\mathcal{H}_{(2)}^n (\ti M^{2n})$. The answer is positive for rotationally symmetric metric (cf. J. Dodziuk \cite{Dod1}) and symmetric spaces (cf. A. Borel \cite{Bor}).
\vskip 0.1cm
For non-elliptic K\"ahler manifolds, we have the following theorem due to X. Cheng, Jost-Zuo, M. Gromov (\cite{Chg,JZ,Gr2}, K\"ahler hyperbolic case) and Cao-Xavier, N. Hitchin (\cite{CX,Hit}, K\"ahler parabolic case):

\bthm \label{cc}
Let $(M^{2n}, \w, J,g_J)$ be a non-elliptic K\"ahler manifold of dimension $2n$, and $(\ti M^{2n}, \ti\w, \ti J, \ti g_J)$ its universal covering. Then $\mathcal{H}^p_{(2)} (\ti M^{2n}) =0$ for $p\not= n$. Therefore,
for K\"ahler hyperbolic case
$$ (-1)^n \chi (M^{2n}) =(-1)^n \chi_{(2)} (\ti M^{2n}) > 0;$$
for K\"ahler parabolic case
$$ (-1)^n \chi (M^{2n}) =(-1)^n \chi_{(2)} (\ti M^{2n}) \geq 0.$$
\ethm

Note that the proof of the above theorem is based on the identity $$[L_{\ti\w}, \Delta_{\ti g_J}]=0$$ which implies the hard Lefschetz condition of K\"ahler manifolds. Where Lefschetz maps \cite{GH,Yan} are defined as follows:
\begin{align*}
L_\w:  \W^p(M^{2n}) &\to \W^{p+2} (M^{2n}) , \quad \b\mapsto \w\wedge\b; \\
L_{\ti\w}:  \W^p( \ti M^{2n}) &\to \W^{p+2} (\ti M^{2n}) , \quad \g\mapsto \ti\w \wedge\g.
\end{align*}
But, in general, non-elliptic almost K\"ahler manifold $(M^{2n}, \w, J, g_J)$ is not K\"ahler which has no identity $[L_{\w}, \Delta_{g_J}]=0.$

When $(M^{2n}, \w, J,g_J)$ is a K\"ahler manifold, it is well known that
$$ [L_\w, \Delta_{g_J}] =0, \quad [L_{\ti \w}, \Delta_{\ti g_J}]=0, $$
and the induced maps
\begin{align*}
L^k_\w: H^{n-k}_{{\rm dR\/}} (M^{2n}) &\to H^{n+k}_{{\rm dR\/}} (M^{2n}) ; \\
L^k_{\ti\w}: H^{n-k}_{(2)} (\ti M^{2n}) &\to H^{n+k}_{(2)} ( \ti M^{2n}),
\end{align*}
where $0\leq k\leq n$, are isomorphisms by the hard Lefschetz condition due to the identity $[L_\w, \Delta_{g_J}] =0$ (resp. $[L_{\ti \w}, \Delta_{\ti g_J}]=0$) (cf. Griffith-Harris \cite[p.122]{GH} or D. Yan \cite[Corollary 2.9]{Yan}).
It is well known that for a closed symplectic manifold $(M^{2n}, \w)$ with the hard Lefschetz property, every de Rham cohomology
$H_{dR}^*(M^{2n})$ class contains a symplectic harmonic form (cf. D. Yan \cite[Theorem 0.1]{Yan}). Tan-Wang-Zhou \cite{TWZ2} used symplectic cohomology introduced by Tseng-Yau \cite{TY} to study symplectic parabolic
manifold $(M^{2n}, \w)$ which satisfies the hard Lefschetz condition. Based on this, they obtained that if  $(M^{2n}, \w)$ is a symplectic parabolic manifold with the hard Lefschetz property, then
its signed Euler characteristic satisfies the inequality $(-1)^n\chi(M^{2n})\geq 0$. Note that all
solvmanifolds are aspherical, there are many symplectic solvmanifolds with the hard Lefschetz property, no K\"ahler structure, for example, $M^6(c)$, $N^6(c)$, $P^6(c)$,
more details see \cite{FMS}. In this paper, Theorem \ref{main} does not use the hard Lefschetz condition, this answers a question of Tan-Wang-Zhou \cite[Question 1.7]{TWZ2}.

\vskip 0.2 cm
Recall that if $M^{2n}$ is a compact manifold of dimension $2n$, the set $\mathfrak{M}$ of all Riemannian metrics may be given a Riemannian metric (cf. D. Ebin \cite{Ebin}): For symmetric tensor fields
$S, T$ of type (2,0)
$$\li S,T\ri_g=\int_{M^{2n}}S_{ij}T_{kl}g^{ik}g^{jl}dvol_g.$$
The geodesics in $\mathfrak{M}$ are the curves $g_t=ge^{St}$ which is computed by $g_t(X,Y)=g(X,e^{St}Y),$ where $S$ is a symmetric tensor field of type (2,0) and $e^{St}$ acts on $Y$ as a tensor field of type (1,1). Let $\mathfrak{H}$
be the set of all metrics with the same volume element which is totally geodesics in the set $\mathfrak{M}$. Let $(M^{2n},\w)$ be a compact symplectic manifold of dimension $2n$, where $\w$
is a given symplectic form. Let $\mathfrak{A}(\w)$ be the set of all associated metrics with respect to the given symplectic form on $M^{2n}$. Then $\mathfrak{A}(\w)$ is a totally geodesic in
$\mathfrak{H}$ in the sense that if $h$ is $J$-anti-invariant symmetric (2,0)-tensor field, where $J$ is an $\w$-compatible almost complex structure on $M^{2n}$ and $g(\cdot,\cdot)=\w(\cdot, J\cdot)$, the $g_t=ge^{ht}$ lies in $\mathfrak{A}(\w)$. In particular, two metrics in $\mathfrak{A}(\w)$ may be jointed by a unique geodesic (see D. E. Blair \cite[Propositions p.304,307]{Blair}).

\vskip 0.2 cm
In order to investigate non-elliptic symplectic manifolds without the hard Lefschetz property as done in Section 2, we will use the global deformation of $\w$-compatible almost complex structures to K\"ahler flat structures on an open dense submanifold of the given symplectic manifold. We can obtain a family of (measurable) Lipschitz almost K\"ahler structures. With respect to those constructed (measurable) Lipschitz almost K\"ahler metrics, the corresponding $L^2_k$-norms of $p$-forms are quasi isometry (that is, with equivalent norms).

Recall that suppose that $(M^{2n}, \w, J, g_J)$ is a closed almost K\"ahler manifold with an infinite fundamental group $\G=\pi_1(M^{2n})$.
We can construct a special fundamental domain, $F$, using Darboux coordinate charts (cf. McDuff-Salamon \cite{MS}), such that the $\w$-compatible almost complex structure $J$ on $F$ can be deformed to an $\w$-compatible integrable complex structure $J_0$ on $F$.

\vskip 0.2 cm
We now construct the fundamental domain of $\G$-manifold $\ti M^{2n}$. For a closed almost K\"ahler manifold $(M^{2n}, \w, J, g_J)$, as done in Section 2, one can find a disjoint partition $P(\w,J)$: Choose a finite Darboux's coordinate subatlas $\{U_j,\varphi_j\}_{1\leq j\leq N}$.
Let $W_1=U_1$, $T_1=U_1$, and as $j\geq2$, $$W_j =U_j -\bigcup\limits_{i=1}^{j-1} \bar U_i, \quad T_j =U_j -\bigcup\limits_{i=1}^{j-1}U_i.$$
Where $W_j$, $1\leq j\leq N$, are open sets of $M^{2n}$,
$$W_j\subset T_j\subset U_j,\quad T_j\cap T_k=\emptyset\quad \text{for}\quad j\neq k.$$
Hence, $\{T_j\}_{1\leq j\leq N}$ constitutes a disjoint partition of $M^{2n}$, that is, $M^{2n}=\bigcup\limits_{j=1}^N T_j$. Then, let
$$\overset{\circ}{M} {}^{2n}=\overset{\circ}{M} {}^{2n}(\w,J,P)\coloneqq\bigcup_{j=1}^N W_j.$$
$M^{2n}\backslash\overset{\circ}{M} {}^{2n}$ has Hausdorff dimension $\leq2n-1$ with Lebesgue measure zero. By (\ref{a17}), let
$$g_t=\{g_j(t)|_{T_j}\}_{1\leq j\leq N},\quad 0\leq t\leq 1,$$
be a family of (measurable) Lipschitz almost K\"ahler metrics on $M^{2n}$ with (measurable) Lipschitz almost K\"ahler structures $(\w, J_t, g_t)$, $t\in[0, 1]$, and Lipschitz condition \eqref{a19}-\eqref{a20}, where $(\w, J_0, g_0)$ is a K\"ahler flat structure on $\overset{\circ}{M} {}^{2n}$, and
$(\w, J_1, g_1)=(\w, J, g_J)$ on $M^{2n}$.

Since $\pi:\ti M^{2n}\to M^{2n}$ is the universal covering of $M^{2n}$ with deck transformation group $\Gamma=\pi_1(M^{2n})$, without loss of generality, one may assume that each open component of $\pi^{-1}(U_j)\subset \ti M^{2n}$ is diffeomorphic to $U_j$, $1\leq j\leq N$. For each $U_j$, $1\leq j\leq N$, choose an open component $\ti U_j$ of $\pi^{-1}(U_j)$ such that $\pi(\ti U_j)=U_j$. Define
$$\widetilde{W}_j=\pi^{-1}(W_j)\cap \ti U_j, \ti T_j =\pi^{-1}(T_j)\cap \ti U_j,\quad \text{for}\quad j\neq k, \ti T_j\cap \ti T_k=\emptyset.$$
Let $F=\bigcup\limits_{j=1}^N \widetilde{W}_j$  as a fundamental domain of the universal covering $$\pi: \ti M^{2n} \to M^{2n}.$$ Define
$$\overset{\circ}{\ti M} {}^{2n}=\overset{\circ}{\ti M} {}^{2n}(\w,J,P)\coloneqq\bigcup_{\gamma\in \Gamma} \gamma F$$
which is an open dense $\G$-submanifold of $\ti M^{2n}$ and $\ti M^{2n}\backslash\overset{\circ}{\ti M} {}^{2n}$ has Lebesgue measure zero. We can make a global deformation of $\ti\w$-compatible almost complex structures on $\ti M^{2n}$ off a Lebesgue
measure zero subset $\ti M^{2n}\backslash\overset{\circ}{\ti M} {}^{2n}$, define
$$\ti g_t=\pi^*g_t,\quad 0\leq t\leq 1,$$
is a family of (measurable) Lipschitz almost K\"ahler $\G$-metrics on $\ti M^{2n}$.

On $$\overset{\circ}{\ti M} {}^{2n}\coloneqq\bigcup\limits_{\g\in \G} \g F \subset \ti M^{2n},$$ we have the lifted K\"ahler flat structure $(\ti\w, \ti J_0, \ti g_{J_0})$ on $\overset{\circ}{\ti M}{}^{2n}$. Note that
$$\ti M^{2n}=\bigcup_{\g\in \G} \g \bar F,$$ where $$\ti M^{2n}\backslash\overset{\circ}{\ti M} {}^{2n} =\bigcup\limits_{\g\in \G} \g  \p F =\p \overset{\circ}{\ti M} {}^{2n}$$ is a noncompact,  piecewise smooth hypersurface in $\ti M^{2n}$.

We have two different almost K\"ahler structures on the open dense manifold $\overset{\circ}{\ti M} {}^{2n}$. One is  $(\overset{\circ}{\ti M} {}^{2n}, \ti\w, \ti J, \ti g_J)$, which is inherited from the universal covering $(\ti M^{2n}, \ti\w, \ti J, \ti g_J)$ as an open dense $\G$-submanifold. The other one is induced from the standard K\"ahler structure on $\R^{2n}\cong\C^n$, that is $(\overset{\circ}{\ti M} {}^{2n}, \ti\w, \ti J(0), \ti g(0))$ which is a K\"ahler flat $\G$-manifold. The latter complex structure $\ti J(0)$ in general has no continuous extension to $\ti M^{2n}$.

By equation \eqref{a17} and the definition of $\ti g(t)$, as done in Section 2 we have
\beq \label{ck}
\ti g_{J} |_{\g F} =(\ti g(0) |_{\g F}) \, e^h,
\eeq
where $h=\pi^*h_j$ on $\ti U_j$ which is a symmetric $\ti J$-anti-invariant $(2,0)$ tensor. Hence on $\overset{\circ}{\ti M} {}^{2n}$, by \eqref{ck} we can construct a family of $\G$-invariant almost K\"ahler structures
\beq \label{cl}
\ti g (t) |_{\g F} \coloneqq(\ti g(0) |_{\g F} )\, e^{t\,h},\ \text{for}\ t\in [0,1], \ \g\in \G.
\eeq
Hence, we have a family of $\G$-invariant almost K\"ahler structures $(\ti\w, \ti J(t), \ti g(t))$ on $\overset{\circ}{\ti M} {}^{2n}$, $t\in [0, 1]$. In particular, $t=1$, $\ti g(1)=\ti g_J$, $\ti J(1)=\ti J$, on $\overset{\circ}{\ti M} {}^{2n}$. Since $(\ti\w, \ti J, \ti g_J)$ is a smooth $\G$-invariant almost K\"ahler structure on $\ti M {}^{2n}$, thus $(\ti\w, \ti J(1), \ti g(1))$ can be regarded as a smooth $\G$-invariant almost K\"ahler structure on $\ti M {}^{2n}$. If $t=0$, $(\ti\w, \ti J(0), \ti g(0))$ is a $\G$-invariant K\"ahler flat structure on the open dense manifold $\overset{\circ}{\ti M} {}^{2n}$,
hence $(\ti\w, \ti J(0), \ti g(0))$ can be regarded as a $\G$-invariant (measurable) Lipschitz K\"ahler flat structure on $\ti M^{2n}$ (but, in general, $\ti g(0)$ is not Lipschitz metric).
Thus $(\ti\w, \ti J(t), \ti g(t))$, $t\in[0,1]$ can be regarded as $\G$-invariant (measurable) Lipschitz almost K\"ahler structures on $\ti M^{2n}$ since $\ti M {}^{2n}\backslash\overset{\circ}{\ti M} {}^{2n}$ has Lebesgue measure zero (see C. Drutu and M. Kapovich \cite[p.34]{DK}), in fact, $(\ti\w, \ti J(t), \ti g(t))$ can be regarded as lifting of $(\w, J(t), g(t))$ on $M^{2n}$ defined in Section 2 by universal covering $\pi:\ti M^{2n}\to M^{2n}$; $(\ti\w, \ti J(1), \ti g(1))=(\ti\w, \ti J, \ti g_J)$ is a smooth $\G$-invariant almost K\"ahler metric on $\ti M^{2n}$; $(\ti\w, \ti J(0), \ti g(0))$ is a $\G$-invariant (measurable) Lipschitz K\"ahler flat metric on $\ti M^{2n}$ which is regarded as a singular K\"ahler metric, and a lifting of a K\"ahler structure ($\w, J(0), g(0)$) on the underlying manifold $\pi(F)\subset M^{2n}$ (see Section 2 or N. Teleman \cite[\S3]{Tele3} and Drutu-Kapovich \cite[Chapter 1,2]{DK}).
Since$\{\pi(\ti U_j)\}_{1\leq j\leq N}$ is an open cover of $M^{2n}$, and by \eqref{a18}-\eqref{a20} due to \eqref{cl} which is similar to \eqref{a15}, for $t\in[0, 1]$, $||\ti J(t)|_{\ti U_j}||_{C^k(\ti U_j, \ti g_J)}$, $||\ti g(t)|_{U_j}||_{C^k(U_j, \ti g_J)}$, and $||\ti g^{-1}(t)|_{\ti U_j}||_{C^k(\ti U_j, \ti g_J)}$ are bounded from above by constants $C(M^{2n}, \w, J, k)$ where $k$ is a nonnegative integer, and $\ti g(t)$ is a family of $\G$-invariant (measurable) Lipschitz metrics on $\ti M^{2n}$ which are with equivalent $\G$-invariant norms (cf. Section 2). It follows that for $t\in[0,1]$,
\begin{align} \label{cb}
&||\ti J(t)||_{C^k(\overset{\circ}{\ti M} {}^{2n}, \ti g_J)}+||\ti g(t)||_{C^k(\overset{\circ}{\ti M} {}^{2n}, \ti g_J)}+||\ti g^{-1}(t)||_{C^k(\overset{\circ}{\ti M} {}^{2n}, \ti g_J)}\nonumber\\
&<C(M^{2n}, \w, J, k).
\end{align}
\vskip 6pt
Recall the construction of $\ti g(t)$, $t\in[0, 1]$, notice that $\{\pi(\ti U_j)\}_{1\leq j\leq N}$ is a finite open cover of $M^{2n}$, we have the following estimates on $\overset{\circ}{\ti M}{}^{2n}$ by \eqref{cb}:
\beq \label{cn}
C^{-1}\ti g(t)(x)(X,X)\leq\ti g(1)(x)(X,X)\leq C\ti g(t)(x)(X,X), t\in [0, 1],
\eeq
where $\forall x\in \overset{\circ}{\ti M} {}^{2n}$, $\forall X\in T_x\overset{\circ}{\ti M} {}^{2n}$, $C>1$ is a constant depending only on $M^{2n}$, $\w$, $J$ and $k$, $k$ is a nonnegative integer. It is easy to see that on $\overset{\circ}{\ti M} {}^{2n}$, $\ti g(1)$ and $\ti g(t)$, $t\in[0, 1)$, are quasi isometric and  with $\G$-invariant Lipschitz condition \eqref{cn}. It follows that on $\overset{\circ}{\ti M} {}^{2n}$, $$\text{coefficients of }\ti J(t), \ti g(t)\text{ are in }C^\infty(\overset{\circ}{\ti M} {}^{2n})\cap L_k^\infty(\ti M^{2n}, \ti g_J),$$ where $k$ is a nonnegative integer.
By using $\ti g(t)$, $t\in [0, 1]$, one can define an inner product $(\cdot,\cdot)_{\ti g(t)}$ on $\Lambda^pT^*\overset{\circ}{\ti M} {}^{2n}$, $0\leq p\leq 2n$, with respect to  the metric $\ti g(t)$, $t\in[0, 1]$ as follows:
$$
(\a, \b)_{\ti g(t)}dvol_{\ti g(t)}\coloneqq\a \wedge \ti *_{t} \b \in\Lambda^{2n}T^*\overset{\circ}{\ti M} {}^{2n}, \quad \forall \a,\b\in \Omega^p(\overset{\circ}{\ti M} {}^{2n}),
$$
where $\ti *_{t}$ is the Hodge star operator with respect to the metric $\ti g(t)$, $t\in[0, 1]$ (cf. \cite{Cha}). Notice that $$dvol_{\ti g(t)}=\frac{\ti \w^n}{n!},\quad t\in[0, 1],$$
which are smooth $\G$-invariant volume form on $\ti M^{2n}$. In particular, $\ti g(0)$
and $\ti g(1)$ are equivalent $\G$-invariant (measurable) Lipschitz Riemannian metrics on $\ti M^{2n}$, that is, with $\G$-invariant Lipschitz condition \eqref{cn} (cf. N. Teleman \cite[\S3]{Tele3}).  Since $$\text{coefficients of }\ti g(t)\text{ are in }C^\infty(\overset{\circ}{\ti M} {}^{2n})\cap L_k^\infty(\ti M^{2n}),t\in[0,1],$$ where $k$ is a nonnegative integer, $\ti M {}^{2n}\backslash\overset{\circ}{\ti M} {}^{2n}$ has Hausdorff dimension $\leq2n-1$ with Lebesgue measure zero,
hence, we can define an inner product on
$\Omega^p(\ti M {}^{2n})$, $\forall \a,\b\in \Omega_c^p(\ti M {}^{2n})$, $0\leq p\leq 2n$,
$$
\li\a, \b\ri_{\ti g(t)}\coloneqq \int_{\overset{\circ}{\ti M}{}^{2n}}(\a, \b)_{\ti g(t)}\frac{\ti \w^n}{n!}, \quad t\in [0, 1].
$$
Thus, by the same reason, we can define $L^2_k(t)$-norm on $\Omega_c^p(\ti M {}^{2n})$, $k$ nonnegative integer, with respect to the metric $\ti g(t)$ on $\overset{\circ}{\ti M} {}^{2n}$, $t\in[0, 1]$, $\forall\a\in\Omega_c^p(\ti M {}^{2n})$
\beq \label{cn1}
\|\a\|^2_{L^2_k(t)}\coloneqq \sum_{i=0}^{k}\li(\ti\nabla(t))^i\a|_{\overset{\circ}{\ti M} {}^{2n}}, (\ti\nabla(t))^i\a|_{\overset{\circ}{\ti M} {}^{2n}}\ri_{\ti g(t)},
\eeq
where $\ti\nabla(t)$ is the Levi-Civita connection with respect to the metric $\ti g(t)$ on $\overset{\circ}{\ti M} {}^{2n}$, $0\leq p\leq 2n$, and $t\in[0, 1]$. Hence we are able to define Hilbert spaces of $\Omega^p(\ti M {}^{2n})$,
$0\leq p\leq 2n$, $t\in[0, 1]$ as follows:
\bdefn
Let $(\ti M {}^{2n}, \ti \w)$ be the universal covering of $(M {}^{2n}, \w)$. Then there exists a family of smooth almost K\"ahler structures $(\ti\w, \ti J(t), \ti g(t))$ on $\overset{\circ}{\ti M} {}^{2n}$ (or (measurable) Lipschitz almost K\"ahler structures $(\ti\w, \ti J(t), \ti g(t))$ on $\ti M^{2n}$), $t\in [0, 1]$.
$L^2_k\Omega^p(\ti M^{2n})(t)$ is the completion of $\Omega^p_c(\ti M^{2n})\otimes_\mathds{R}\mathds{C}$ with respect to the norm $\|\cdot\|_{L^2_k(t)}$.
\edefn

For a family of $\G$-invariant (measurable) Lipschitz Riemannian metrics $\ti g(t)$, $t\in[0, 1]$, on $\ti M^{2n}$, if $t=1$, $\ti g(1)=\ti g_J$ is a smooth $\G$-invariant Riemannian metric on $\ti M^{2n}$.
We have the following lemma:
\blem  [cf. Lemma \ref{lem2.6}]\label{lem7}
(1) The Hilbert spaces $$L^2_k\Omega^p(\ti M^{2n})(t), t\in[0, 1), \quad and \quad L^2_k\Omega^p(\ti M^{2n})(1)$$ are quasi-isometric (that is, with Lipschitz condition \eqref{cn}) for $0\leq p\leq 2n$, where $k$ is a nonnegative integer;\\
(2) the coefficients of Hodge star operator $\ti *_t$ with respect to
the metric $\ti g(t)$, $t\in[0, 1]$, are in $L^\infty_k(\ti M^{2n}, \ti g(1))\cap C^\infty(\overset{\circ}{\ti M} {}^{2n})$, where $k$ is a nonnegative integer.
\elem
\begin{proof}
It is clear that $\G$-invariant Riemannian metrics $\ti g(t)$, $t\in[0,1]$ on $\overset{\circ}{\ti M} {}^{2n}$ which are quasi-isometric can be regarded as $\G$-invariant metrics on $\Lambda^pT^*\overset{\circ}{\ti M} {}^{2n}$, $0\leq p\leq 2n$. Hence
the inequalities \eqref{cn} still hold. Consider $\a\in L^2\W^1(\ti M^{2n})(1)$. By \eqref{cn} we have
$$C^{-1}\|\a\|_{L^2\W^1(\ti M^{2n})(t)}\leq \|\a\|_{L^2\W^1(\ti M^{2n})(1)}\leq C\|\a\|_{L^2\W^1(\ti M^{2n})(t)},$$
where $t\in[0,1)$, $C>1$ depends only on $\w, J, M^{2n}$. Let $\ti\nabla(t)$ be the Levi-Civita connection on $\overset{\circ}{\ti M} {}^{2n}$ with respect to the metric $\ti g(t)$, $t\in[0,1]$ (cf. I. Chavel \cite{Cha}). $\forall q\in \overset{\circ}{\ti M} {}^{2n}$, there exists coordinate chart $\{(x_1, \cdots, x_{2n})\}\subset \R^{2n}$ on the neighborhood of $q$. For the Levi-Civita connection $\ti\nabla(t)$, $t\in[0,1]$, the corresponding
Christoffel symbols have following
$$
\ti\nabla(t)_{\frac{\partial}{\partial x_k}}\frac{\partial}{\partial x_j}=\Gamma^l_{jk}(t)\frac{\partial}{\partial x_l}, \quad \Gamma^l_{jk}(t)=\Gamma^l_{kj}(t),
$$
and
$$
\Gamma^l_{jk}(t)=\frac{1}{2}\ti g^{lr}(t)\left[\frac{\partial}{\partial x_j}\ti g^{rk}(t)+\frac{\partial}{\partial x_k}\ti g^{rj}(t)-\frac{\partial}{\partial x_r}\ti g^{jk}(t)\right].
$$
By \eqref{cb}, we have
\beq \label{cn3}
\|\Gamma^l_{jk}(t)\|_{C^k(\overset{\circ}{\ti M} {}^{2n},\ti g_J)}\leq C(M^{2n}, \w, J, k),
\eeq
where $k$ is a nonnegative integer and $C(M^{2n}, \w, J, k)>1$ depends only on $M^{2n}, \w, J$, and $k$. In local coordinate chart $\{(x_1, \cdots, x_{2n})\}$, 1-form $\a$ can be written as $\a=\a_idx_i$.
By \eqref{cn} and \eqref{cn3}, for $t\in[0,1]$, we have
\begin{align*}
|\ti\nabla(t)\a|^2_{\ti g(t)}&\leq C_1(|\ti\nabla(1)\a|^2_{\ti g(1)}+|(\ti\nabla(t)-\ti\nabla(1))\a|^2_{\ti g(1)})\\
&=C_1(|\ti\nabla(1)\a|^2_{\ti g(1)}+|(\Gamma^l_{jk}(t)-\Gamma^l_{jk}(1))\a_l|^2_{\ti g(1)})\\
&\leq C_1(|\ti\nabla(1)\a|^2_{\ti g(1)}+|\Gamma^l_{jk}(t)-\Gamma^l_{jk}(1)|^2_{\ti g(1)}|\a_l|^2_{\ti g(1)})\\
&\leq C_1(|\ti\nabla(1)\a|^2_{\ti g(1)}+|\a|^2_{\ti g(1)}),
\end{align*}
Similarly,
\begin{align*}
|\ti\nabla(1)\a|^2_{\ti g(1)}\leq C_1(|\ti\nabla(t)\a|^2_{\ti g(t)}+|\a|^2_{\ti g(t)}).
\end{align*}
where $C_1>1$ depends only on $\w, J, M^{2n}$. By \eqref{cn1}, \eqref{cn3} and Definition 3.4 we have
$$C^{-1}\li\a, \a\ri_{L^2_1(t)}\leq \li\a, \a\ri_{L^2_1(1)}\leq C\li\a, \a\ri_{L^2_1(t)}, t\in[0, 1)$$where $C>1$ depends only on $\w, J, M^{2n}$. From the discussion above, by the induction,
it is clear that the Hilbert space $L^2_k\Omega^p(\ti M^{2n})(t)$, $t\in [0, 1]$, are equivalent, for $0\leq p\leq2n$. More precisely, $L^2_k\Omega^p(\ti M^{2n})(t)$, $t\in [0, 1)$,
and $L^2_k\Omega^p(\ti M^{2n})(1)$ are quasi-isometric or with equivalent norms for $0\leq p\leq 2n$. Hence, we have
\beq \label{cn2}
C(M^{2n}, \w, J, k)^{-1}\li\a, \a\ri_{L^2_k(t)}\leq\li\a, \a\ri_{L^2_k(1)}\leq C(M^{2n}, \w, J, k)\li\a, \a\ri_{L^2_k(t)},
\eeq
for any $ \a\in\Omega^p_c(\ti M^{2n})$, $t\in [0, 1)$, $0\leq p\leq 2n$, where $k$ is a nonnegative integer, and $C(M^{2n}, \w, J, k)>1$ is a constant depending only on $k, M^{2n}, \w$, and $J$.

According to the discussion before, for the Hodge star operator $\ti*_t=*_{\ti g(t)}$ defined by $\ti g(t)$ on $\ti U_j$, $1\leq j\leq N$, by \eqref{a17}, \eqref{a18} and \eqref{cn3}, we have the following estimates:
$$
||\ti *_t||_{C^k(\ti U_j, \ti g_J)}\leq C(M^{2n}, \w, J, k),
$$
where $C(M^{2n}, \w, J, k)$ depends only on $M^{2n}, \w, J, k$, $k$ is a nonnegative integer, and $t\in [0, 1]$.
Notice that $\{\pi(\ti U_j)\}_{1\leq j\leq N}$ is a finite open cover of $M^{2n}$. Hence, the coefficients of the Hodge star operator $\ti *_t$ are in $$L^\infty_k(\ti M^{2n}, \ti g_J)\cap C^\infty(\overset{\circ}{\ti M} {}^{2n})$$ with respect to the metric $\ti g_J$ and its Levi-Civita connection.

This completes the proof of Lemma \ref{lem7}.
\end{proof}

Since $(\ti\w, \ti J(1), \ti g(1))=(\ti\w, \ti J, \ti g_J)$ is an almost K\"ahler structure on $\ti M^{2n}$ which is a complete smooth $\G$-invariant Riemannian metric on $\ti M^{2n}$ with bounded geometry, it is not hard to see that $\Omega^p_c(\ti M^{2n})\otimes_\mathds{R}\mathds{C}$ is dense in $L^2_k\Omega^p(\ti M^{2n})(1)$ (cf. T. Aubin \cite[2.6 Theorem and 2.7 Remark]{Au} or \cite{Au0}). By \eqref{cn2}, it implies that $\Omega^p_c(\ti M^{2n})\otimes_\mathds{R}\mathds{C}$ is dense in $L^2_k\Omega^p(\ti M^{2n})(t)$, $t\in [0, 1]$.
It is well known that on
$\Omega^p_c(\ti M^{2n})\otimes_\mathds{R}\mathds{C}$, $0\leq p\leq 2n$, $d^{\ti*_1}\coloneqq-\ti*_1d \ti*_1$, $\ti*_1$ being the Hodge star operator with respect to the metric $\ti g(1)$, is $L^2$-adjoint operator of $d$, and $\Delta_1\coloneqq dd^{\ti*_1}+d^{\ti*_1}d$ is $L^2$-self-adjoint operator on $\Omega^p_c(\ti M^{2n})$ (cf. \cite{Ch,Luck,Shu}).

\vskip 0.15 cm
In order to study elliptic operators on complete Riemannian manifolds, we need the following elements of linear operator theory:

\vskip 0.15 cm
Let $H$ be a Hilbert space and $T:dom(T)\to H$ be a (not necessarily bounded) linear operator defined on a dense linear subspace $dom(T)$ which is called (initial) domain. We call $T$
closed if its graph $$gr(T)\coloneqq \{(u,T(u))|u\in dom(T)\}\subset H\times H$$ is closed. We say that $S:dom(S)\to H$ is an extension of $T$ and write $T\subset S$ if $dom(T)\subset dom(S)$
and $S(u)=T(u)$ holds for all $u\in dom(T)$. We write $T=S$ if $dom(T)=dom(S)$ and $S(u)=T(u)$. We call $T$ closable if and only if $T$ has a closed extension. Since the intersection
of an arbitrary family of closed sets is closed again, a closable unbounded densely defined operator $T$ has a unique minimal closure, also called minimal closed extension, that is, a
closed operator $T_{\min}:dom(T_{\min})\to H$ which $T\subset T_{\min}$ such $T_{\min}\subset S$ holds for any closed extension $S$ of $T$. Explicitly $dom(T_{\min})$ consists of elements
$u\in H$ for which there exist a sequence $\{u_n\}_{n\geq0}$ in $dom(T)$ and an element $v$ in $H$ satisfying $\lim\limits_{n\to\infty}u_n=u$ and $\lim\limits_{n\to\infty}T(u_n)=v$. Then
$v$ is uniquely determined by this property and we put $T_{\min}(u)=v$. Equivalently $dom(T_{\min})$ is the Hilbert space completion of $dom(T)$ with respect to the inner product
$$\li u, v\ri_{gr}=\li u, v\ri_H+\li T(u), T(v)\ri_H.$$
If not stated otherwise we always use the minimal closed extension as the closed extension of a closable unbounded densely defined linear operator.

The formal adjoint of $T$ is the operator $T^*:dom(T^*)\to H$ whose domain consists of elements $v\in H$ for which there is an element $u$ in $H$ such that $\li u', u\ri_H=\li T(u'), v\ri_H$
holds for all $u'\in dom(T)$. Then $u$ is uniquely determined by this property and we put $T^*(v)=u$. Notice that $T^*$ may not have a dense domain in general. If $T$ is closable, then
$T^*_{\min}=T^*$ and $T_{\min}=(T^*)^*$. We call $T$ symmetric if $T\subset T^*$ and self-adjoint if $T=T^*$. Any self-adjoint operator is necessarily closed and symmetric. A bounded operator $T:H\to H$ is always closed and is self-adjoint if and only if it is symmetric. We call $T$ essentially self-adjoint if $T_{\min}$ is self-adjoint. The maximal closure $T_{\max}$
of $T$ is defined by the formal adjoint of $(T^*)_{\min}$. For any closure $\bar{T}$ of $T:dom(T)\to H$ we have $T_{\min}\subset\bar{T}\subset T_{\max}$. Hence if $T$ is essentially self-adjoint, then
$T_{\min}=T_{\max}$. For more details, see \cite{At,Ch,Luck,Shu}.
\vskip 0.1 cm
Let $E_1$ and $E_2$ be Hermitian vector bundles over a complete Riemannian manifold without boundary. Let $D:C^{\infty}_c(E_1)\to C^{\infty}_c(E_2)$ be an elliptic differential operator where $C^{\infty}_c(E_1)$ is the space
of smooth sections with compact support. Our main examples are Dirac operator $d+d^*$ where $$d:\Omega^p_c(M)\to\Omega^{p+1}_c(M), d^*:\Omega^p_c(M)\to\Omega^{p-1}_c(M),$$ and Hodge Laplacian operator $$\Delta=dd^*+d^*d:\Omega^p_c(M)\to\Omega^p_c(M).$$ Notice that there is a formally
adjoint operator $D^*:C^{\infty}_c(E_2)\to C^{\infty}_c(E_1)$ which is uniquely determined by the property that $\li D(u), v\ri_{L^2}=\li u, D^*(v)\ri_{L^2}$ holds for all $u,v\in C^{\infty}_c(E_1)$. It is again an elliptic
differential operator. The minimal closure $D_{\min}$ of $D:C^{\infty}_c(E_1)\to L^2C^{\infty}(E_2)$ has been defined above where $L^2C^{\infty}(E_2)$ is the Hilbert space
completion of $C^{\infty}_c(E_2)$. The maximal closure $D_{\max}$ of $D$ is defined by the formal adjoint of $(D^*)_{\min}$. Indeed, for any closure $\bar{D}$ of $D:C^{\infty}_c(E_1)\to L^2C^{\infty}(E_2)$ we have
$D_{\min}\subset\bar{D}\subset D_{\max}$. One can also describe $dom(D_{\max})$ as the space of $u\in L^2C^{\infty}(E_1)$ for which the distribution $D(u)$ actually lies in $L^2C^{\infty}(E_2)$ (cf. \cite{At,Ch,Luck,Shu}).
\vskip 0.1 cm
I. M. Singer \cite{Si} exposed a comprehensive program aimed to extend theory of elliptic operators and their index to more general situations.

N. Teleman \cite{Tele1,Tele2,Tele3} obtained results for a Hodge theory on PL manifolds or Lipschitz manifolds, and J. Cheeger \cite{Cheeger} produced a Hodge theory on a very general class of pseudo-manifolds.
\vskip 0.1 cm

Recall that by Lemma \ref{lem7} the coefficients of $\ti *_t$ are in $L^{\infty}_1(\ti M^{2n}, \ti g_J)\cap C^\infty(\overset{\circ}{\ti M} {}^{2n})$ due to the construction of $\ti g(t)$,  $\ti M {}^{2n}\backslash\overset{\circ}{\ti M} {}^{2n}$ has Hausdorff dimension $\leq2n-1$ with Lebesgue measure zero, thus
$$
d:\Omega^{p-1}_c(\ti M^{2n})\otimes_\mathds{R}\mathds{C}\longrightarrow L^2\Omega^p(\ti M^{2n})(t),
$$
$$
-\ti *_td\ti *_t:\Omega^p_c(\ti M^{2n})\otimes_\mathds{R}\mathds{C}\longrightarrow L^2\Omega^{p-1}(\ti M^{2n})(t),
$$
are well defined, where $\ti *_t=*_{\ti g(t)}$ is Hodge star operator with respect to the metric $\ti g(t)$, since $\Omega^{p-1}_c(\ti M^{2n})\otimes_\mathds{R}\mathds{C}$ and $\Omega^p_c(\ti M^{2n})\otimes_\mathds{R}\mathds{C}$ are dense in $L^2\Omega^{p-1}(\ti M^{2n})(t)$ and $L^2\Omega^p(\ti M^{2n})(t)$ respectively, $0\leq p\leq 2n$, $t\in[0, 1]$. Let $d^{\ti *_t}$ be $L^2$-adjoint operator of $d$ in the sense of distributions with respect to the metric $\ti g(t)$, $t\in[0, 1]$.
Analogous to Tan-Wang-Zhou \cite[Lemma 2.6]{TWZ2}, as done in smooth case, we have the following lemma:
\blem {\rm(}cf. M. Gromov \cite[Lemma 1.1.A]{Gr2}, Lemma \ref{lem2.7}{\rm)} \label{lem1}
Let $d^{\ti *_t}$ be $L^2$-adjoint operator of $d$ with respect to the metric $\ti g(t)$, then $d^{\ti *_t}=-\ti *_td\ti *_t, (d^{\ti *_t})^{\ti *_t}=d$ in the sense of distributions. Therefore, $d+d^{\ti *_t}$, $t\in[0,1]$, are essentially self-adjoint operators on $\ti M^{2n}$.
\elem
\begin{proof}
Notice that the Hilbert spaces
$$
L^2_k\W^p(\ti M^{2n})(t), t\in[0,1),\quad and\quad L^2_k\W^p(\ti M^{2n})(1)
$$
are quasi-isometry for $0\leq p\leq 2n$ (see Lemma \ref{lem7} and \eqref{cn2}). Hence $L^2\Omega^p(\ti M^{2n})(t)$, $0\leq p\leq 2n$, $t\in [0,1]$, can be approximated by $\Omega^p_c(\ti M^{2n})$, it is enough to choose $\a\in \Omega^{p-1}_c(\ti M^{2n})\otimes_\mathds{R}\mathds{C}$ and $\b\in \Omega^{p}_c(\ti M^{2n})\otimes_\mathds{R}\mathds{C}$. Then
$$
\a\wedge \ti *_t\b\in \Omega^{2n-1}(\overset{\circ}{\ti M} {}^{2n})\bigcap L^{\infty}(\Lambda^{2n-1}\ti M^{2n}).
$$
Since the coefficients of $\ti *_t$ are in $L^{\infty}_1(\ti M^{2n})$, $\ti M {}^{2n}\backslash\overset{\circ}{\ti M} {}^{2n}$ has Lebesgue measure zero,
$$\a\in L^2_1\Omega^{p-1}(\ti M^{2n})(1),\quad \ti *_t\b\in L^2_1\Omega^{2n-p}(\ti M^{2n})(1),$$ hence,
$$\a\wedge \ti *_t\b\in L^1\Omega^{2n-1}(\ti M^{2n})(1),$$ $$d(\a\wedge \ti*_t\b)\in L^1\Omega^{2n}(\ti M^{2n})(1).$$ By Gromov's result (cf. \cite[Lemma 1.1.A]{Gr2}, also see W. L\"uck \cite{Luck}), we have
$$
\int_{\ti M^{2n}}d(\a\wedge \ti*_t\b)=0.
$$
In fact, $\ti*_t\b\in L^2_1\Omega^{2n-p}(\ti M^{2n})(1)$ since $\ti*_t\in L^\infty(\ti M^{2n}, \ti g(1))$, $L^2_1\Omega^{2n-p}(\ti M^{2n})(t)$ and $L^2_1\Omega^{2n-p}(\ti M^{2n})(1)$ are quasi-isometry (that is, with equivalent norms, see \eqref{cn2}) by Lemma \ref{lem7}, hence there is a sequence $\ti\b_k\in\Omega^{2n-p}_c(\ti M^{2n})$ such that $$\|\ti*_t\b-\ti\b_k\|_{L^2_1(1)}\to0\quad \text{as}\quad k\to\infty.$$ Therefore, it follows from integration by parts that
$$
\int_{\ti M^{2n}}d(\a\wedge \ti*_t\b)=\lim\limits_{k\to\infty}\int_{\ti M^{2n}}d(\a\wedge \ti\b_k)=0.
$$
Thus
\begin{align}\label{co}
0=&\int_{\ti M^{2n}}d(\a\wedge \ti*_t\b) \nonumber\\
=&\int_{\ti M^{2n}}d\a\wedge \ti*_t\b+(-1)^{p-1}\int_{\ti M^{2n}}\a\wedge d(\ti*_t\b) \nonumber\\
=&\int_{\ti M^{2n}}(d\a, \b)_{\ti g(t)}dvol_{\ti g(t)}-\int_{\ti M^{2n}}(\a,  -\ti*_td\ti*_t\b)_{\ti g(t)}dvol_{\ti g(t)} \nonumber\\
=&\li d\a, \b\ri_{\ti g(t)}-\li\a,  -\ti*_td\ti*_t\b\ri_{\ti g(t)}.
\end{align}
Therefore, by \eqref{co}, formal adjoint $d^{\ti*_t}$ of $d$ is $-\ti*_td\ti*_t$, and $(d^{\ti*_t})^{\ti*_t}=d$ in Hilbert spaces. Thus, $d+d^{\ti*_t}$, $t\in[0,1]$, are essentially self adjoint operators on $\ti M^{2n}$. This completes the proof of Lemma \ref{lem1}.
\end{proof}
As in classical analysis, one can define Hodge Laplacian:
\beq \label{cp}
\Delta_t\coloneqq  dd^{\ti *_t}+d^{\ti *_t}d:L^2_2\Omega^p(\ti M^{2n})(t)\to L^2\Omega^p(\ti M^{2n})(t),
\eeq
where $0\leq p\leq 2n$ and $t\in[0, 1]$. Notice that, by Lemma \ref{lem1}, for $t\in[0, 1]$, $\Delta_t$ are essentially self-adjoint elliptic operators on $L^2\Omega^p(\ti M^{2n})(t)$ (cf. \cite{At,Ch,Luck,Shu}).

\vskip 0.1cm
In the remainder of this section, we will consider Hodge-Kodaira decomposition with respect to the (measurable) Lipschitz metrics $g(t)$, $t\in[0,1]$ (cf. \cite{Gr2,Luck,Shu,TWZ2}). Let $\a$, $d\a$, and $d^{\ti*_t}\a\in L^2\Omega^p(\ti M^{2n})(t)$,  $L^2\Omega^{p+1}(\ti M^{2n})(t)$, and $L^2\Omega^{p-1}(\ti M^{2n})(t)$ respectively, $0\leq p\leq 2n$, $t\in[0, 1]$. Then $\a$ is called ($d, d^{\ti*_t}$)-harmonic if $d\a=0$,$d^{\ti*_t}\a=0$ in the sense of distributions, $t\in[0, 1]$. Similarly, $\a\in L^2\Omega^p(\ti M^{2n})(t)$, $\a$ is called $\Delta_t$-harmonic if $\Delta_t\a=0$ in the sense of distributions, $0\leq p\leq 2n$, $t\in[0, 1]$. It is not hard to see that ($d, d^{\ti *_t}$)-harmonic implies $\Delta_t$-harmonic. By Gaffney cutoff trick, analogous to Tan-Wang-Zhou \cite[Lemma 2.7]{TWZ2}, as done in smooth case, we can obtain the following lemma (cf. M. Gromov \cite[Lemma 1.1.B]{Gr2} and W. L\"uck \cite{Luck}):
\blem [cf. Lemma \ref{lem2.7}]\label{lem2}
If an $L^2$ $p$-form $\a$ is $\Delta_t$-harmonic, then $\a$ is ($d, d^{\ti *_t}$)-harmonic in the sense of distributions, where $0\leq p\leq 2n$ and $t\in[0, 1]$.
\elem
\begin{proof}
First, we construct a family of cutoff functions $a_K$, $K\in\mathds{R}^+$, such that
$$a_K\geq 0, |\ti\nabla^1(1)a_K|_{\ti g(1)}\leq \frac{1}{K}(a_K)^{\frac{1}{2}},
|(\ti\nabla^1(1))^2a_K|_{\ti g(1)}\leq \frac{1}{K^2},$$
and the subsets $a_K^{-1}(1)$ exhaust $\ti M^{2n}$ as $K\rightarrow+\infty$ on complete noncompact manifold $\ti M^{2n}$ ($\G$-manifold), where $\ti\nabla^1(1)$ is the second canonical connection with respect to the metric $\ti g(1)$. Here, we only give the case on $\mathds{R}$ . Let
\beq \label{cq}
f(x)=\left\{
             \begin{array}{lr}
             \exp(\frac{-1}{x}), & x>0, \\
             0, & x\leq 0.
             \end{array}
\right.
\eeq
Define $$\psi(x)=\frac{f(x)}{f(x)+f(1-x)}.$$ Note that $0\leq \psi(x) \leq 1$ for $0<x<1$, $\psi(x)=0$ if $x\leq0$, $\psi(x)=1$ if $x\geq1$, and $\psi$, $\psi'$ and $\psi''$ are all bounded. Finally, for $x\geq0$, let
$$a_K=\psi^2(2-\frac{x}{K}).$$ Clearly, $a_K(x)=1$ on $[0, K]$ and $a_K(x)=0$ on $[2K, +\infty)$. For $K<x<2K$, we have
$$
a_K'=-\frac{2}{K}\psi(2-\frac{x}{K})\psi'(2-\frac{x}{K}).
$$
Since $\psi'$ is bounded, we see that $$-\frac{1}{K}C_1\sqrt{a_K}\leq a_K'\leq0$$ for some constant $C_1$. Moreover,
$$
a_K''=\frac{2}{K^2}\psi'^2(2-\frac{x}{K})+\frac{2}{K^2}\psi(2-\frac{x}{K})\psi''(2-\frac{x}{K}).
$$
Since $\psi$, $\psi'$ and $\psi''$ are bounded, we have $$|a_K''(x)|\leq \frac{C_2}{K^2}$$ for some constant $C_2$.

Suppose that $\a$ is an $L^2$ $p$-form on $\ti M^{2n}$ satisfying $$0=\Delta_t\a=(dd^{\ti*_t}+d^{\ti*_t}d)\a$$ in the sense of distributions, $0\leq p\leq 2n$, $t\in[0, 1]$. Then
\begin{align}\label{cr}
0=&\li \Delta_t\a, a_K\a\ri_{\ti g(t)} \nonumber\\
=&\int_{\ti M^{2n}}(dd^{\ti*_t}+d^{\ti*_t}d)\a\wedge a_K\ti *_t\a \nonumber\\
=&\int_{\ti M^{2n}}(d^{\ti *_t}\a, d^{\ti*_t}(a_K\a))_{\ti g(t)}\frac{\ti\w^n}{n!}+\int_{\ti M^{2n}}(d\a,  d(a_K\a))_{\ti g(t)}\frac{\ti\w^n}{n!} \nonumber\\
=&a_K\int_{\ti M^{2n}}(d^{\ti*_t}\a, d^{\ti*_t}\a)_{\ti g(t)}\frac{\ti\w^n}{n!}+\int_{\ti M^{2n}}d^{\ti*_t}\a\wedge(-\ti*_t^2(da_K\wedge\ti*_t\a)) \nonumber\\
&+a_K\int_{\ti M^{2n}}(d\a,  d\a))_{\ti g(t)}\frac{\ti\w^n}{n!}+\int_{\ti M^{2n}}d\a\wedge\ti*_t(da_K\wedge\a) \nonumber\\
=&I_1(K)+I_2(K),
\end{align}
where $$I_1(K)=a_K\{\li d\a, d\a\ri_{\ti g(t)}+\li d^{\ti*_t}\a, d^{\ti*_t}\a\ri_{\ti g(t)}\},$$
\beq \label{cs}
I_2(K)=\int_{\ti M^{2n}}d\a\wedge\ti*_t(da_K\wedge\a)+\int_{\ti M^{2n}}d^{\ti*_t}\a\wedge(-\ti*_t^2(da_K\wedge\ti*_t\a))
\eeq
and
\beq \label{ct}
I_1(K)+I_2(K)=0.
\eeq
Hence,
\begin{align*}
|I_2(K)|\leq &C_0\int_{\ti M^{2n}}|da_K|_{\ti g(t)}(|d\a|_{\ti g(t)}+|d^{\ti*_t}\a|_{\ti g(t)})|\a|_{\ti g(t)}\frac{\ti\w^n}{n!} \\
\leq &C_0\int_{\ti M^{2n}}|da_K|_{\ti g(t)}\{(d\a, d\a)_{\ti g(t)}+(d^{\ti*_t}\a, d^{\ti*_t}\a)_{\ti g(t)}+(\a, \a)_{\ti g(t)}\}\frac{\ti\w^n}{n!},
\end{align*}
where $$|da_K|^2_{\ti g(t)}=(da_K, da_K)_{\ti g(t)}, |d\a|^2_{\ti g(t)}=(d\a, d\a)_{\ti g(t)},$$ and
$$|d^{\ti*_t}\a|^2_{\ti g(t)}=(d^{\ti *_t}\a, d^{\ti *_t}\a)_{\ti g(t)}.$$
Since $$a_K\geq0, |\ti\nabla^1(1)a_K|_{\ti g(1)}\leq \frac{1}{K}(a_K)^{\frac{1}{2}},$$ and the subsets $a_K^{-1}(1)\subset\ti M^{2n}$ exhaust $\ti M^{2n}$ as $K\to +\infty$. Then
$$|I_2(K)|\to0\text{ as }K\to +\infty.$$ Since $|I_2(K)|=I_1(K)$, hence $$I_1(K)\to0\text{ as }K\to +\infty.$$ On the other hand,
$$
I_1(K)=\int_{\ti M^{2n}}a_K(|d\a|^2_{\ti g(t)}+|d^{\ti*_t}\a|^2_{\ti g(t)})\frac{\ti\w^n}{n!}.
$$
As $K\to +\infty$,
$$I_1(\infty)=\int_{\ti M^{2n}}(|d\a|^2_{\ti g(t)}+|d^{\ti*_t}\a|^2_{\ti g(t)})\frac{\ti\w^n}{n!}.$$
Therefore, $\li d\a, d\a\ri_{\ti g(t)}=0$, that is,
$d\a=0=d^{\ti*_t}\a$ in the sense of distributions. We complete the proof of Lemma \ref{lem2}.
\end{proof}
Let $L^2\Omega^p(\ti M^{2n})(t)$ be $L^2$-space of $p$-forms on $\ti M^{2n}$ with respect to the metric $\ti g(t)$, $0\leq p\leq 2n$, $t\in[0, 1]$. By Lemma \ref{lem2}, as done in smooth case, we have Hodge-Kodaira
decomposition (reduced)
\begin{align} \label{cu}
L^2\Omega^p(\ti M^{2n})(t)=&\mathcal{H}^p_{(2)}(\ti M^{2n})(t)\nonumber\\
&\oplus\overline{dL^2\Omega^{p-1}(\ti M^{2n})(t)}\oplus\overline{d^{\ti*_t}L^2\Omega^{p+1}(\ti M^{2n})(t)},
\end{align}
$0\leq p\leq 2n$, $t\in[0, 1]$, where $$\mathcal{H}^p_{(2)}(\ti M^{2n})(t)=\ker\Delta_t|_{L^2\Omega^p(\ti M^{2n})(t)} ,$$
$\overline{dL^2\Omega^{p-1}(\ti M^{2n})(t)}$ and $\overline{d^{\ti*_t}L^2\Omega^{p+1}(\ti M^{2n})(t)}$ are
closure of $d\Omega^{p-1}_c(\ti M^{2n})(t)$ and $d^{\ti*_t}\Omega^{p+1}_c(\ti M^{2n})(t)$ with respect to $L^2(t)$-norm respectively.

For the Hodge theory on non-compact smooth manifolds, see J. Dodziuk \cite{Dod}, M. Gromov \cite{Gr2} and N. Hitchin \cite{Hit}. See also N. Teleman \cite{Tele2,Tele3} for the Hodge theory on closed
PL manifolds or Lipschitz Riemannian manifolds.

Since $\ti g(1)=\ti g_J$ on $\ti M^{2n}$ is a smooth $\G$-invariant complete metric with bounded geometry, it is clear that $$\mathcal{H}^p_{(2)}(\ti M^{2n})(1)\subset \Omega^p(\ti M^{2n})\cap L^2\Omega^p(\ti M^{2n})(1).$$
For any $ t\in [0, 1)$, $\ti g(t)$ is smooth on
$\overset{\circ}{\ti M} {}^{2n}$, $$\ti M {}^{2n}\backslash\overset{\circ}{\ti M} {}^{2n}=\bigcup\limits_{\g\in\G}\g(\partial F)$$ is piecewise smooth of Hausdorff dimension $\leq2n-1$ with Lebesgue measure zero and $C^0(\bar{F})$ is dense in $L^2(\bar{F})$ (cf. T. Aubin \cite[2.9 Theorem]{Au} or \cite{Au0}). Notice that for $t\in[0, 1)$, $-\ti*_td\ti*_t$ involves differentiating the $\ti*_t$ operator and in general, for (measurable) Lipschitz almost K\"ahler structure $(\ti\w, \ti J(t), \ti g(t))$,  the coefficients of $\ti *_t$ are at best bounded measurable, with no control of its regularity. But, in our situation, underlying manifold $(\ti M^{2n}, \ti g(1))=(\ti M^{2n}, \ti g_J)$ is a complete $\G$-invariant Riemannian manifold with bounded geometry,
$L^2_k\Omega^p(\ti M^{2n})(t)$ can be approximated by $\Omega^p_c(\ti M^{2n})$ with respect to the Riemannian metric $\ti g(1)=\ti g_J$, and $L^2_k\Omega^p(\ti M^{2n})(t)$ and $L^2_k\Omega^p(\ti M^{2n})(1)$ are quasi-isometry (See Lemma \ref{lem7}). Hence, we can prove that
$$\mathcal{H}^p_{(2)}(\ti M^{2n})(t)\subset \Omega^p(\ti M {}^{2n})\cap L^2\Omega^p(\ti M^{2n})(1),  t\in [0, 1).$$
Define
\beq \label{cv}
b^p_{(2)}(\ti M^{2n})(t)\coloneqq\dim_{\G}\mathcal{H}^p_{(2)}(\ti M^{2n})(t), 0\leq p\leq 2n, t\in[0, 1].
\eeq
Notice that on $\ti M^{2n}$, $t\in[0,1]$, we defined $L^2_k\Omega^p(\ti M^{2n})(t)$ for a nonnegative integer $k$ be the $k$-th Hilbert space of $p$-forms on $\ti M^{2n}$ with respect to the
$\G$-invariant (measurable) Lipschitz metric $\ti g(t)$ (see \eqref{cn2}), that is,
the Hilbert square completion of $\Omega^p_c(\ti M^{2n})$ with respect to the inner product or norm induced by the $\G$-invariant (measurable) Lipschitz metric $\ti g(t)$ which is equivalent
to the following definition:
\begin{align*}
\li \a, \b\ri_{L^2_k(t)}\coloneqq&\li (1+\Delta_t)^{\frac{k}{2}}\a, (1+\Delta_t)^{\frac{k}{2}}\b\ri_{L^2(t)}\\
=&\li \a, (1+\Delta_t)^{k}\b\ri_{L^2(t)},
\end{align*}
So,
\beq \label{cw}
\|\a\|_{L^2_k(t)}\coloneqq\sqrt{\li \a, (1+\Delta_t)^{k}\a\ri_{L^2(t)}}.
\eeq
Hence, one can define a Hilbert $\G$-cochain complex $L^2_{l-*}\Omega^*(\ti M^{2n})(t)$ as follows (cf. \cite{At,BL,Luck,Shu}):
\begin{align} \label{cx}
0\to L^2_{l}\Omega^0(\ti M^{2n})(t)&\xrightarrow{d}L^2_{l-1}\Omega^1(\ti M^{2n})(t)\xrightarrow{d}\cdots\nonumber\\
&\xrightarrow{d}L^2_{l-2n}\Omega^{2n}(\ti M^{2n})(t)\to 0,\quad l\geq2n.
\end{align}

By \eqref{cx}, define the reduced $p$-th $L^2$-cohomology as follows:
 $$\mathrm{H}^p_{(2)}(\ti M^{2n})(t)\coloneqq\left. \ker d|_{L^2\Omega^p(\ti M^{2n})(t)}\middle/\overline{dL^2\Omega^{p-1}(\ti M^{2n})(t)}\right..$$ It is easy to see that
\beq \label{cy}
\mathrm{H}^p_{(2)}(\ti M^{2n})(t)\cong\mathcal{H}^p_{(2)}(\ti M^{2n})(t), 0\leq p\leq 2n, t\in[0, 1].
\eeq
For more details, see L\"uck \cite[Chapter 1]{Luck}.

In summary, analogous to Tan-Wang-Zhou \cite[Proposition 2.8]{TWZ2}, since $\ti g(t)$, $t\in[0, 1]$, is a family of $\G$-invariant (measurable) Lipschitz metrics on $\ti M^{2n}$ which are with equivalent norms (with Lipschitz condition \eqref{cn}), by the definition and properties of the $\G$-dimension we have the following theorem:
\bthm [cf. Theorem \ref{thm2.9}]\label{bi}
Suppose that $(M^{2n}, \w, J, g_J)$ be a closed almost K\"ahler manifold with an infinite fundamental group $\G=\pi_1(M^{2n})$. Let $(\ti M^{2n}, \ti\w, \ti J, \ti g_J)$ be
the universal covering of $(M^{2n}, \w, J, g_J)$. Suppose that $F$ is a special fundamental domain such that
$$\ti M {}^{2n} =\bigcup\limits_{\g\in \G} \g  \bar{F},
\overset{\circ}{\ti M} {}^{2n} =\bigcup\limits_{\g\in \G} \g  F.$$
On ${\ti M}^{2n}$ we can construct a family of $\G$-invariant (measurable) Lipschitz almost K\"ahler structures $(\ti\w, \ti J(t), \ti g(t))$, $t\in[0, 1]$. When $t=1$, $(\ti\w, \ti J(1), \ti g(1))$$=$
$(\ti\w, \ti J, \ti g_J)$ which is smooth as a $\G$-invariant almost K\"ahler structure; $t=0$, $(\ti\w, \ti J(0), \ti g(0))$ is a $\G$-invariant (measurable) Lipschitz K\"ahler flat structure. Thus  $(\ti\w, \ti J(t), \ti g(t))$, $t\in[0, 1]$, is a family of $\G$-invariant (measurable) Lipschitz almost K\"ahler structures on $\ti M^{2n}$ which are with equivalent norms (with Lipschitz condition \eqref{cn}). For the family of $\G$-invariant (measurable) Lipschitz metrics
$\ti g(t)$, $t\in[0, 1]$, we have the following Hodge-Kodaira decomposition
$$
L^2\Omega^p(\ti M^{2n})(t)=\mathcal{H}^p_{(2)}(\ti M^{2n})(t)\oplus\overline{dL^2\Omega^{p-1}(\ti M^{2n})(t)}\oplus\overline{d^{*_t}L^2\Omega^{p+1}(\ti M^{2n})(t)},
$$
$\mathcal{H}^p_{(2)}(\ti M^{2n})(t)$$\cong$$\mathrm{H}^p_{(2)}(\ti M^{2n})(t)$, $\mathcal{H}^p_{(2)}(\ti M^{2n})(t)\subset\Omega^p(\ti M^{2n})\cap L^2\Omega^p(\ti M^{2n})(1)$,
$b^p_{(2)}(\ti M^{2n})(t)=b^p_{(2)}(\ti M^{2n})(1)<\infty$, and $$b^p_{(2)}(\ti M^{2n})(t)=b^{2n-p}_{(2)}(\ti M^{2n})(t), 0\leq p\leq 2n, t\in[0, 1],$$ where
$$b^p_{(2)}(\ti M^{2n})(t)\coloneqq\dim_{\G}\mathcal{H}^p_{(2)}(\ti M^{2n})(t).$$ In particular,
\begin{align*}
\chi_{(2)}(\ti M^{2n})(t)&=\chi_{(2)}(\ti M^{2n})(1),  t\in[0,1)\\
&=\chi(M^{2n}),
\end{align*}
where $$\chi_{(2)}(\ti M^{2n})(t)\coloneqq\sum_{p=0}^{2n}(-1)^pb_{(2)}^p(\ti M^{2n})(t).$$
\ethm
\begin{proof}
Since $\Delta_t\coloneqq dd^{\ti*_t}+d^{\ti*_t}d$ and $d+d^{\ti*_t}$, $t\in[0,1]$, are essentially self-adjoint elliptic operators on $\ti M^{2n}$. By Lemmas \ref{lem1} and \ref{lem2},
$$\mathcal{H}^p_{(2)}(\ti M^{2n})(t)=\ker\Delta_t|_{L^2\W^p(\ti M^{2n})(t)},$$ $$\ker d|_{L^2\W^p(\ti M^{2n})(t)}=d^{*_t}(\W^{p+1}_c(\ti M^{2n}))^\perp,$$ $$\overline{dL^2\Omega^{p-1}(\ti M^{2n})(t)}=\overline{d\W^{p-1}_c(\ti M^{2n})},$$ with respect to the metric $\ti g(t)$, for $t\in[0,1]$, $0\leq p\leq 2n$. Hence
$$\mathcal{H}^p_{(2)}(\ti M^{2n})(t)\cong H^p_{(2)}(\ti M^{2n})(t), \quad t\in[0,1],0\leq p\leq 2n.$$ By Lemma \ref{lem7},
$\ti g(t)$, $t\in[0, 1)$ and $\ti g(1)$ are quasi-isometry (that is, with equivalent norms), hence  $$H^p_{(2)}(\ti M^{2n})(t)\cong H^p_{(2)}(\ti M^{2n})(1).$$ Since  $\ti g(t)$, $t\in[0, 1]$ are $\G$-invariant (measurable) Lipschitz metrics on $\ti M^{2n}$ and
$\ti g(1)=\ti g_J$ is a $\G$-invariant smooth metric on $\ti M^{2n}$, thus
$$b^p_{(2)}(\ti M^{2n})(t)=\dim_{\G} H^p_{(2)}(\ti M^{2n})(t)=\dim_{\G}H^p_{(2)}(\ti M^{2n})(1)=b^p_{(2)}(\ti M^{2n})(1),$$ for $t\in[0,1]$, $0\leq p\leq 2n$.

Since $\ti g(1)=\ti g_J$ is a $\G$-invariant smooth metric on $\ti M^{2n}$, this allows us to define the analytic $p$-th $L^2$-Betti number by  the von Neumann dimension of the
$\G$-module $\mathcal{H}^p_{(2)}(\ti M^{2n})(1)$. This can be interpreted in terms of the heat kernel by the following expression (cf. Atiyah \cite{At} and L\"uck \cite{Luck})
$$
b^p_{(2)}(\ti M^{2n})(1)=\lim_{s\to+\infty}\int_F {\rm tr\/}_\C(e^{-s\Delta^p_1}(x,x))\frac{\ti \w^n}{n!},
$$
where $\Delta^p_1=\Delta_1|_{L^2\W^p(\ti M^{2n})(1)}$ and $F$ is the fundamental domain for the $\G$-action. Since $M^{2n}=\ti M^{2n}/ \G$ is a closed Riemannian manifold, $$vol(M^{2n}, g_J)<\infty.$$ Hence $b^p_{(2)}(\ti M^{2n})(1)<\infty$ (cf. W. L\"uck \cite{Luck} or Cheeger-Gromov \cite{CG2,CG3}). In fact, the coefficients of $\Delta_1$ are bounded smooth $\G$-invariant on $\ti M^{2n}$.

Similarly, we have $$\mathrm{H}^p_{(2)}(\ti M^{2n})(t)\cong\mathrm{H}^{2n-p}_{(2)}(\ti M^{2n})(t),\quad t\in[0,1],0\leq p\leq 2n.$$ Notice that $(\ti M^{2n},\ti g(1))=(\ti M^{2n},\ti g_J)$
is a $\G$-invariant Riemannian manifold, $(M^{2n}, g_J)=(\ti M^{2n}, \ti g_J)/ \G$ is a closed Riemannian manifold of dimension $2n$ which is orientable. Then
$$\mathcal{H}^p_{(2)}(\ti M^{2n})(1)\cong\mathrm{H}^p_{(2)}(\ti K),
$$
where $\ti K$ is a $\G$-equivariant smooth triangulation of $\ti M^{2n}$ which consists of a simplicial complex $\ti K$ with simplicial $\G$-action. By basic properties of cellular $L^2$-Betti numbers, $b^p_{(2)}(\ti M^{2n})(1)=b^{2n-p}_{(2)}(\ti M^{2n})(1)$ (cf. Dodziuk \cite{Dod}, L\"uck \cite{Luck}),
$$b^p_{(2)}(\ti M^{2n})(t)=b^{2n-p}_{(2)}(\ti M^{2n})(t)\quad \text{for}\; t\in[0,1],0\leq p\leq 2n.$$

For any $$\a\in\mathcal{H}^p_{(2)}(\ti M^{2n})(t)=\ker\Delta_t|_{L^2\Omega^p{\ti M^{2n}}(t)},$$ by the definition of the norm $\|\cdot\|_{L^2_k(t)}$ (see \eqref{cw} and Definition 3.4), it is easy to get
$$\|\a\|^2_{L^2_k(t)}=\li \a, (1+\Delta_t)^{k}\a\ri_{L^2(t)}=\|\a\|^2_{L^2(t)}.$$
By Lemma \ref{lem7}, $L^2_k\Omega^p(\ti M^{2n})(t)$ and $L^2_k\Omega^p(\ti M^{2n})(1)$ are quasi-isometry (that is, with equivalent norms), hence $\a\in L^2_k\Omega^p(\ti M^{2n})(1)$, for any $k\in\N$. By Sobolev embedding theorem, $$\a\in C^l(\Lambda^pT^*\ti M^{2n})\cap L^2\Omega^p(\ti M^{2n})(1),$$ for any $l\in\N$ such that $k-l>n$, hence $$\a\in\Omega^p(\ti M^{2n})\cap L^2\Omega^p(\ti M^{2n})(1).$$ Therefore, $$\mathcal{H}^p_{(2)}(\ti M^{2n})(t)\subset\Omega^p(\ti M^{2n})\cap L^2\Omega^p(\ti M^{2n})(1).$$

Since $\ti g(t)$, $t\in[0,1]$, is a family of $\G$-invariant (measurable) Lipschitz metircs on $\ti M^{2n}$, and $L^2_k(t)$, $t\in[0,1)$ and $L^2_k(1)$ are quasi-isometry where $k$ is nonnegative integer. Then
$$
\chi_{(2)}(\ti M^{2n})(t)=\chi_{(2)}(\ti M^{2n})(1), t\in[0,1)
$$
(cf. Atiyah \cite{At}, L\"uck \cite{Luck} and Shubin \cite{Shu}).
\end{proof}
\brem  Notice that $\ti g(t)$, $t\in[0,1)$, are $\G$-invariant smooth (resp. (measurable) Lipschitz) metrics on $\overset{\circ}{\ti M} {}^{2n}$ (resp. $\ti M^{2n}$). So by Lemma \ref{lem7}, $\Delta_t$, $t\in[0,1)$, are elliptic operators with bounded smooth (resp. $L^\infty_k$, $k\in\N$) coefficients on $\overset{\circ}{\ti M} {}^{2n}$ (resp. $\ti M^{2n}$). But, as in closed symplectic manifolds, $\mathcal{H}^p_{(2)}(\ti M^{2n})(t)$, $t\in[0,1]$, $0\leq p\leq 2n$, are smooth harmonic forms on $\ti M^{2n}$ with respect to the $\G$-invariant (measurable) Lipschitz metrics $\ti g(t)$ (cf. Section 2).
\erem
\vskip 0.1cm

\section{Vanishing of $L^2$-Betti numbers of symplectic parabolic manifolds}
For non-elliptic K\"ahler manifolds, we have a vanishing theorem due to M. Gromov (\cite{Gr2}, K\"ahler hyperbolic case) and N. Hitchin (\cite{Hit}, K\"ahler parabolic case).
Notice that the proof of the vanishing theorem is based on the identity $[L_{\ti\w}, \Delta_{\ti g}]=0$ which implies the hard Lefschetz condition of K\"ahler manifolds. In Tan, Wang and Zhou \cite{TWZ2}, for a non-elliptic
almost K\"ahler manifold $(M^{2n}, \w, J, g_J)$ with hard Lefschetz condition, then we have a vanishing theorem for the universal cover $(\ti M^{2n}, \ti\w, \ti J, g_{\ti J})$ of $(M^{2n}, \w, J, g_J)$.
But, in general, non-elliptic complete almost K\"ahler manifold $(M^{2n}, \w, J, g_J)$ has no vanishing theorem. R. Hind and A. Tomassini \cite{HT} have the following result:
\beg
\emph{
By using methods of contact geometry, starting with a contact manifold $M$ having an exact symplectic filling (see Hind-Tomassini \cite[Definition 3.1, Example 3.11]{HT}), R. Hind and A. Tomassini constructed $d$(bounded) complete almost K\"ahler manifold $Y$
satisfying $\mathcal{H}^1_{(2)}(Y)\neq \{0\}.$ For this case, the Lefschetz map has no hard Lefschetz property.}
\eeg
Notice that the complete almost K\"ahler manifold $Y$ in the above example does not admit a cocompact free proper action of a discrete group by isometries. Thus Example 4.1 can not be used to construct counterexample to the
Chern-Hopf conjecture in symplectic version since we only consider the universal covering spaces of closed symplectic manifolds with infinite fundamental group. On the other hand, in general, if ($M^{2n}, \w, J, g_J$) is a complete almost K\"ahler manifold of dimension $2n$, then for any $d$-closed $L^2$ $(n-p)$-form $\a$, $0<p\leq n$,
$\w^p\wedge\a$ is an exact $L^2$ $(n+p)$-form on $M^{2n}$ (cf. N. Hitchin \cite[Theorem 1]{Hit}). But the Lefschetz maps
$$L^p_{\w}:H^{n-p}_{(2)}(M^{2n}, g_J)\to H^{n+p}_{(2)}(M^{2n}, g_J),\quad 0<p\leq n,$$
may not be injective and surjective unless it has the hard Lefschetz property. Note that the hard Lefschetz condition is a sufficient condition such that the Lefschetz maps $$L^p_{\w}:H^{n-p}_{(2)}(M^{2n}, g_J)\to H^{n+p}_{(2)}(M^{2n}, g_J),\quad 0<p\leq n,$$ are injective and surjective. Hence, for a $d$(bounded) complete almost K\"ahler manifold without the hard Lefschetz property, we may not be able to obtain directly a vanishing theorem for
the reduced $L^2$-cohomology.
\vskip 0.1cm
In the remainder of the paper, we only consider complete manifolds with a cocompact free proper action of a discrete group $\G$ by isometries.

Suppose that a $2n$-dimensional closed almost K\"ahler manifold $(M^{2n}, \w, J, \linebreak g_J)$ is symplectic parabolic. Let $$\pi:(\ti M^{2n}, \ti\w, \ti J, \ti g_J)\to(M^{2n}, \w, J, g_J)$$
be the universal covering of $(M^{2n}, \w, J, g_J)$. Then $$ \ti\w=\pi^*\w,\quad\ti J=\pi^*J,\quad\ti g_J=\pi^*g_J$$ are $\G=\pi_1(M^{2n})$ invariant, and $\ti\w=d\b$, where fixed a point
$$x_0\in\ti M^{2n},|\b(x)|_{\ti g_J}\leq C(\ti\rho_1(x_0, x)+1),$$ $C$ is a positive constant, $\ti\rho_1$ is a distance function defined by $\ti g_J$. Where $\ti\rho_1$ is the
infimum of the lengths of all piecewise smooth curves starting from $x_0$ and ending at $x$ with respect to Riemannian metric $\ti g_J$ (cf. \cite{Cha,DK}).

In order to prove Theorem 1.7 without the hard Lefschetz condition, as done in K\"ahler case (cf. Theorem \ref{cc}), we need making a suitable choice of the $\G$-invariant (measurable) Lipschitz K\"ahler structure $(\ti\w, \ti J(0),\ti g(0))$ which is $\G$-homotopy equivalent to $(\ti\w, \ti J, \ti g_J)$, and restricted to $\overset{\circ}{\ti M} {}^{2n}$, $\ti g(1)$ and $\ti g(0)$ are $\G$ quasi isometric by a global deformation of $\ti\w$-compatible almost complex structures on $\ti M^{2n}$ off the Lebesgue measure zero subset $\ti M {}^{2n}\backslash\overset{\circ}{\ti M} {}^{2n}$
(see Section 3) such that the Lefschetz maps
$$L^p_{\ti \w}:H^{n-p}_{(2)}(\ti M^{2n}, \ti g(0))\to H^{n+p}_{(2)}(\ti M^{2n}, \ti g(0)),\quad 0<p\leq n,$$ are injective and surjective.

One can define the Lefschetz map (cf. \cite{GH,Yan})
$$L_{\ti\w}:\a\in\Omega^p(\ti M^{2n})\to\ti\w\wedge\a\in\Omega^{p+2}(\ti M^{2n}).$$
Since coefficients of $\ti *_0$ and $\ti g(0)$ are in $C^\infty(\overset{\circ}{\ti M} {}^{2n})\cap L^\infty_k(\ti M^{2n},\ti g_J)$ for $k$ nonnegative integer, it is not hard to see that on $\overset{\circ}{\ti M} {}^{2n}$
\beq  \label{di}
[L_{\ti\w}, \Delta_0]=0,
\eeq
where $$\Delta_0=dd^{\ti*_0}+d^{\ti*_0}d:L^2\Omega^p(\ti M^{2n})(0)\to L^2\Omega^p(\ti M^{2n})(0), 0\leq p\leq 2n,$$ in the distributional sense. By(\ref{di}), we have the following
theorem:
\bthm \label{bj}{\rm(}Lefschetz, see M. Gromov \cite[1.2.A Theorem]{Gr2}{\rm)}
Suppose that $(\ti M^{2n}, \ti\w, \ti J(0), \ti g(0))$ is a $\G$-(measurable) Lipschitz K\"ahler flat manifold. Restricted to $\overset{\circ}{\ti M} {}^{2n}$, $(\ti\w, \ti J(0), \ti g(0))$
is a K\"ahler flat structure which is $\G$ quasi isometric to $\ti g(1)$. Hence,
$$L^k_{\ti\w}:\mathcal{H}^{n-k}_{(2)}(\ti M^{2n})(0)\to\mathcal{H}^{n+k}_{(2)}(\ti M^{2n})(0),\quad 0\leq k\leq n,$$
is an isomorphism.
\ethm
\begin{proof}
Notice that
$$L^k_{\ti\w}:\Omega^{n-k}(\ti M^{2n})\to\Omega^{n+k}(\ti M^{2n}),\quad 0\leq k\leq n,$$
is an isomorphism (see Yan \cite[p.148]{Yan}) and $(\ti\w, \ti J(0), \ti g(0))$ is a K\"ahler flat structure on $\overset{\circ}{\ti M} {}^{2n}$. Hence, by Lemma \ref{lem7}, $[L^k_{\ti\w}, \Delta_0]=0$ on $\overset{\circ}{\ti M} {}^{2n}$, where $$\Delta_0=dd^{\ti*_0}+d^{\ti*_0}d,d^{\ti*_0}=-*_{\ti g(0)}d*_{\ti g(0)},$$ and $\ti*_0=*_{\ti g(0)}$ is the Hodge star operator with respect to the metric
$\ti g(0)$ on $\overset{\circ}{\ti M} {}^{2n}$. Therefore, $[L^k_{\ti\w}, \Delta_0]=0$, $0\leq k\leq n$, on $\ti M^{2n}$ in the sense of distributions, since $$\text{coefficients of }\ti g(0)\text{ are in } C^\infty(\overset{\circ}{\ti M} {}^{2n})\cap L_k^\infty(\ti M^{2n}, \ti g_J)$$ for $k$ being a nonnegative integer and $\ti M {}^{2n}\backslash\overset{\circ}{\ti M} {}^{2n}$ has Lebesgue measure zero.
Recall that by Theorem \ref{bi}, $$\mathcal{H}^{p}_{(2)}(\ti M^{2n})(0)\subset\Omega^p(\ti M^{2n})\cap L^2\Omega^p(\ti M^{2n})(1)$$ and $\overset{\circ}{\ti M} {}^{2n}\subset\ti M^{2n}$ is an open and dense submanifold of $\ti M^{2n}$.
Thus, it immediately yields the ``injective" part of the theorem.

To prove the surjectivity we invoke the adjoint operator with respect to the metric $\ti g(0)$
$$\Lambda_k\coloneqq(L^k_{\ti\w})^{\ti*_0}:\Omega^{n+k}(\overset{\circ}{\ti M} {}^{2n})\to\Omega^{n-k}(\overset{\circ}{\ti M} {}^{2n}).
$$
Since the Lefschetz map $L^k_{\ti\w}$ is parallel on $\overset{\circ}{\ti M} {}^{2n}$ with respect to $\ti g(0)$, $\Lambda_k$ is also parallel and
$$\Lambda_k:\mathcal{H}^{n+k}_{(2)}(\ti M^{2n})(0) \to \mathcal{H}^{n-k}_{(2)}(\ti M^{2n})(0).$$
Hence $\Lambda_k$ is also injective and therefore $L^k_{\ti\w}$ has dense image. On the other hand, the Lefschetz map $L^k_{\ti\w}$ is a quasi-isometry, i.e. for $\varphi\in\W^{n-k}(\ti M^{2n})$,
$$C^{-1}\|\varphi\|_{L^2(\ti M^{2n})(1)}\leq\|L^k_{\ti\w}(\varphi)\|_{L^2(\ti M^{2n})(1)}\leq C\|\varphi\|_{L^2(\ti M^{2n})(1)}.$$
The injectivity of $L^k_{\ti\w}$  and dense image property (that is, $L^k_{\ti\w}(\mathcal{H}^{n-k}_{(2)}(\ti M^{2n})(0))$ is dense in $\mathcal{H}^{n+k}_{(2)}(\ti M^{2n})(0)$) show that
$$L^k_{\ti\w}:\mathcal{H}^{n-k}_{(2)}(\ti M^{2n})(0)\to\mathcal{H}^{n+k}_{(2)}(\ti M^{2n})(0)$$ is bijective.

This completes the proof of Theorem \ref{bj}.
\end{proof}

\brem Theorem \ref{bj} shows that the $\G$-invariant (measurable) Lipschitz K\"ahler flat manifold $(\ti M^{2n},\ti\w$, $\ti J(0),\ti g(0))$ has the hard Lefschetz property (cf. Griffith-Harris\cite{GH}).
But, in general, the Lefschetz map $$L^k_{\ti\w}:\mathcal{H}^{n-k}_{(2)}(\ti M^{2n})\to\mathcal{H}^{n+k}_{(2)}(\ti M^{2n}),\quad 0\leq k\leq n,$$ is not an isomorphism. Since the $\G$-invariant (measurable) Lipschitz K\"ahler structure is of the hard Lefschetz property, it is natural to ask the following question:

For any closed aspherical manifold of even dimension, does there exist a $\G$-invariant (measurable) Lipschitz K\"ahler structure? This is equivalent to the existence of a $\G$-invariant (measurable) Lipschitz symplectic form on its universal
covering space (cf. Section 5 for the nonpositive curvature case).
\erem
As done in K\"ahler case, analogous to Tan-Wang-Zhou \cite[Theorem 1.3]{TWZ2} (cf. N. Hitchin \cite[Theorem 1,2]{Hit}, or M. Gromov \cite{Gr2}), we have the following vanishing theorem:
\bthm \label{bk}
Suppose that $(M^{2n}, \w, J, g_J)$ is a closed symplectic parabolic manifold of dimension $2n$ and $(\ti M^{2n}, \ti\w, \ti J, \ti g_J)$ is its universal covering.
Let $(\ti M^{2n}, \ti\w, \ti J(t)$, $\ti g(t))$ be a family of $\G$-invariant (measurable) Lipschitz almost K\"ahler structures which are with equivalent norms (with Lipschitz condition \eqref{cn}) by the deformation of $\ti\w$-compatible almost complex structure on $\ti M^{2n}$ off the Lebesgue measure zero subset $\ti M {}^{2n}\backslash\overset{\circ}{\ti M} {}^{2n}$ constructed in Section 3, where $$(\ti M^{2n}, \ti\w, \ti J(1), \ti g(1))=(\ti M^{2n}, \ti\w, \ti J, \ti g_J)$$
is a $\G$-invariant smooth almost K\"ahler manifold and  $(\ti M^{2n}, \ti\w, \ti J(0), \ti g(0))$ is a $\G$-invariant (measurable) Lipschitz K\"ahler flat manifold. Then for $p\not= n$
$$b^p_{(2)}(\ti M^{2n})(1)=b^{p}_{(2)}(\ti M^{2n})(0)=0.$$
\ethm
\begin{proof}
Since $(M^{2n}, \w)$ is symplectic parabolic in the sense of Definition 1.1, then there is a 1-form $\b$ on $\ti M^{2n}$ such that $\ti \w=d\b$ and for a fixed $x_0\in\ti M^{2n}$,
$$|\b(x)|_{\ti g(1)}\leq C(\ti\rho_1(x_0, x)+1),$$
where $\ti\rho_1$ is a distance function defined by $\ti g(1)=\ti g_J$ and $C$ is a positive constant which depends on $x_0$, $M^{2n}$, $\w$, and $J$.

Let $B_K$ be the ball in $\ti M^{2n}$ with center $x_0$ and radius $K$. Take a smooth function $a_K:\ti M^{2n}\to \mathds{R}^+$ (defined in Section 3) with
\begin{align*}
&0\leq a_K(x)\leq 1,\\
&a_K(x)=1, \quad for\, x\in B_K,\\
&a_K(x)=0, \quad for\, x\in \ti M^{2n}\backslash B_{2K},\\
&|da_K|_{\ti g(1)}\leq \frac{\sqrt{a_K}}{K}.
\end{align*}
Recall that $$\ti \w=d\b\quad \text{and}\quad |\b(x)|_{\ti g(1)}\leq C(\ti\rho_1(x_0, x)+1).$$ For every $\eta\in\mathcal{H}^{p}_{(2)}(\ti M^{2n})(0)$, $0\leq p\leq n-1$, then, by Theorem \ref{bi},
$$\eta\in \Omega^p(\ti M^{2n})\cap L^2\Omega^p(\ti M^{2n})(1).$$
Notice that $L^2_k\Omega^p(\ti M^{2n})(1)$ is equivalent to $L^2_k\Omega^p(\ti M^{2n})(t)$ (that is, they are quasi-isometric or with equivalent norms), $\forall t\in[0, 1]$, $0\leq p\leq 2n$.
$d(a_K\b\wedge\eta)$ has compact support. So $d(a_K\b\wedge\eta)\in L^2\Omega^{p+2}(\ti M^{2n})(1)$. We want to show that as $K\to\infty$, these forms converge in $L^2$ to $\ti \w\wedge\eta.$
Consider
\beq \label{dj}
d(a_K\b\wedge\eta)=da_K\wedge\b\wedge\eta+a_K\ti\w\wedge\eta.
\eeq
As $a_K=1$ on $B_K$, $a_K\ti\w\wedge\eta$ converges pointwise to $\ti\w\wedge\eta$. Moreover, $|\ti\w|_{\ti g(1)}=Constant$ since $\nabla^1(1)\ti\w=0$, where $\nabla^1(1)$
is the second canonical connection with respect to the $\G$-invariant metric $\ti g(1)$ and $\eta\in L^2\Omega^p(\ti M^{2n})(1)$, so  $\ti\w\wedge\eta\in L^2\Omega^{p+2}(\ti M^{2n})(1)$. Hence
$$a_K\ti\w\wedge\eta\to\ti\w\wedge\eta \quad \text{in}\quad L^2\Omega^{p+2}(\ti M^{2n})(1)$$ as $K\to\infty$.

Now $da_K$ vanishes on $B_K$ and outside $B_{2K}$, $$|da_K|_{\ti g(1)}\leq \frac{\sqrt{a_K}}{K}, \quad and\quad |\b(x)|_{\ti g(1)}\leq C(\ti\rho_1(x_0, x)+1).$$ Thus
\begin{align} \label{dk}
\int_{\ti M^{2n}}|da_K\wedge\b\wedge\eta|^2_{\ti g(1)}\frac{\ti \w^n}{n!}&\leq Const.\int_{B_{2K}\backslash B_{K}}|\eta|^2_{\ti g(1)}\frac{\ti \w^n}{n!}\nonumber\\
&\leq Const.\int_{\ti M^{2n}\backslash B_{K}}|\eta|^2_{\ti g(1)}\frac{\ti\w^n}{n!}.
\end{align}
This converges to zero as $K\to\infty$. We thus have convergence of both terms on the right-hand side of (\ref{dj}) and consequently $d(a_K\b\wedge\eta)$ converges in $L^2$ to
$\ti \w\wedge\eta.$ Hence $\ti \w\wedge\eta$ lies in the closure of $dL^2\Omega^{p+1}(\ti M^{2n})(1)\cap L^2\Omega^{p+2}(\ti M^{2n})(1)$, and its $L^2$-cohomology class
with respect to the $\G$-invariant (measurable) Lipschitz K\"ahler flat metric $\ti g(0)$ vanishes since $[L_{\ti\w}, \Delta_0]=0$
on $\ti M^{2n}$ in the sense of distributions. By Theorem \ref{bj}, $0\leq p\leq n-1$, it implies $\eta=0$. By Theorem \ref{bi}, therefore, $$b^p_{(2)}(\ti M^{2n})(1)=b^{p}_{(2)}(\ti M^{2n})(0)=0\quad for \quad p\neq n.$$ This completes the proof of Theorem \ref{bk}.
\end{proof}
\brem
(1) Suppose that $(M^{2n}, \w, J, g_J)$ is symplectic hyperbolic, we still have $$ b^p_{(2)}(\ti M^{2n})(1)=b^p_{(2)}(\ti M^{2n})(0)=0\quad for\quad p\not=n.$$

(2)Teng Huang in \cite[Theorem 1.1]{Huang1} studied $L^2$ vanishing theorem for some twisted elliptic operators, the corresponding elliptic operators have no hard Lefschetz property. Using Huang's method
(cf. \cite{Huang2}), by Lefschetz decomposition into primative forms (cf. L. S. Tseng and S. T. Yau \cite{TY}), we have the following result:

Suppose that $(\ti M^{2n}, \ti\w, \ti J, \ti g_J)$ is a smooth complete almost K\"ahler structure on $\ti M^{2n}$. If $\ti\w$ is $d$(sublinear), then the space $\mathcal{H}^{i,j}_{(2)}(\ti M^{2n}, \ti g_J)$ vanishes for $i+j\neq n$, where $$\mathcal{H}^{i,j}_{(2)}(\ti M^{2n}, \ti g_J)\coloneqq\{\a\in L^2\W^{i,j}(\ti M^{2n}, \ti g_J)|d\a=0=d^{*}\a\},$$ where $d^*=-*_{\ti g_J}d*_{\ti g_J}$.

Notice that $$L^2\W^{k}(\ti M^{2n}, \ti g_J)=\underset{i+j=k}\oplus L^2\W^{i,j}(\ti M^{2n}, \ti g_J).$$ If $J$ is integrable, then $$\mathcal{H}^{i,j}_{(2)}(\ti M^{2n}, \ti g_J)=0, i+j=k,$$ imply that $\mathcal{H}^{k}_{(2)}(\ti M^{2n}, \ti g_J)=0$ since $$\mathcal{H}^{k}_{(2)}(\ti M^{2n}, \ti g_J)=\underset{i+j=k}\oplus \mathcal{H}^{i,j}_{(2)}(\ti M^{2n}, \ti g_J).$$ If $J$ is not integrable, then $$\underset{i+j=k}\oplus \mathcal{H}^{i,j}_{(2)}(\ti M^{2n}, \ti g_J)\subseteq\mathcal{H}^{k}_{(2)}(\ti M^{2n}, \ti g_J).$$ It is natural to ask the following question:

Can one directly prove  $\mathcal{H}^{k}_{(2)}(\ti M^{2n}, \ti g_J)=0$ for $k\neq n$ if $(\ti M^{2n},\ti\w,\ti J, \ti g_J)$ is $\G$-invariant and $d$(sublinear)?
\erem
By Theorems \ref{bi}, \ref{bk} and Atiyah's $\G$-index theorem (cf. \cite{At}, also see Theorem \ref{bh}), we obtain the following result:
\bthm \label{bl}
Suppose that $(M^{2n}, \w, J, g_J)$ is symplectic parabolic. Then
$$b^p_{(2)}(\ti M^{2n})(1)=b^{p}_{(2)}(\ti M^{2n})(0)=0, p\not= n,$$
$$ (-1)^n\chi(M^{2n}, g_J)=(-1)^n\chi_{(2)}(\ti M^{2n})(1)\geq 0.$$
\ethm
\vskip 0.1cm
It is easy to see that Theorem \ref{bl} gives the second part of Theorem \ref{main}.
\section{Symplectic hyperbolic manifolds}
In this section, we will study a non-vanishing theorem for $L^2$-cohomology and complete the proof of Theorem \ref{main}.
Suppose that $(M^{2n}, \w, J, g_J)$ is symplectic hyperbolic. Let $$\pi:(\ti M^{2n}, \ti\w, \ti J, \ti g_J)\to(M^{2n}, \w, J, g_J)$$
be the universal covering of $(M^{2n}, \w, J, g_J)$. As done in Section 3 and 4, construct a family of $\ti\w$-compatible almost complex structures $\ti J(t)$, $t\in[0, 1]$, on $\overset{\circ}{\ti M} {}^{2n}$ which is an open dense submanifold of $\ti M^{2n}$ such that
$\ti J(1)=\ti J$ and $\ti J(0)$ is integrable. Set $\ti g(t)=\ti\w(\cdot, \ti J(t)\cdot)$, $t\in[0, 1]$. Thus, we construct a family of $\G$-invariant (measurable) Lipschitz almost K\"ahler structures on $\ti M^{2n}$ which are with equivalent norms (that is, with Lipschitz condition \eqref{cn}). If $t=1$,
$(\ti\w, \ti J(1)=\ti J, \ti g(1)=\ti g_J)$ is a smooth $\G=\pi_1(M^{2n})$-invariant almost K\"ahler structure on $\ti M^{2n}$; if $t=0$, $(\ti\w, \ti J(0), \ti g(0))$ is a $\G$-invariant (measurable) Lipschitz K\"ahler flat structure on $\ti M^{2n}$. Since $(M^{2n}, \w)$ is symplectic hyperbolic in the sense of Definition 1.1, then there is a 1-form $\b$ on $\ti M^{2n}$ such that $\ti\w=d\b$ which is $d$(bounded), that is, $||\b||_{L^\infty}<C$, where $C$ is a constant on $\ti M^{2n}$ depending only on $(M^{2n}, \w, J)$.

Notice that $(\ti\w, \ti J(0), \ti g(0))$ is a $\G$-invariant (measurable) Lipschitz K\"ahler structure on $\ti M^{2n}$ which is equivalent to $(\ti\w,\ti J(1),\ti g(1))$ (that is, with Lipschitz condition \eqref{cn}), $\ti M {}^{2n}\backslash\overset{\circ}{\ti M} {}^{2n}$ has Lebesgue measure zero with respect to the metric $\ti g_J$, and $(\ti\w, \ti J(0), \ti g(0))$ is
a smooth $\G$-invariant K\"ahler flat structure on $\overset{\circ}{\ti M} {}^{2n}$. As done in the K\"ahler case, since $\ti\w$ is d(bounded) on $\ti M^{2n}$, we have the following theorem (cf. M. Gromov \cite[1.4.A Theorem]{Gr2},
or P. Pansu \cite[Theorem 8.2]{Pan}):

\bthm \label{et1}
Let $(\ti M^{2n}, \ti\w, \ti J, \ti g_J)$ be a complete  $\G$-invariant almost K\"ahler manifold of dimension $2n$. We can construct a family of $\G$-invariant (measurable) Lipschitz almost K\"ahler structures $(\ti\w, \ti J(t), \ti g(t))$ on $\ti M^{2n}$ (that is, $\ti g(t)$ is measurable metric on $\ti M^{2n}$ with Lipschitz condition \eqref{cn}), where if $t=1$, $(\ti\w, \ti J(1), \ti g(1))=(\ti\w, \ti J, \ti g_J)$, if $t=0$, $(\ti M^{2n}, \ti \w, \ti J(0), \ti g(0))$ is a $\G$-invariant (measurable) Lipschitz K\"ahler flat manifold. If
$\ti\w$ is d(bounded), as in K\"ahler case, by Theorem \ref{bi}, then $L^2\Omega^*(\ti M^{2n})(0)$ splits orthogonally as $$L^2\Omega^*(\ti M^{2n})(0)=\mathcal{H}^n_{(2)}(\ti M^{2n})(0)\oplus E(*)$$ and the Hodge-Laplacian with respect to
$\ti g(0)$ is invertible on $E(*)$ in the sense of distributions, where $$L^2\W^p(\ti M^{2n})(0)=E(p)\; \text{for}\;p\neq n,$$ $$L^2\Omega^n(\ti M^{2n})(0)=\mathcal{H}^n_{(2)}(\ti M^{2n})(0)\oplus E(n).$$ In particular, by vanishing theorem (cf. Theorem \ref{bk}), $\mathcal{H}^p_{(2)}(\ti M^{2n})(0)=0$ for $p\neq n$.
\ethm
\begin{proof}
It is similar to the proof of Theorem 8.2 in P. Pansu \cite{Pan} by Lemma \ref{lem1}, the proof is omitted here.
\end{proof}

\brem
By the construction of the family of $\G$-invariant (measurable) Lipschitz almost K\"ahler structures $(\ti\w, \ti J(t), \ti g(t))$, $t\in[0, 1]$, and Lemma \ref{lem7}, $\ti g(1)$ and $\ti g(0)$ are quasi-isometry on $\overset{\circ}{\ti M} {}^{2n}$ (cf. \eqref{cn}), thus by Theorem \ref{bi} Theorem \ref{et1} still holds for the smooth $\G$-invariant almost K\"ahler manifold $(\ti M^{2n}, \ti\w, \ti J(1), \ti g(1))$.
\erem
Hence, similar to the K\"ahler case, we have the following theorem:
\bthm \label{et2}
Let $(M^{2n}, \w, J, g_J)$ be a closed almost K\"ahler manifold of dimension $2n$, where the symplectic form $\w$ is symplectic hyperbolic. Let
$$\pi:(\ti M^{2n}, \ti \w, \ti J, \ti g_J)\to(M^{2n}, \w, J, g_J)$$ be the universal covering, then $(\ti\w, \ti J, \ti g_J)$ is a smooth $\G$-invariant almost K\"ahler structures on $\ti M^{2n}$, $$L^2\Omega^*(\ti M^{2n}, \ti g_J)=\mathcal{H}^n_{(2)}(\ti M^{2n}, \ti g_J)\oplus E(*)$$ and the Hodge-Laplacian with respect to
$\ti g_J$ is invertible on $E(*)$, where $$L^2\W^p(\ti M^{2n}, \ti g_J)=E(p) \quad for\quad p\neq n,$$ $$L^2\Omega^n(\ti M^{2n}, \ti g_J)=\mathcal{H}^n_{(2)}(\ti M^{2n}, \ti g_J)\oplus E(n).$$ In particular, $\mathcal{H}^p_{(2)}(\ti M^{2n}, \ti g_J)=0$ for $p\neq n$.
\ethm
\begin{proof}
Since $\ti g(1)$ and $\ti g(0)$ are quasi-isometry on $\overset{\circ}{\ti M} {}^{2n}$ by Lemma \ref{lem7}, hence, by Theorem \ref{bi}, Theorem \ref{et2} is a direct consequence of Theorem \ref{et1}.
\end{proof}
Notice that for $$L^2\Omega^n(\ti M^{2n}, \ti g_J)=\mathcal{H}^n_{(2)}(\ti M^{2n}, \ti g_J)\oplus E(n),$$ we can give a direct proof (cf. P. Pansu \cite[Theorem 8.2]{Pan}).

The Hodge decomposition
\begin{align*}
L^2\Omega^n(M^{2n}, \ti g_J)=&\mathcal{H}^n_{(2)}(M^{2n}, \ti g_J)\oplus\overline{dL^2\Omega^{n-1}(M^{2n}, \ti g_J)}\\
&\oplus\overline{d^{*}L^2\Omega^{n+1}(M^{2n}, \ti g_J)}
\end{align*}
holds in general on complete manifolds. Here it becomes
$$E(n)=\overline{dL^2\W^{n-1}(\ti M^{2n}, \ti g_J)}\oplus \overline{*_{\ti g_J}dL^2\W^{n-1}(\ti M^{2n}, \ti g_J)}.$$
Since the Hodge star operator $*_{\ti g_J}$ and the differential $d$ commute with the Hodge-Laplacian $\Delta_{\ti g_J}$, the inequality
$$||\a||^2_{L^2}\leq Const.\li\Delta_{\ti g_J}\a,\a\ri$$ holds on $E(n)$. This implies that $\Delta_{\ti g_J}$ is invertible on $E(n)$.
\brem
As done in T. Huang \cite[Theorem 1.3]{Huang1}, we have the following result (cf. T. Huang \cite{Huang2}):

Let $(\ti M^{2n}, \ti \w, \ti J, \ti g_J)$ be the universal covering of closed almost K\"ahler manifold $(M^{2n}, \w, J, g_J)$, where $\ti\w=d\ti\b$ is a $d$(bounded) symplectic form,
$\ti\b$ is a bounded 1-form. Then any $\a\in L^2\W^{i,j}(\ti M^{2n}, \ti g_J)$ of degree $k=i+j<n$ satisfies the inequality
$$C(n,k,\w,J)||\ti\b||^{-2}_{L^\infty(\ti M^{2n})}||\a||^2_{L^2(\ti M^{2n},\ti g_J)}\leq||d\a||^2_{L^2(\ti M^{2n},\ti g_J)}+||d^{\ti*}\a||^2_{L^2(\ti M^{2n},\ti g_J)},
$$
where $d^{\ti*}=-*_{\ti g_J}d*_{\ti g_J}$, $C(n,k,\w,J)$ is a constant which depends only on $n, k, \w,$ and $J$.

It is natural to ask the following question:

 Can one directly prove Theorem \ref{et2}?
\erem

Following M. Gromov \cite{Gr2}, similar to the K\"ahler case, we shall prove that, if $(M^{2n}, \w)$ is symplectic hyperbolic, then $\chi(M^{2n})\neq0$. In order to prove this, we use the signature operator (cf. Atiyah-Singer \cite[\S 6]{AS}):

Suppose that $(M^{2n}, g)$ is a closed oriented Riemannian manifold of dimension $2n$. Then $\Delta_g=dd^*+d^*d$ is the Hodge-Laplacian operator, where $d^*=-*_gd*_g$, $*_g$ is the Hodge star operator with respect to the metric $g$.
For any $\a\in\Omega^p(M^{2n})$, $0\leq p\leq 2n$, $$*_g^2\a=(-1)^p\a.$$ We now introduce a map on differential forms defined by
\beq \label{ed1}
\star \a=(\sqrt{-1})^{p(p-1)+n}*_g\a,\quad \a\in\Omega^p(M^{2n}).
\eeq
Note that $\star$ is real if $n$ is even, and imaginary if $n$ is odd. Then $$\star^2 \a=(\sqrt{-1})^\varepsilon*_g^2\a,$$ where
\begin{align*}
\varepsilon=&p(p-1)+n+(2n-p)(2n-p-1)+n\\
=&4n^2-4np+2p^2\equiv2p  \mod 4.
\end{align*}
Since $*_g^2\a=(-1)^p\a$, it follows that $\star^2 \a=\a$, and so $\star$ is an involution (cf. Atiyah-Singer \cite{AS}).
We can therefore decompose the space $$\Omega^*(M^{2n}, \C)=\overset{2n}{\underset{p=0}{\oplus}}\Omega^p(M^{2n}, \C)$$ into the $\pm$1-eigenvalues $\Omega^{\pm}$ of $\star$. Using
$$*_g^2\a=(-1)^p\a, \a\in\Omega^p(M^{2n}, \C)$$ and $d^*\a=-*_gd*_g\a$, one then verifies that
\beq \label{ed2}
(d+d^*)\star=-\star(d+d^*),\star\Delta_g=\Delta_g\star,  \Delta_g=(d+d^*)^2.
\eeq
So $d+d^*$ maps $\Omega^+$ into $\Omega^-$ and $\Omega^-$ into $\Omega^+$. Where let $\Lambda^{\pm}T^*M^{2n}\otimes\C$ be the subspaces of $\Lambda^*T^*M^{2n}\otimes\C$ satisfying
$$\star\Lambda^{\pm}T^*M^{2n}\otimes\C=\pm\Lambda^{\pm}T^*M^{2n}\otimes\C,$$ hence
\beq \label{ed21}
\Lambda^*T^*M^{2n}\otimes\C=\Lambda^+T^*M^{2n}\otimes\C\oplus\Lambda^-T^*M^{2n}\otimes\C,
\eeq
\beq \label{ed22}
\Omega^+\coloneqq\Gamma(\Lambda^+T^*M^{2n}\otimes\C),\quad\Omega^-\coloneqq\Gamma(\Lambda^-T^*M^{2n}\otimes\C).
\eeq
We denote by $D^{\pm}$ the restriction of $d+d^*$ to the appropriates subspace of $\Omega^*$. Thus
\beq \label{ed23}
D^+:\Omega^+\to\Omega^-,\quad D^-:\Omega^-\to\Omega^+.
\eeq
Each is elliptic, and they are formal adjoint of each other. $D^{\pm}$ are called the signature operators. Moreover, the solution of $D^+u=0$ are just the $\Delta_g$-harmonic forms in $\Omega^+$, and similarly for $\Omega^-$. We
shall denote by $H^+$ and $H^-$ these subspaces of $\Delta_g$-harmonic forms, and we put
$$h^+\coloneqq\dim H^+,\quad h^-\coloneqq\dim H^-.$$
Thus we have $\rm index D^+=h^+-h^-$. Since $\star\Delta_g=\Delta_g\star$, it follows that $\star$ introduces an involution on the $\Delta_g$-harmonic forms, the $\pm$ eigenspaces are
precisely $H^+$ and $H^-$. Thus,
\beq \label{ed24}\dim\ker\Delta_g=h^++h^-.\eeq

Let $L$ be a complex line bundle over $(M^{2n}, g)$, and $\nabla$ be a connection on $L$. Then we can define the twisted signature operator
\beq \label{ed3}
D^+\otimes\nabla:\Gamma(\Lambda^+T^*M^{2n}\otimes L)\to\Gamma(\Lambda^-T^*M^{2n}\otimes L),
\eeq
where $$\star\Lambda^{\pm}T^*M^{2n}\otimes \C=\pm\Lambda^{\pm}T^*M^{2n}\otimes \C.$$ Hence we have the local index theorem (cf. Atiyah-Bott-Patodi \cite[\S 6]{ABP} or Y. L. Yu \cite[\S 8.4]{Yu}):
$$
\rm loc. index(D^+\otimes\nabla)=\mathscr{L}\cup \exp\frac{\sqrt{-1}R^\nabla}{2\pi},
$$
where $R^\nabla$ is the curvature of $\nabla$, $\mathscr{L}=\sum\limits_kL_k$ is Hirzebruch $L$-class, the polynomials $L_k(p_1, \cdots, p_k)$ defining by the generating functions
$$\sum L_k\coloneqq\Pi\frac{x_j}{\tanh x_j},$$ where $p_k$ is the Pontrjagin class, $\sum p_k=\Pi(1+x_j^2)$. In particular, $$L_0=1,\quad L_1=\frac{1}{3}p_1, \quad and \quad L_2=\frac{1}{45}(7p_2-p_1^2).$$ For more details, see F. Hirzebruch \cite[\S 1.4]{Hir}.

Suppose that $(M^{2n}, \w, J, g_J)$ is a closed almost K\"ahler manifold of dimension $2n$, where $\w$ is symplectic hyperbolic. Then the universal covering $$\pi:(\ti M^{2n}, \ti \w, \ti J, \ti g_J)\to(M^{2n}, \w, J, g_J)$$ is $\G$-invariant, where $\G=\pi_1(M^{2n})$ is an infinite discrete group. Since $\ti\w=d\b$, $\b$ is a bounded 1-form on $\ti M^{2n}$, we may assume that $\ti\w$ is the Chern class of a trivial complex line bundle $\ti L$ over $\ti M^{2n}$, or $\ti\w$ can be regarded as the curvature
of the connection $d+\sqrt{-1}\b$ on $\ti L\cong \ti M^{2n}\times \mathds{C}$. S. K. Donaldson gave an approach \cite[p.667]{Don}: for a general symplectic manifold $(M, \w)$ is to choose a compatible almost
complex structure $J$ on $M$, and then to extend familiar results about positive line bundles suitably formulated to the almost complex case. Hence, we can assume that de Rham cohomology class
$$[\frac{\w}{2\pi}]\in H^2(M, \R)$$ lies in the integral lattice $H^2(M, \Z)/Torsion$.

As done in the K\"ahler case (cf. M. Gromov \cite{Gr2}) let $$E^{\pm}\coloneqq\Gamma(\Lambda^\pm T^*\ti M^{2n}\otimes\C).$$ The signature operator $D^+:E^+\to E^-$ is the restriction $d+d^*$ to the $E^+$, where $d^*=-*_{\ti g_J}d*_{\ti g_J}$, $*_{\ti g_J}$ is the Hodge star operator defined by the metric $\ti g_J$.
Consider the perturbation $D^+\otimes\nabla$ of $D^+$
\beq \label{ee}
D^+\otimes\nabla:E^+\otimes\ti L\to E^-\otimes\ti L,
\eeq where $D^+\otimes\nabla$ is called a twisted signature operator.
Since $\ti L$ is trivial and $\b$ is bounded, formally we define $$\nabla^{(k)}=d+\frac{1}{k}\sqrt{-1}\b, k\in \mathds{Z},$$ which is a unitary connection on the trivial line bundle $\ti L=\ti M^{2n}\times\C$. One can try to make it $\G$-invariant by changing to a non trivial action of $\G$ on $\ti M^{2n}\times\C$, i.e., setting, for $\g\in\G$, $$\g(\ti x,z)=(\g\ti x, e^{\sqrt{-1}u(\g,\ti x)}z).$$ We want $\g^*\nabla^{(k)}=\nabla^{(k)}$, i.e., $du=-(\g^*\b-\b)$. Since $$d(\g^*\b-\b)=\g^*\ti\w-\ti\w=0,$$ there always exists a solution $u(\g,\cdot)$ well defined up to a constant. Hence we have the following definition:
\bdefn  Let $G_k$ be the subgroup of ${\rm Diff}(\ti M^{2n}\times\mathds{C})$ formed by maps $g$ which are linear unitary on fibres, preserve connection $\nabla^{(k)}$ and cover an element of $\G=\pi_1(M^{2n})$.
\edefn
For more details, see M. Gromov \cite[Chapter 2]{Gr2}, also see P. Pansu \cite{Pan}. Notice that the locally compact Lie group $G_k$ is unimodular \cite[p.59]{DK}. By the construction we have an exact sequence
\beq \label{ef}
1\to U(1)\to G_k\to \G\to 1.
\eeq
The exact sequence (\ref{ef}) is called a central co-extension of $\G$ by $U(1)$ \cite[p.86]{DK}.
Since sections of the line bundle $\ti M^{2n}\times\mathds{C}\to \ti M^{2n}$ can be viewed as $U(1)$ equivariant functions on $\ti M^{2n}\otimes U(1)$, the operator $D^+\otimes\nabla^{(k)}$ can be viewed as
a $G_k$ invariant operator on the Hilbert space $\mathcal{H}$ of $U(1)$ equivariant basic $L^2$-differential forms on $\ti M^{2n}\times\mathds{C}$.

One can define a projective von Neumann dimension for $G_k$-invariant subspace in a projective representative, and state $G_k$-index theorem (that is a local index theorem which gives the ``index per unit volume"
of the twisted signature operator in the question) (cf. \cite{BGV,Luck,Pan}).
The following theorem is a particular case of $L^2$ index theorem for $G$ invariant operator, where $G$ is a unimodular Lie group (with countably many components) acting properly and freely on $\ti M^{2n}$ with a compact
quotient. Although they do not state theorem in this generality, A. Connes and H. Moscovici \cite[5.2 Theorem]{CM} provide all the necessary ingredients.
\bthm {\rm(}cf. M. Gromov \cite[2.3.A Theorem]{Gr2}{\rm)}\label{bn} The twisted signature operator $D^+\otimes\nabla^{(k)}$ has a finite projective $L^2$-index given by
\begin{align*}
L^2{\rm index}_{G_k}(D^+\otimes\nabla^{(k)})&=\int_{M^{2n}}\mathscr{L}_{M^{2n}}\cup \exp\frac{1}{2k\pi}[\w]\\
&=\int_{\pi(F)}\mathscr{L}_{M^{2n}}\cup \exp\frac{1}{2k\pi}[\w],
\end{align*}
where $[\w]$ is an equivalence class of symplectic form $\w$, $\mathscr{L}_{M^{2n}}$ is the Hirzebruch class (cf. \cite{Hir}), $\mathscr{L}_{M^{2n}}=1+\cdots+e(M^{2n})$, and $e(M^{2n})$ is the Euler class of $M^{2n}$.
\ethm
\brem
In fact, N. Teleman \cite[Theorem 1.1, Theorem 6.3]{Tele4} extended the Hirzebruch-Atiyah-Singer signature theorem and the index theorem for the abstract elliptic operators on topological manifolds. He showed that the Hirzebruch-Atiyah-Singer signature theorem and the index theorem for the abstract elliptic operators
with value in continuous vector bundles over topological manifolds can be expressed in topological terms by the same formula as in the smooth case.
\erem
The proof of Theorem \ref{bn} follows from the heat equation method of M. Atiyah, R. Bott and V. Patodi \cite{ABP}. For details, see \cite{ABP,BGV,Pan}.
In Theorem \ref{bn}, the index is a polynomial in $\frac{1}{k}$ whose highest degree form is
$$\int_{\pi(F)}\left(\frac{\w}{2\pi}\right)^n=\int_{M^{2n}}\left(\frac{\w}{2\pi}\right)^n\not=0,$$ thus for $\frac{1}{k}$ small enough, $D^+\otimes\nabla^{(k)}$ has a non zero $L^2$-kernel. Therefore, $D^+\otimes\nabla^{(k)}$ is
not invertible. By the construction and Theorem \ref{et2} (In particular, the Hodge-Laplacian operator on $L^2\W^p(\ti M^{2n},\ti g_J)$ is invertible for $p\neq n$), $D^+\otimes\nabla^{(k)}$ is a  $\frac{1}{k}$-small
perturbation of $$D^+=d+d^*:E^+\to E^-,$$ so by Lemma \ref{lem4} below $d+d^*$ is not invertible too. Hence $D^+=d+d^*$ has a non zero $L^2$-kernel. On the other hand, by Theorem \ref{et2} it implies that $\mathcal{H}^p_{(2)}(\ti M^{2n}, \ti g_J)=0$ for $p\neq n$, thus $\mathcal{H}^n_{(2)}(\ti M^{2n}, \ti g_J)\not=0$. Since $\frac{1}{k}\b\to0$ as $k\to\infty$, we have the following lemma:
\blem {\rm(}cf. M. Gromov \cite[\S2.4]{Gr2}{\rm)}\label{lem4}
The operator $D^+\otimes\nabla^{(k)}$ will converge to  $D^+=d+d^*$ as $k\to \infty$.
\elem
By Theorem \ref{bn} and Lemma \ref{lem4}, we have the following corollary:
\bcor \label{bo}
Let $\pi:(\ti M^{2n}, \ti\w, \ti J, \ti g_J)\to(M^{2n}, \w, J, g_J)$ be the Riemannian covering of $(M^{2n}, \w, J, g_J)$. Assume that $(M^{2n}, \w, J, g_J)$ is symplectic hyperbolic.
Then $\ti M^{2n}$ admits non zero $L^2$-$(d, d^{\ti*})$-harmonic $n$-forms.
\ecor
In terms of Corollary \ref{bo}, by Theorem \ref{bi} and Atiyah's $\G$-index theorem, we have:
\bthm \label{bp}
Let $(M^{2n}, \w, J, g_J)$ be a closed symplectic hyperbolic manifold, then
$$
(-1)^n\chi(M^{2n})=(-1)^n\chi_{(2)}(\ti M^{2n})(1)=(-1)^n\chi_{(2)}(\ti M^{2n})(0)>0.
$$
\ethm
Therefore, by Theorem \ref{bl} and Theorem \ref{bp}, we obtain one of our main results Theorem \ref{main}.
\vskip 0.2cm

\section{The Proof of Theorem 1.9}
This section is devoted to proving Theorem \ref{mc}. We consider even dimensional closed Riemannian manifolds with non-positive curvature.

If $(M^{2n}, g)$ is a $2n$-dimensional closed Riemannian manifold with infinite fundamental group $\pi_1(M^{2n})$, then the universal covering
$$\pi:(\ti M^{2n}, \ti g)\to(M^{2n}, g)$$ is a complete simply connected manifold with deck transformation group $\G=\pi_1(M^{2n})$.

For $(M^{2n}, g)$, as done in Section 2, we can construct a family of almost K\"ahler structures on an open dense submanifold of $M^{2n}$ which are quasi isometric. More precisely, $\forall p\in M^{2n}$, there exists a coordinate chart
$$B(p,3)\supset\supset U_p\supset B(p,1)=\{x\in U_p:\rho_g(p,x)<1\},$$ $$\psi_p:U_p\to\psi_p(U_p)\subset\R^{2n}\cong\C^n,$$ where $\rho_g$ is distance function with respect to the Riemannian metric $g$ such that a family of almost K\"ahler structures
$(\w'_p,J'_p(t),g'_p(t))$, $t\in[0,1]$, where $\w'_p=\psi^*_p\w_0$, $\w_0$ is the standard K\"ahler form on $\R^{2n}\cong\C^n$, $J'_p(1)=J_{g,p}$ which is defined by $\w'_p$ and $g|_{U_p}$,
$$\w'_p(\cdot,J_{g,p}\cdot)=g'_p(1)=g|_{U_p},$$ (notice that for the triple  $(\w, J, g)$ any one of the pairs $(\w,J)$, $(J,g)$, or $(g,\w)$ determines the other two.) $J'_p(0)=\psi^*_pJ_0$,
$J_0$ is the standard complex structure on $\R^{2n}\cong\C^n$, $g'_p(0)=\w'_p(\cdot, J'_p(0)\cdot)$, $(\w'_p,J'_p(0),g'_p(0))$ is a K\"ahler flat structure on $U_p$ (cf. Section 2, McDuff-Salamon \cite{MS}).

Let $S'_p=-J'_p(0)\circ J_{g,p}=-J'_p(0)\circ J'_p(1)$ which is a symmetric positive definite and symplectic matrix on $U_p$ (cf. McDuff-Salamon \cite[Chapter 2, Lemma 2.5.5]{MS}), that is, $S'_pJ'_p(0)S'_p=J'_p(0)$, and $J'_p(1)=J'_p(0)S'_p$. Let $h'_p=g'_p(0)\ln S'_p$. As in Section 2, the $C^k$-norms of $S'_p$ and $h'_p$ with respect to $g$ and its Levi-Civita connection are bounded from above by $C(M^{2n}, g, k)$ which depends only on  $M^{2n}$, $p$, $g$, and $k$. Define $J'_p(t)=J'_p(0)e^{th'_p}$, $t\in[0,1]$ on $U_p$.
Let $g'_p(t)(\cdot,\cdot)=\w'_p(\cdot, J'_p(t)\cdot)$. If $U_{p_1}\cap U_{p_2}\neq\emptyset$, for $x\in U_{p_1}\cap U_{p_2}$, $g'_{p_1}(1)(x)=g'_{p_2}(1)(x)=g(x)$. As in Section 2 (cf. \eqref{a18}), for $t\in[0,1]$,
\beq \label{ea12}
||J'_p(t)||_{C^k(U_p,g)},||g'_p(t)||_{C^k(U_p,g)},||g'^{-1}_p(t)||_{C^k(U_p,g)}\leq C(M^{2n},p,g,k),
\eeq
where $C(M^{2n},p,g,k)$ depends only on $M^{2n}$, $p$, $g$, and $k$. Let $*'_{p,t}$ be the Hodge star operator with respect to the metric $g'_p(t)$ on $U_p$, therefore
$||*'_{p,t}||_{C^k(U_p,g)}$, $t\in[0,1]$, are bounded by $C(M^{2n},p,g,k)<+\infty$.

Similarly, as in Section 2 (cf. \eqref{a19}-\eqref{a20}), since $g'_p(1)=g$, for $x\in U_p$, $X\in T_xU_p$,
\beq \label{ea13}
g'_p(t)(x)(X,X)\leq C(M^{2n},p,g)g(x)(X,X),
\eeq
and
\beq \label{ea14}
g(x)(X,X)\leq C(M^{2n},p,g)g'_p(t)(x)(X,X),
\eeq
where $t\in [0,1]$ and $C(M^{2n},p,g)$ depending only on $M^{2n}$, $p$ and $g$.

As done in Section 2, for $(M^{2n},g)$, one can find a partition $P(M^{2n},g)$ of $M^{2n}$: Choose a finite subatlas $\{\psi_i, U_i\}_{1\leq i\leq N}$ of $M^{2n}$. Define $W_1=U_1$,  $T_1=U_1$,
$$W_i=U_i \setminus\bigcup^{i-1}_{j=1}\overline{U_j},\quad T_i=U_i \setminus\bigcup^{i-1}_{j=1}U_j,\quad 2\leq i\leq N.$$
Then it is easy to see that $W_i\subset T_i\subset U_i$, $1\leq i\leq N$. $W_i$ is an open subset of $M^{2n}$, $$T_i\cap T_j=\emptyset,\quad i\neq j,\quad M^{2n}=\bigcup\limits^N_{i=1}T_i.$$ Let
$$\overset{\circ}{M}{}^{2n}=\overset{\circ}{M}{}^{2n}(g, P)\coloneqq\bigcup\limits^N_{i=1}W_i$$ be an open dense submanifold of $M^{2n}$, $M^{2n}\setminus\overset{\circ}{M}{}^{2n}$ has Lebesgue measure zero. One can define a family of almost K\"ahler structures on each coordinate system $\{\psi_i, U_i\}$   $$(\w'_i,J'_i(t),g'_i(t)),\quad 1\leq i\leq N, t\in[0,1].$$
Where $\w'_i=\psi'^{*}_i \w_0$, $g'_i(t)(\cdot, \cdot)=\w'_i(\cdot, J'_i(t)\cdot)$, and $g'_i(1)(\cdot, \cdot)=g(\cdot, \cdot)$, where
$$B(p_i,3)\supset\supset U_i\supset B(p_i,1), 1\leq i\leq N, \{p_1, ...,p_N\}\subset M^{2n}. $$
$\{U_i\}_{1\leq i\leq N}$ is an open cover of $M^{2n}$. For $x\in U_i$, $ X\in T_xU_i$, $1\leq i\leq N,$
\begin{align}\label{ea10}
C^{-1}g'(1)(x)(X,X)\leq g'(t)(x)(X,X)\leq C g'(1)(x)(X,X),\nonumber\\ t\in[0,1),g'(1)=g,
\end{align}
on the coordinate system $\{\psi_i, U_i\}$ (cf. Section 2). Hence, one can define a family of measurable almost K\"ahler metrics with Lipschitz condition \eqref{ea10}.
More precisely, for partition, $\{T_i\}_{1\leq i\leq N}$, of $M^{2n}$, transfer on each $T_i$ a family of almost K\"ahler metrics on $\psi_i(U_i)\subset$$\R^{2n}$ ($(\w'_i,J'_i(t),g'_i(t))$, $t\in[0,1]$) via diffeomorphism $$\psi_i:U_i\to\R^{2n},$$
$$(\w',J'(t),g'(t))\coloneqq\{(\w'_i|_{T_i},J'_i(t)|_{T_i},g'_i(t)|_{T_i})\}_{1\leq i\leq N},$$
where $$g'_i(1)(x)=g'_j(1)(x)=g(x),\quad \forall x\in U_i\cap U_j\neq\emptyset.$$

Restricted to $\overset{\circ}{M}{}^{2n}$, $\{(\w',J'(t),g'(t))|t\in[0,1]\}$ is a family of almost K\"ahler structures which are quasi isometry (see \eqref{ea10}), $(\w',J'(0),g'(0))$ is a K\"ahler flat structure. Similar to the symplectic case in Section 2, $g'(t)$, $t\in[0,1]$, are called (measurable) Lipschitz metircs on $M^{2n}$ which are with equivalent norms (see \eqref{ea10}).

Hence one can define a family of (measurable) Lipschitz almost K\"ahler structures $(\w', J'(t), g'(t))$, $t\in[0,1]$, on $M^{2n}$ off the Lebesgue measure zero subset $M^{2n}\backslash\overset{\circ}{M} {}^{2n}$ (cf. Section 2), where $g'(1)=g$ is smooth on $M^{2n}$, restricted to  $\overset{\circ}{M}{}^{2n}$, $(\w',J'(t),$ $g'(t))$, $t\in[0,1]$, are almost K\"ahler structures which are quasi isometry (cf. \eqref{ea10}), $(\w', J'(0), g'(0))$ is a K\"ahler flat structure on $\overset{\circ}{M}{}^{2n}$ which can be regarded as a singular K\"ahler structure on $M^{2n}$. Hence $\{g'(t),$ $t\in[0,1]\}$ is a family of (measurable) Lipschitz Riemannian metrics on $M^{2n}$ with equivalent norms. Note that, in general, $\w'$, $J'(t)$, and $g'(t)$, $t\in[0,1]$, can not be continuously extended to $M^{2n}$. In particular, in general $g'(t)$ is not a Lipschitz metric, but a measurable metric on $M^{2n}$. Restricted to the open dense submanifold $\overset{\circ}{M}{}^{2n}$, $$\text{coefficients of }g'(t)\text{ are in }C^\infty(\overset{\circ}{M}{}^{2n})\cap L^\infty_k(M^{2n},g),$$ where $k$ is a nonnegative integer.

As done in Section 3, for $\pi:\ti M^{2n}\to M^{2n}$, one can choose $$\ti W_i\subset \ti T_i\subset \ti U_i,\quad 1\leq i\leq N$$ such that $$\pi(\ti W_i)=W_i, \pi(\ti T_i)=T_i, \pi(\ti U_i)=U_i,\quad 1\leq i\leq N.$$ Define
$F=\bigcup\limits_{j=1}^N \ti W_j$  as a fundamental domain of $$\pi: (\ti M^{2n},\ti g) \to (M^{2n},g).$$ Let
$$\overset{\circ}{\ti M} {}^{2n}=\overset{\circ}{\ti M} {}^{2n}(g,P)\coloneqq\bigcup_{\gamma\in \Gamma} \gamma F.$$
Then $$\ti M^{2n}=\bigcup\limits_{\gamma\in \Gamma} \gamma \overline{F},$$ $\overset{\circ}{\ti M} {}^{2n}$ is an open dense submanifold of $\ti M^{2n}$ and $\ti M^{2n}\backslash\overset{\circ}{\ti M} {}^{2n}$ has Lebesgue measure zero with respect to the metric $g$. Let
\beq\label{ea0}
(\ti\w', \ti J'(t), \ti g'(t))=\pi^*(\w',J'(t),g'(t)),\quad t\in [0,1],
\eeq
be a family of $\G$-invariant (measurable) Lipschitz almost K\"ahler structures with equivalent norms, in particular, restricted to $\overset{\circ}{\ti M} {}^{2n}$, $(\ti\w', \ti J'(t), \ti g'(t))$, $t\in [0,1]$, are $\G$-invariant quasi isometry satisfying Lipschitz condition (cf. \eqref{cn}):
\begin{align}\label{ea11}
C^{-1}\ti g'(t)(x)(X,X)\leq \ti g'(1)(x)(X,X)\leq C \ti g'(t)(x)(X,X),\nonumber\\\forall x\in\overset{\circ}{\ti M} {}^{2n}, X\in T_x\overset{\circ}{\ti M} {}^{2n},t\in[0,1),
\end{align}
where $C>1$ is a constant depending only on $M^{2n}$ and $g$, on $M^{2n}$ $\ti g'(1)=\ti g$, on $\overset{\circ}{\ti M} {}^{2n}$
$(\ti\w', \ti J'(0),$ $\ti g'(0))$ is K\"ahler flat, and
\beq \label{ea1}
dvol_{\ti g'(t)}=\frac{\ti \w'^n}{n!}=dvol_{\ti g'(1)}=dvol_{\ti g}
\eeq
is a $\G$-invariant smooth volume form on $\overset{\circ}{\ti M} {}^{2n}$. Notice that $\ti g'(0)$ is a measurable metric only with Lipschitz condition \eqref{ea11}. As in Section 3, we have
$$
||\ti J'(t)||_{C^k(\overset{\circ}{\ti M} {}^{2n},\ti g)}+||\ti g'(t)||_{C^k(\overset{\circ}{\ti M} {}^{2n},\ti g)}+||\ti g'^{-1}(t)||_{C^k(\overset{\circ}{\ti M} {}^{2n},\ti g)}\leq C(M^{2n},g,k).
$$
Hence $$\text{coefficients of }\ti g'(t), \ti*'_t\text{ are in }C^\infty(\overset{\circ}{\ti M} {}^{2n})\cap L^\infty_k(\ti M^{2n}, \ti g),$$ where $\ti *'_t$ is the Hodge star operator with respect to the metric $\ti g'(t)$, $t\in[0,1]$, and $k$ is a nonnegative integer.

In summary, we have the following proposition (cf. Section 3):
\bprop\label{bi1}
Suppose that $(M^{2n}, g)$ is a $2n$-dimensional closed Riemannian manifold with infinite fundamental group $\pi_1(M^{2n})$ of $M^{2n}$. On the universal covering space $(\ti M^{2n}, \ti g)$ of $(M^{2n}, g)$ one can construct
a family of $\G=\pi_1(M^{2n})$ invariant (measurable) Lipschitz almost K\"ahler structures $(\ti\w', \ti J'(t), \ti g'(t))$, $t\in[0, 1]$ with equivalent norms (that is, with Lipschitz condition \eqref{ea11}) on $\overset{\circ}{\ti M} {}^{2n}$. If $t=1$, $$\ti\w'(\cdot, \ti J'(1)\cdot)=\ti g'(1)=\ti g=\pi^*g$$ is smoothly extended to $\ti M^{2n}$. Restricted to the $\G$-invariant open dense submanifold $\overset{\circ}{\ti M} {}^{2n}$ of $\ti M^{2n}$, $(\ti\w', \ti J'(t), \ti g'(t))$, $t\in[0, 1]$, is a family of $\G$-invariant almost K\"ahler structures which are quasi isometric satisfying \eqref{ea11}. If $t=0$,  $(\ti\w', \ti J'(0), \ti g'(0))$ is a $\G$-invariant K\"ahler flat structure on $\overset{\circ}{\ti M} {}^{2n}$.
In particular, $$\text{coefficients of }\ti g'(t), \ti*'_t\text{ are in }C^\infty(\overset{\circ}{\ti M} {}^{2n})\cap L^\infty_k(\ti M^{2n}, \ti g),$$ where $k$ is a nonnegative integer.
Note that $\ti M {}^{2n}\backslash\overset{\circ}{\ti M} {}^{2n}$ has Lebesgue measure zero with respect to the metric $\ti g$. $$dvol_{\ti g'(t)}=\frac{\ti\w'^n}{n!}=dvol_{\ti g}$$ is smooth $\G$-invariant on $\ti M^{2n}$, $t\in[0, 1]$.
\eprop

Suppose that $(M^{2n},g)$ is a $2n$-dimensional closed manifold with nonpositive (resp. strictly negative) sectional curvature, then the universal covering $$\pi: (\ti M^{2n},\ti g) \to (M^{2n},g)$$ is a complete simply connected manifold with an infinite deck transformation group $\G=\pi_1(M^{2n})$.
Notice that $\ti g$ is $\G$-invariant Riemannian metric on $\ti M^{2n}$ with bounded geometry. By the construction of $(\ti\w',\ti J'(1),\ti g'(1))$ (where $\ti g'(1)=\ti g$) it is easy to see that $||\ti \w'||_{L^\infty_k(\overset{\circ}{\ti M} {}^{2n},\ti g)}$ is finite (cf. Section 3), where $k$ is a nonnegative integer. It is worth to remarking that when $(M^{2n},g)$ has non-positive curvature, $\ti M^{2n}\cong\R^{2n}\cong\C^n$. For this case, in general, one can not construct $\G$-invariant $\ti\w'$, $\ti J'(t)$, and $\ti g'(t)$, $t\in[0,1)$ on $\ti M^{2n}$.

\vskip 0.2cm
As done in Lemma \ref{lem7}, by Proposition \ref{bi1}, it is easy to see that the operator $\ti*'_t=*_{\ti g'(t)}$ with respect to the $\G$-invariant (measurable) Lipschitz metric $\ti g'(t)$, $t\in [0,1]$, its coefficients are in $L^\infty_k(\ti M^{2n},\ti g)$ with respect to the metric $\ti g'(1)=\ti g$ since $\{\ti U_i\}_{1\leq i\leq N}$ is an open cover of $\bar{F}$ due to the construction of $\ti g'(t)$, where $k$ is a nonnegative integer. Thus for the family of $\G$-invariant (measurable) Lipschitz
almost K\"ahler structures $(\ti\w', \ti J'(t), \ti g'(t))$, $t\in[0,1]$, which are with equivalent norms satisfying Lipschitz condition \eqref{ea11}, as done in Lemma \ref{lem7}, we can define $\G$-invariant Hilbert space $L^2_k(\ti g'(t))$ on $\ti M^{2n}$ (which is well-defined), $t\in[0,1]$, and it is similar to \eqref{cn2} we have
\beq \label{cy0}
C(M^{2n}, g, k)^{-1}\li\a, \a\ri_{L^2_k(\ti g'(t))}\leq\li\a, \a\ri_{L^2_k(\ti g'(1))}\leq C(M^{2n}, g, k)\li\a, \a\ri_{L^2_k(\ti g'(t))},
\eeq
for any $ \a\in\Omega^p_c(\ti M^{2n})$, $t\in [0, 1)$, $0\leq p\leq 2n$, where $k$ is a nonnegative integer, and $C(M^{2n}, g, k)>1$ is a constant depending only on $k, M^{2n}$, and $g$.
The Hodge star operator $\ti*'_t$ with respect to the $\G$-invariant (measurable) Lipschitz metric $\ti g'(t)$ on $\ti M^{2n}$, $t\in [0, 1)$, coefficients of $\ti*'_t$ are in $L^\infty_k(\ti M^{2n}, \ti g)$, where $k$ is a nonnegative integer.
Since $(\ti M^{2n}, \ti g)$ is a complete $\G$-invariant Riemannian manifold, $$L^2_k\Omega^p(\ti M^{2n}, \ti g'(t)), t\in [0, 1), 0\leq p\leq 2n,$$ $k$ being nonnegative integer, are
defined as the completion of $\Omega^p_c(\ti M^{2n})$ with respect to the norms $||\cdot||_{L^2_k(\ti g'(t))}$(cf. Section 3). Hence, as in Section 3 by \eqref{cy0}, $$L^2_k\Omega^p(\ti M^{2n}, \ti g'(t)), t\in [0, 1),\quad and\quad L^2_k\Omega^p(\ti M^{2n}, \ti g)$$ are quasi isometry or with equivalent norms. Therefore, it is similar to Lemma \ref{lem1}, we have the following lemma:
\blem \label{lem5}
Let $(\ti\w', \ti J'(t), \ti g'(t))$, $t\in[0,1]$, be a family of $\G$-invariant Lipschitz almost K\"ahler structures which are with equivalent norms (that is, with Lipschitz condition \eqref{ea11}), where $(\ti M^{2n}, \ti g'(1))$$=$$(\ti M^{2n}, \ti g)$ is a complete smooth $\G$-invariant Riemannian manifold,
and $(\ti\w', \ti J'(0), \ti g'(0))$ is a $\G$-invariant (measurable) Lipschitz K\"ahler flat structure on $\ti M^{2n}$. Restricted to $\overset{\circ}{\ti M} {}^{2n}$, $\ti g'(t)$, $t\in [0,1)$, and $\ti g'(1)=\ti g$ are quasi-isometric. On $\ti M^{2n}$
$$L^2_k\Omega^p(\ti M^{2n}, \ti g'(t)), t\in[0,1)\quad\text{and}\quad L^2_k\Omega^p(\ti M^{2n}, \ti g), 0\leq p\leq 2n$$ are quasi-isometric or with equivalent norms.
Let $d^{\ti*'_t}$ be $L^2$-adjoint operator of $d$ on $\ti M^{2n}$ with respect to the (measurable) Lipschitz almost K\"ahler metric $\ti g(t)$ in the sense of distributions, then $d^{\ti*'_t}=-\ti*'_td\ti*'_t$ and $(d^{\ti*'_t})^{\ti*'_t}=d$ in the sense of distributions, where $\ti*'_t=*_{\ti g'(t)}$ is the Hodge star operator with respect to the (measurable) Lipschitz metric $\ti g'(t)$ on $\ti M^{2n}$.
\elem
By Lemma \ref{lem5}, similar to the smooth case, the Hodge-Laplacian $\Delta'_t=dd^{\ti*'_t}+d^{\ti*'_t}d$ is also well-defined for $t\in[0, 1]$ in the sense of distributions. It is easy to see that Lemma \ref{lem7} and \ref{lem1} still hold for this situation. If an $L^2$ $p$-form $\a$ on $\ti M^{2n}$ is $\Delta'_t$-harmonic, then $\a$ is $(d, d^{\ti*'_t})$-harmonic in the sense of distributions, where
$0\leq p\leq 2n$ and $t\in [0, 1]$.

By Lemma \ref{lem5}, for $\G$-invariant (measurable) Lipschitz almost K\"ahler manifold $$(\ti M^{2n}, \ti\w', \ti J'(t), \ti g'(t)),\quad t\in[0, 1],$$
\beq \label{cy1}
\mathcal{H}^p_{(2)}(\ti M^{2n}, \ti g'(t))\coloneqq\ker \Delta'_t|_{L^2\Omega^p(\ti M^{2n}, \ti g'(t))}, t\in[0, 1].
\eeq
It is not hard to see that (cf. Section 3) $$\mathcal{H}^p_{(2)}(\ti M^{2n}, \ti g'(t))\subset\Omega^p(\ti M^{2n})\cap L^2(\ti M^{2n}, \ti g), 0\leq p\leq 2n, t\in[0, 1].$$ One also defines a Hilbert $\G$-cochain complex $L^2_{l-*}\Omega^*(\ti M^{2n}, \ti g'(t))$ as follows:
\begin{align} \label{cy2}
0\to L^2_{l}\Omega^0(\ti M^{2n}, \ti g'(t))&\xrightarrow{d}L^2_{l-1}\Omega^1(\ti M^{2n}, \ti g'(t))\xrightarrow{d}\cdots\nonumber\\
&\xrightarrow{d}L^2_{l-2n}\Omega^{2n}(\ti M^{2n}, \ti g'(t))\to 0,\quad l\geq2n.
\end{align}
Set
\beq \label{cy3}
\mathrm{H}^p_{(2)}(\ti M^{2n}, \ti g'(t))\coloneqq\left. \ker d|_{L^2\Omega^p(\ti M^{2n}, \ti g'(t))}\middle/\overline{dL^2\Omega^{p-1}(\ti M^{2n}, \ti g'(t))}\right..
\eeq
Analogous to Theorem \ref{bi}, similar to the smooth case, we have the following theorem:
\bthm \label{bi3}
Suppose that $(M^{2n}, g)$ is a closed Riemannian manifold with infinite fundamental group $\G=\pi_1(M^{2n})$ of $M^{2n}$.
We can construct a family of $\G$-invariant (measurable) Lipschitz almost K\"ahler structures $(\ti\w', \ti J'(t), \ti g'(t))$, $t\in[0,1]$, on $\ti M^{2n}$ which are with equivalent norms satisfying Lipschitz condition \eqref{ea11}. Where $(\ti\w'$, $\ti J'(0)$, $\ti g'(0))$ is a $\G$-invariant (measurable) Lipschitz K\"ahler flat structure on $\ti M^{2n}$, $(\ti\w'$, $\ti J'(0)$, $\ti g'(0))$ is a $\G$-invariant K\"ahler flat structure on $\overset{\circ}{\ti M} {}^{2n}$, and $$\ti g'(1)=\ti g, \frac{\ti\w'^n}{n!}=dvol_{\ti g'(t)}=dvol_{\ti g'(1)}=dvol_{\ti g},t\in [0,1),$$ are smooth and $\G$-invariant on $\ti M^{2n}$.
$$L^2_k\Omega^p(\ti M^{2n},\ti g'(t)), t\in[0,1),\quad \text{and}\quad L^2_k\Omega^p(\ti M^{2n},\ti g), 0\leq p\leq 2n$$ are quasi-isometric, where $k$ is a nonnegative integer.
For the $\G$-invariant (measurable) Lipschitz almost K\"ahler manifold $(\ti M^{2n}$, $\ti\w'$, $\ti J'(t)$, $\ti g'(t))$ we have Hodge-Kodaira decomposition for $0\leq p\leq 2n$, $t\in[0, 1]$
\begin{align*}
L^2\Omega^p(\ti M^{2n}, \ti g'(t))=\mathcal{H}^p_{(2)}(\ti M^{2n}, \ti g'(t))&\oplus\overline{dL^2\Omega^{p-1}(\ti M^{2n}, \ti g'(t))}\\
&\oplus\overline{d^{\ti *'_t}L^2\Omega^{p+1}(\ti M^{2n}, \ti g'(t))},
\end{align*}
$$b^p_{(2)}(\ti M^{2n}, \ti g'(t))\coloneqq\dim_{\G}\mathcal{H}^p_{(2)}(\ti M^{2n}, \ti g'(t))=\dim_{\G}\mathcal{H}^p_{(2)}(\ti M^{2n}, \ti g)<\infty,$$ and $\mathrm{H}^p_{(2)}(\ti M^{2n}, \ti g'(t))$
and $\mathcal{H}^p_{(2)}(\ti M^{2n}$, $\ti g'(t))$ are isomorphic, $$\mathcal{H}^p_{(2)}(\ti M^{2n}, \ti g'(t))\subset \W^p(\ti M^{2n})\cap L^2\W^p(\ti M^{2n},\ti g).$$ In particular,
$$\chi_{(2)}(\ti M^{2n},\ti g'(t))=\chi_{(2)}(\ti M^{2n},\ti g'(1))=\chi_{(2)}(\ti M^{2n},\ti g)=\chi(M^{2n}, g),$$ where $$\chi_{(2)}(\ti M^{2n},\ti g'(t))\coloneqq\sum_{p=0}^{2n}(-1)^pb_{(2)}^p(\ti M^{2n},\ti g'(t)).$$
\ethm
Suppose that $(M^{2n},g)$ is a $2n$-dimensional closed Riemannian manifold with strictly negative (resp. nonpositive) sectional curvature. As done in K\"ahler case, we have the following key theorem:
\bthm\label{bi4}
Suppose that $(M^{2n}, g)$ is a $2n$-dimensional closed Riemannian manifold which has strictly negative (resp. nonpositive) sectional curvature. Then the universal cover $(\ti M^{2n}, \ti g)$ of $(M^{2n}, g)$ is a complete
$\G=\pi_1(M)$-invariant Riemannian manifold of bounded geometry. By a disjoint partition for $M^{2n}$, we construct a $\G$-invariant open dense submanifold $\overset{\circ}{\ti M} {}^{2n}$ of $M^{2n}$. Hence we can construct a family of $\G$-invariant (measurable) Lipschitz almost K\"ahler structures $(\ti \w',\ti J'(t),$ $\ti g'(t))$, $t\in[0,1]$, restricted to $\overset{\circ}{\ti M} {}^{2n}$, which are with equivalent norms satisfying Lipschitz condition \eqref{ea11} and $$\text{coefficients of }\ti g(t)\in C^\infty(\overset{\circ}{\ti M} {}^{2n})\cap L^\infty_k(\ti M^{2n}, \ti g)$$ for $k$ nonnegative integer. Where $(\ti \w',\ti J'(0), \ti g'(0))$ is a $\G$-invariant K\"ahler flat structure on $\overset{\circ}{\ti M} {}^{2n}$, and $\ti g'(1)=\ti g$ which is a $\G$-invariant Riemannian metric on
$\ti M^{2n}$. Therefore
$\ti\w'$ is $d$(bounded) (resp. d(sublinear)) on $\overset{\circ}{\ti M} {}^{2n}$ in the sense of Definition 1.1, that is, there is a 1-forms $\ti\b'$ such that $\ti\w'=d\ti\b'$ which is $d$(bounded) (resp. d(sublinear))
in the $L^\infty(\ti M^{2n},\ti g)$ sense.
\ethm
We will prove Theorem \ref{bi4} later. By Theorem \ref{bi4} we now give a proof of Theorem \ref{mc}.
\begin{proof}[Proof of Theorem \ref{mc}]
By Theorem \ref{bi4}, Theorem \ref{mc} can be regarded as a corollary of Theorem \ref{main} since its proof is same to the proof of Theorem \ref{main}.

Suppose that $(M^{2n}, g)$ is a closed $2n$-dimensional Riemannian manifold with nonpositive (resp. strictly negative) sectional curvature. Then
the universal covering, $(\ti M^{2n}, \ti g)$, of $(M^{2n}, g)$ has bounded geometry, $\ti g$ is smooth $\G$-invariant on $\ti M^{2n}$, $\G=\pi_1(M^{2n})$.
By Proposition \ref{bi1}, one can construct a family of $\G$-invariant (measurable) Lipschitz almost K\"ahler structures $(\ti\w', \ti J'(t), \ti g'(t))$, $t\in [0, 1]$, on $\ti M^{2n}$ which are with equivalent norms (that is, with Lipschitz condition \eqref{ea11}). Restricted to $\overset{\circ}{\ti M} {}^{2n}$, $\ti g'(t))$, $t\in [0, 1]$, are quasi isometry. Here $\ti g'(1)=\ti g$, for $t\in [0,1]$,
$$dvol_{\ti g}=dvol_{\ti g'(t)}=\frac{\ti\w'^n}{n!}$$ and $\ti\w'(\cdot, \ti J'(1)\cdot)=\ti g'(1)=\ti g$ is smooth $\G$-invariant on $\overset{\circ}{\ti M} {}^{2n}$ which can be extended to $\ti M^{2n}$, in particular, $(\ti M^{2n}, \ti g'(1))=(\ti M^{2n}, \ti g)$ is a smooth complete $\G$-invariant Riemannian manifold with nonpositive (resp. strictly negative) sectional curvature which has bounded geometry, and $(\ti\w', \ti J'(0), \ti g'(0))$ is a $\G$-invariant (resp. (measurable) Lipschitz) K\"ahler flat structure on $\overset{\circ}{\ti M} {}^{2n}$ (resp. $\ti M^{2n}$). Note that in general it does not exist any $\G$-invariant almost complex structure on $\ti M^{2n}$, but $(\ti\w', \ti J'(t), \ti g'(t))$, $t\in[0, 1]$, is a family of $\G$-invariant almost K\"ahler structures on $\overset{\circ}{\ti M} {}^{2n}$ which are quasi isometric (cf. \eqref{ea11}).
As done in Section 3, $L^2_k\Omega^p(\ti M^{2n}, \ti g'(t))$, $t\in[0, 1]$, which are well-defined and quasi isometric (that is, are with equivalent norms due to \eqref{cy0}), can be approximated by $\Omega^p_c(\ti M^{2n})$. Since $(\ti\w', \ti J'(t), \ti g'(t))$, $t\in[0, 1]$, is a family of $\G$-invariant (measurable) Lipschitz almost K\"ahler structures
on $\ti M^{2n}$, then one can define Hilbert $\G$-cochain complex on $\ti M^{2n}$ with respect to the $\G$-invariant (measurable) Lipschitz metric $\ti g'(t)$, $t\in[0, 1]$ and obtain $L^2$ Hodge-Kodaira decomposition.
By Theorem \ref{bi3}, for $0\leq p\leq 2n$, $t\in[0, 1]$,
\begin{align}\label{}
L^2\Omega^p(\ti M^{2n}, \ti g'(t))=\mathcal{H}^p_{(2)}(\ti M^{2n}, \ti g'(t))&\oplus\overline{dL^2\Omega^{p-1}(\ti M^{2n}, \ti g'(t))}\nonumber\\
&\oplus\overline{d^{\ti *'_t}L^2\Omega^{p+1}(\ti M^{2n}, \ti g'(t))},
\end{align}
and $\mathcal{H}^p_{(2)}(\ti M^{2n}, \ti g'(t))$ is smooth on $\ti M^{2n}$. Moreover,
\begin{align}
b^p_{(2)}(\ti M^{2n}, \ti g)&=\dim_{\G}\mathcal{H}^p_{(2)}(\ti M^{2n}, \ti g)=b^p_{(2)}(\ti M^{2n}, \ti g'(t))\nonumber\\
&=\dim_{\G}\mathcal{H}^p_{(2)}(\ti M^{2n}, \ti g'(t))<\infty,
\end{align}
and $\mathrm{H}^p_{(2)}(\ti M^{2n}, \ti g'(t))$ and $\mathcal{H}^p_{(2)}(\ti M^{2n}$, $\ti g'(t))$ are isomorphic for $0\leq p\leq 2n$, $t\in[0, 1]$. Hence, we define Lefschetz map
$$L^k_{\ti\w'}:\Omega^{n-k}(\overset{\circ}{\ti M} {}^{2n})\to\Omega^{n+k}(\overset{\circ}{\ti M} {}^{2n}),\quad 0\leq k\leq n.$$
Since $\ti M {}^{2n}\backslash\overset{\circ}{\ti M} {}^{2n}$ has Lebesgue measure zero, Lefschetz map can be extended to $L^2_m\Omega^{n-k}(\ti M {}^{2n})$. As done in symplectic case (see Section 4),
by Theorem \ref{bi3} and Atiyah $\G$-index theorem, we have
\beq \label{eg}
\chi(M^{2n}, g)=\chi_{(2)}(\ti M^{2n}, \ti g)=\chi_{(2)}(\ti M^{2n}, \ti g'(1))=\chi_{(2)}(\ti M^{2n}, \ti g'(0)).
\eeq
Note that
$$ \text{coefficients of}\ \ti g'(0),\ti *'_0\ \text{are in}\ C^\infty(\overset{\circ}{\ti M} {}^{2n})\cap L^\infty_k(\ti M^{2n}, \ti g)$$
 for $k$ nonnegative integer, it is easy to see that $[L_{\ti\w'}, \Delta'_0]=0$ on $\overset{\circ}{\ti M} {}^{2n}$, where
$$\Delta'_0\coloneqq dd^{\ti*'_0}+d^{\ti*'_0}d,\quad \text{and}\quad d^{\ti*'_0}=-*_{\ti g'(0)}d*_{\ti g'(0)}.$$
Since $\ti g'(0)$ is a (measurable) Lipschitz K\"ahler flat metric on $\ti M^{2n}$ (resp. K\"ahler flat metric on $\overset{\circ}{\ti M} {}^{2n}$) which is regarded as a singular K\"ahler metric on $\ti M^{2n}$, similar to Theorem \ref{bj}, for $0\leq k\leq n$, we have
\beq\label{eg1}
L^k_{\ti\w'}:\mathcal{H}^{n-k}_{(2)}(\overset{\circ}{\ti M} {}^{2n}, \ti g'(0))\to\mathcal{H}^{n+k}_{(2)}(\overset{\circ}{\ti M} {}^{2n}, \ti g'(0))
\eeq
are isomorphisms. Hence, $\mathcal{H}^{n-k}_{(2)}(\ti M^{2n}, \ti g'(0))\cong\mathcal{H}^{n+k}_{(2)}(\ti M^{2n}, \ti g'(0))$.

If  $(M^{2n}, g)$ has nonpositive sectional curvature, then it is similar to Theorem \ref{bk}, by Theorem \ref{bi4} and \eqref{eg1}, we have
$$b^p_{(2)}(\ti M^{2n}, \ti g)=b^p_{(2)}(\ti M^{2n}, \ti g'(0))=0 \quad \text{for}\quad p\not= n.$$ Hence by \eqref{eg} we have
\beq \label{eh}
(-1)^n\chi(M^{2n}, g)=(-1)^n\chi_{(2)}(\ti M^{2n}, \ti g)=(-1)^n\chi_{(2)}(\ti M^{2n}, \ti g'(0))\geq 0.
\eeq

If  $(M^{2n}, g)$ has strictly negative sectional curvature, as done in symplectic hyperbolic case, we can construct the signature operator in the sense of distributions,
$$D'^+=d+d^{*_0}:E^+\to E^-,$$ where $$E^{\pm}=\Gamma(\Lambda^\pm T^*\ti M^{2n}\otimes\C).$$ $\ti\w'=d\ti \b'$ is the Chern class of a trivial complex line bundle $\ti L$ over $\overset{\circ}{\ti M} {}^{2n}$.
Let $$\nabla'^{(k)}=d+\frac{1}{k}\sqrt{-1}\ti \b', k\in\N,$$ where $\ti \b'$ is bounded in the $L^\infty(\ti M^{2n}, \ti g)$ sense. Thus $$D'^+\otimes\nabla'^{(k)}:E^+\otimes\ti L\to E^-\otimes\ti L$$ is $\frac{1}{k}$-small perturbation of $D'^+$. Note that $$\frac{\ti\w'^n}{n!}=dvol_{\ti g}=\pi^*dvol_g$$ is smoothly extended
to $\ti M^{2n}$ and $\G$-invariant. Hence, $$\int_{F}(\ti \w')^n=\int_{\pi(F)}(\pi_*\ti \w')^n=n!vol_g(M^{2n})\neq0.$$
It is the same as Theorem \ref{bp}, by Theorems \ref{bi3}, \ref{bi4} and \eqref{eg1}-\eqref{eh}, we have
$$\mathcal{H}^n_{(2)}(\ti M^{2n}, \ti g)\neq0\quad \text{and} \quad\mathcal{H}^p_{(2)}(\ti M^{2n}, \ti g)=0 \quad \text{for} \quad p\neq n.$$ Hence,
\beq \label{ei}
(-1)^n\chi(M^{2n}, g)=(-1)^n\chi_{(2)}(\ti M^{2n}, \ti g)=(-1)^n\chi_{(2)}(\ti M^{2n}, \ti g'(0))> 0.
\eeq
Therefor, by \eqref{eh} and \eqref{ei} we complete the proof of Theorem \ref{mc}.
\end{proof}
\brem
M. T. Anderson \cite{Ad1,Ad2} had constructed some simply connected complete Riemannian manifold with negative sectional curvature such that the space of harmonic forms $\mathcal{H}^p_{(2)}(M)$ is nontrivial for some $p$ with
$2p\neq \dim(M)$. This manifold $M$ does not admit a cocompact free proper action of a discrete group $\G$ by isometries. The example shows the necessity of the vanishing condition of cocompactness for the vanishing theorem to be true. Thus the example can not to be used to construct counterexample to the Chern-Hopf conjecture in the nonpositive curvature case.
\erem

Finally, for sake of completeness, we give a proof of Theorem \ref{bi4}.
\begin{proof}[Proof of Theorem \ref{bi4}]
We will derive the conclusion by using the pull-back metric, as J. Cao and F. Xavier did in \cite[Theorem 1]{CX} (nonpositivly curved case), X. Cheng \cite[Theorem 3.1]{Chg} (strictly negative case).

Fix $p\in \overset{\circ}{\ti M}{}^{2n}$ and denote by $\exp_p:T_p\ti M^{2n}\to \ti M^{2n}$ the exponential map based at $p$. As done in the K\"ahler case, we have the following lemma:
\blem{\rm(}cf. Cao-Xavier \cite[Lemma 1]{CX}{\rm)}\label{lem8}
If  $(\ti M^{2n}, \ti g)$ has nonpositive sectional curvature $K$, then $\exp_p:T_p\ti M^{2n}\to \ti M^{2n}$ is a diffeomorphism. Consider the map $\tau^1_t:\ti M^{2n}\to \ti M^{2n}$, given by $x\to \exp_p(t\exp_p^{-1}(x))$,
where $0\leq t\leq 1$. Then
\beq\label{ea2}
|(\tau^1_t)_*\xi|\leq t|\xi|
\eeq
for every tangent vector $\xi$.
\elem
Also we have the following elementary result:
\blem{\rm(}cf. Cao-Xavier \cite[Lemma 2]{CX}{\rm)}\label{lem9}
Let $\Psi$ be a closed 2-form in $\R^{2n}$. Then the 1-form $\Phi$ defined by
$$ \Phi(x)=r\int_0^1[(\tau^1_t)^*(\iota_{\frac{\partial}{\partial r}}\Psi)](x)dt$$
satisfies $d\Phi=\Psi$, where $$\frac{\partial}{\partial r}=\sum_{i=1}^{2n}\frac{x_i}{r}\frac{\partial}{\partial x_i},r=(\sum_{i=1}^{2n}x_i^2)^{\frac{1}{2}}\quad and \quad\tau^1_t(x)=tx.$$
\elem
Now for the nonpositive sectional curvature case, we give a proof of Theorem \ref{bi4}, that is, construct a 1-form $\ti\b'$ on $\overset{\circ}{\ti M} {}^{2n}$ such that $\ti\w'=d\ti\b'$ and $\ti\b'$ is d(sublinear) in the $L^\infty(\ti M^{2n},\ti g)$ sense (cf. the proof of Theorem 1 in J. Cao and F. Xavier \cite{CX}).

Let $(x_1,\cdots, x_{2n})$ be coordinates in $\R^{2n}\cong T_p\ti M^{2n}$ and consider the pull-back $h$ of the metric $\ti g$ under $\exp_p:T_p\ti M^{2n}\to \ti M^{2n}$. By Lemma \ref{lem8}, $\R^{2n}\cong T_p\ti M^{2n}$ is diffeomorphic to $\ti M^{2n}$. There is a way to interpret the map $\tau^1_t$.
This comes from Lemma \ref{lem8} with $(\ti M^{2n}, \ti g)$ being replaced by $(\R^{2n}\cong T_p\ti M^{2n}, h)$.

Let $\Psi=(\exp_p)^*\ti\w'$ on $\overset{\circ}{\R} {}^{2n}$, where $\overset{\circ}{\R} {}^{2n}=(\exp_p)^*\overset{\circ}{\ti M} {}^{2n}$. Hence $\overset{\circ}{\R} {}^{2n}$ is an open dense submanifold of $\R^{2n}$,
$\R^{2n}\setminus\overset{\circ}{\R} {}^{2n}$ has Lebesgue measure zero. $\Psi$ is a closed 2-form on $\overset{\circ}{\R} {}^{2n}$, it is easy to see that $||\Psi||_{L^\infty(\R^{2n})}\leq C$ being finite. Let $\Phi$
be given by Lemma \ref{lem9}. By Lemma \ref{lem8}, we have for $x\in\overset{\circ}{\R} {}^{2n}$,
\beq\label{ea3}
|(\tau^1_t)^*\Phi(x)|_h\leq t|\Phi(\tau^1_t(x)|_h,
\eeq
where $C$ is a positive constant depending only on $h$. From \eqref{ea3} and Lemma \ref{lem9}, it follows that $$|\Phi(x)|_h\leq C_1r\sup_{0\leq t\leq 1}|\Psi(tx)|_h,$$ where $$r=(\sum_{i=1}^{2n}x_i^2)^{\frac{1}{2}},$$ and $C_1$ is
a positive constant depending only on $h$.

In particular,
$$
|\Phi(x)|_h\leq C_1\rho_h(o,x)\sup_{\overset{\circ}{\R} {}^{2n}}|\Psi|_h\leq C_2\rho_h(o,x).
$$
Hence $\Psi$ is d(sublinear). This means that $\ti\w'$ is d(sublinear) in the $L^\infty(\ti M^{2n}, \ti g)$ sense, that is, $$\ti\w'=d\ti\b'(g),||\ti\b'(g)||_{L^\infty(\ti M^{2n}, \ti g)}\leq C\rho_{\ti g}(x_0,x),$$
where $C$ is a positive constant depending only on $\ti g$ (in fact, on $g$) and point $x_0\in \ti M^{2n}$. This completes the proof of Theorem \ref{bi4} for nonpositive sectional curvature case.
\vskip 6pt
To prove strictly negative sectional curvature case, let $p$ be a fixed point on $\overset{\circ}{\ti M}{}^{2n}$. Let $\gamma(t,x)$ denote the unit speed geodesic satisfying $p=\gamma(-t_0(x))$, $t_0(x)>0$ and $x=\gamma(0)\in \ti M^{2n}$. Let $Y(t)$ denote a perpendicular Jacobi vector field
along with $Y(0)\neq0$ and $Y(-s)$. We will estimate the growth of $Y(t)$ as follows:
When $(\ti M^{2n}, \ti g)$ has sectional curvature bounded above by a negative constant $K\leq-a<0$, we can obtain from Rauch comparison theorem (cf. I. Chavel \cite{Cha}), that every perpendicular Jacobi field along $\gamma(t,x)$
with $Y(0)\neq0$ and $Y(-t_0(x))$ has
\beq\label{ea4}
||Y(t)||_{\ti g}\geq C||Y(0)||_{\ti g}e^{at}\geq 0,
\eeq
where $C$ is a positive constant depending only on $ M^{2n}$ and $g$.

Let $\tau^2_t:\ti M^{2n}\to\ti M^{2n}$ defined by $\tau^2_t(x)=\gamma(t,x)$, $t\in (-\infty, +\infty)$. That is, $\tau^2_t(x)$ is the point at $\gamma(t,x)$ with distance $t$ (with sign) starting from $x$ in the oppositive
direction of $p$. By \eqref{ea4}, we have the following proposition:
\bprop\label{bi5}{\rm(}cf. X. Cheng \cite[Proposition 3.2]{Chg}{\rm)}
Let $(\ti M^{2n}, \ti g)$ be a simply connected complete Riemannian manifold of strictly negative sectional curvature $K\leq-a<0$, then there exists constant $C>0$ depending only on $M^{2n}$ and $g$ such that for all $t>0$,
\beq\label{ea5}
||d\tau^2_t(X)||_{\ti g}\geq C e^{at}||X||_{\ti g},
\eeq
where $X\in T_x\ti M^{2n}$, $X\bot \dot{\gamma}(0, x)$.
\eprop
\begin{proof}
Through $x$, we have a geodesic sphere $S_x(p)\subset \ti M^{2n}$ with center $p$. Choose a smooth curve $$\lambda(s):(-\varepsilon, \varepsilon)\to S_x(p)$$ on $S_x(p)$ satisfying $\lambda(0)=x$, $\lambda'(0)=X\in T_xS_x(p)$.
Join $p$ to each $\lambda(s)$, $s\in(-\varepsilon, \varepsilon)$, with normalized geodesics, denoted by $f_s(t)$, $-t_0(x)\leq t\leq 0$ and extend each geodesic to infinity. Thus, we obtain a variation $f(t,s)$ along
$$\gamma(t, x), \quad t\in[-t_0(x), +\infty),\quad s\in(-\varepsilon, \varepsilon).$$ In fact, $$f(t, s)=\exp_p\frac{t+t_0(x)}{t_0}\exp^{-1}_p\lambda(s).$$

Note that every $f_s(t)$ is a normalized geodesic and $f_s(t)$ is a curve on $S_{\tau^2_t(x)}(p)$, where $S_{\tau^2_t(x)}(p)$ denotes the geodesic sphere with center $p$, passing through $\tau^2_s(x)$. Let
$Y(t)=\frac{\partial f}{\partial s}(t, 0)$ is a perpendicular Jacobi field satisfying $Y(-t_0(x))=0$, $Y(0)=X$. By \eqref{ea4}, we have
$$||Y(t)||_{\ti g}\geq C e^{at}||Y(0)||_{\ti g},\quad t\geq 0.$$
Note that  $f(t, s)=\tau^2_t(\lambda(s))$, where $t\in[-t_0(x), +\infty)$ and $s\in(-\varepsilon, \varepsilon)$. Thus
\begin{align*}
Y(t)&=\frac{\partial f}{\partial s}(t, 0)=d\tau^2_t(\frac{d\lambda}{ds}(0))\\
&=d\tau^2_t(Y(0))=d\tau^2_t(X).
\end{align*}
Hence, $$||d\tau^2_t(X)||_{\ti g}\geq C e^{at}||X||_{\ti g}, \quad t>0.$$
\end{proof}
As done in X. Cheng \cite{Chg}, by Proposition \ref{bi5}, we have the following two propositions:

\bprop\label{}
Let $(\ti M^{2n}, \ti g)$ be as in Proposition \ref{bi5}, then there exists a positive constant $C>0$ depending only on $(M^{2n},g)$ such that for any 2-form $\a$ and $t\geq0$,
\beq\label{ea6}
||(\tau^2_t)^*(\iota_v\a)||_{\ti g(x)}\geq C e^{at}||\iota_v\a||_{\ti g(\tau^2_t(x))};
\eeq
if $0\leq t\leq\rho_{\ti g}(p,x)$, then
\beq\label{ea7}
||(\tau^2_t)^*(\iota_v\a)||_{\ti g(\tau^2_{-t}(x))}\geq C e^{at}||\iota_v\a||_{\ti g(x)},
\eeq
where $\iota_v$ denotes the interior product and $v$ denotes the unit tangent vector of the geodesic $\gamma$ which is the same as Proposition \ref{bi5}.
\eprop
\begin{proof}
This proposition is the same as Proposition 3.3 in X. Cheng \cite{Chg}, hence the proof is ommited.
\end{proof}
\bprop\label{bi6}
Let $(\ti M^{2n}, \ti g)$ be as in Proposition \ref{bi5}, then for any 2-form $\a$ and $0\leq t\leq \rho_{\ti g}(p,x)$,
\beq\label{ea8}
||(\tau^2_{-t})^*(\iota_v\a)||_{\ti g(x)}\leq C e^{-at}||\a||_{\ti g(\tau^2_{-t}(x))},
\eeq
where $\iota_v$ denotes the interior product and $v$ denotes the unit tangent vector of the geodesic $\gamma$ which is the same as Proposition \ref{bi5}.
\eprop
\begin{proof}
This proposition is the same as Proposition 3.4 in X. Cheng \cite{Chg}, hence the proof is ommited.
\end{proof}
Now we can prove the second part of Theorem \ref{bi4}.

Suppose that $(M^{2n}, g)$ is a $2n$-dimensional closed Riemannian manifold of strictly negative sectional curvature $K\leq-a<0$. On the universal covering space $(\ti M^{2n}, \ti g)$ of $(M^{2n}, g)$ we constructed a $\G$-invariant measurable Lipschitz almost K\"ahler structure $(\ti \w', \ti J'(1), \ti g'(1))$, where $\ti g'(1)=\ti g$ is $\G$-invariant smooth Riemannian metric on $\ti M^{2n}$. We consider $(\R^{2n}\cong T_p\ti M^{2n}, h)$ where $p\in\overset{\circ}{\ti M}{}^{2n}$, $h=\exp^*_p\ti g$. Since $\exp_p$ is an isometry, $(\R^{2n}, h)$ is a simply
connected complete Riemannian manifold with strictly negative sectional curvature $K\leq-a<0$. Thus we can apply Proposition \ref{bi6} to $(\R^{2n}, h)$.

Let $\overset{\circ}{\R} {}^{2n}=\exp_p^{-1}(\overset{\circ}{\ti M} {}^{2n})$, then $\overset{\circ}{\R} {}^{2n}$ is an open dense submanifold of $\R^{2n}$, where $\R^{2n}\setminus\overset{\circ}{\R} {}^{2n}$ has Lebesgue measure zero.
Let $\Psi=\exp^*_p\ti\w'$ which is a $\G$-invariant symplectic 2-form on $\overset{\circ}{\R} {}^{2n}$ with $||\Psi||_h<C$, where $C$ is a positive constant depending only on $g$.

Set $$ \Phi(x)=r\int_0^1[(\bar{\tau}^2_s)^*(\iota_{\frac{\partial}{\partial r}}\Psi)](x)ds,$$
where $\bar{\tau}^2_s(x)=sx$, $\frac{\partial}{\partial r}$ and $r$ are given in Lemma \ref{lem9}. We have $d\Phi=\Psi$ on $\overset{\circ}{\R} {}^{2n}$. We will verify $\Phi$ is a bounded 1-form on $\overset{\circ}{\R} {}^{2n}$.
Since $(\ti M^{2n}, \ti g)$ and $(\R^{2n}, h)$ are isometric, the geodesic on $(\R^{2n}, h)$ passing through the origin point $o$ and a point $x\in\overset{\circ}{\R} {}^{2n}$ is just $\bar{\gamma}(s)=sx$, $s\in(-\infty, +\infty)$. Then $\tau^2_{-t}(x)$ ($0\leq t\leq \rho_h(o,x)=||x||_h$) denotes the point between $o$ and $x$ at $\bar{\gamma}(s)$ satisfying the distance with respect to $h$ between $\tau^2_{-t}(x)$ and $x$ is $t$. Exactly, for
$0\leq t\leq \rho_h(p,x)=||x||_h$,
$$
\tau^2_{-t}(x)=\frac{||x||_h-t}{||x||_h}x=sx=\bar{\tau}^2_s(x),
$$
where $$s=\frac{||x||_h-t}{||x||_h},0\leq s\leq 1.$$ Obviously, $$ds=-\frac{dt}{||x||_h},d\tau^2_{-t}=d\bar{\tau}^2_s, (\tau^2_{-t})^*=(\bar{\tau}^2_s)^*\quad and \quad(\frac{\partial}{\partial r})_{sx}=\frac{1}{r}\dot{\bar{\gamma}}(s).$$
By Gauss Lemma (cf. \cite{Cha}),
$$\li\dot{\bar{\gamma}}(s),\dot{\bar{\gamma}}(s)\ri_h=\li x,x\ri_h .$$
Hence, $$(\frac{\partial}{\partial r})_{sx}=\frac{||x||_h}{r}v,$$ where $$v(sx)=\frac{\dot{\bar{\gamma}}(s)}{||x||_h}$$ denotes the unit tangent vector along the geodesic $\gamma$ which is the normalization of $\bar{\gamma}$.
Thus
\begin{align*}
\Phi(x)=&r\int_0^1[(\bar{\tau}^2_s)^*(\iota_{\frac{\partial}{\partial r}}\Psi)](x)ds\\
=&\int_0^{||x||_h}(\tau^2_{-t})^*(\iota_v\Psi)(x)dt.
\end{align*}
By Proposition \ref{bi6},
\beq\label{ea9}
||(\tau^2_{-t})^*(\iota_v\Psi)(x)||_h\leq C e^{-at}||\Psi||_{L^\infty(\R^{2n}, h)}.
\eeq
Hence,
\begin{align}
||\Phi(x)||_h&\leq\int_0^{||x||_h}||(\tau^2_{-t})^*(\iota_v\Psi)||_h(x)dt\nonumber\\
&\leq C\int_0^{||x||_h}e^{-at}dt||\Psi||_{L^\infty(\R^{2n}, h)}\nonumber\\
&\leq\frac{C}{a}||\Psi||_{L^\infty(\R^{2n}, h)}.
\end{align}
Thus, $||\Phi||_{L^\infty(\R^{2n}, h)}<\infty$. Let $\ti\b'=(\exp^{-1}_p)^*\Phi$. Then $\ti\b'$ is a bounded 1-form on $\ti M^{2n}$ in the $L^\infty(\ti M^{2n}, \ti g)$ sense, and $\ti \w'=d\ti\b'$. This completes the proof of Theorem \ref{bi4}.
\end{proof}

Finally, we would like to remark that R. Charney and M. W. Davis posed a topological version of Chern-Hopf conjecture \cite[Conjecture A]{CD}:
\bconj
If $M^{2n}$ is a nonpositively curved, piecewise Euclidean, closed manifold, then
$$(-1)^n\chi(M^{2n})\geq 0.$$
\econj
\vspace{3mm}\par\noindent {\bf{Acknowledgements.}} The authors would like to thank Bo Dai, Youlin Li, Yuxiang Li, Teng Huang, Xiaowei Xu, Zuoqing Wang, Xingwang Xu and Ken Wang for useful conversations.


\end{document}